\documentclass{amsart}
\usepackage{amsmath}
\usepackage{amsfonts,amssymb,amsthm}
\usepackage{mathscinet}
\makeatletter
\newtheorem{theorem}{Theorem}[section]
\newtheorem{lemma}[theorem]{Lemma}
\newtheorem{proposition}[theorem]{Proposition}    
\newtheorem{corollary}[theorem]{Corollary}
\newtheorem*{GST}{The Glimm-Sakai theorem}

\theoremstyle{definition}
\newtheorem{definition}[theorem]{Definition}
\newtheorem{remark}[theorem]{Remark}
\newtheorem{example}[theorem]{Example}
\newtheorem{problem}[theorem]{Problem}

\newtheorem*{NProb}{The Naimark problem}

\newcommand{\A}{\mathfrak{A}}
\newcommand{\B}{\mathfrak{B}}
\newcommand{\g}{\mathfrak{G}}
\newcommand{\p}{\mathfrak{P}}
\newcommand{\s}{\mathfrak{S}}
\newcommand{\gs}{\mathfrak{GS}}
\newcommand{\ps}{\mathfrak{PS}}
\newcommand{\uA}{\widetilde{\mathfrak{A}}}
\newcommand{\HH}{\mathcal{H}}
\newcommand{\KK}{\mathcal{K}}
\newcommand{\core}{\mathop{\mathrm{core}}\nolimits}
\newcommand{\ext}{\mathop{\mathrm{ext}}\nolimits}
\newcommand{\id}{\mathop{\textrm{id}}}
\newcommand{\FCl}{\mathop{\textrm{\bf Cl}_\textrm{\rm F}}}
\newcommand{\FTop}{\mathop{\textrm{\bf Top}_\textrm{\rm F}}}
\newcommand{\Top}{\mathop{\textrm{\bf Top}}}
\newcommand{\Gss}{\mathop{\textrm{\bf Gss}}}

\newcommand{\CGss}{\mathop{C^*\textrm{-}\textrm{\bf Gss}}}

\newcommand{\Sp}{\mathop{\mathrm{sp}}\nolimits}
\newcommand{\re}{\mathop{\mathrm{Re}}\nolimits}
\newcommand{\card}{\mathop{\mathrm{card}}\nolimits}
\newcommand{\Prim}{\mathop{\mathrm{Prim}}\nolimits}
\newcommand{\OG}{\mathop{\hat{\Gamma}}}
\newcommand{\RG}{\mathop{\hat{\Gamma}_0}}
\title[(Primitive) spectra of $C^*$-algebras and the Naimark problem]{Banach space theoretical construction of (primitive) spectra of $C^*$-algebras and the Naimark problem revisited}
\author{Ryotaro Tanaka}
\address{Katsushika Division, Institute of Arts and Sciences, Tokyo University of Science, Tokyo 125-8585, Japan}
\email{r-tanaka@rs.tus.ac.jp}
\thanks{This work was supported by JSPS KAKENHI Grant Number JP19K14561.}


\begin{document}
\begin{abstract}
The Naimark problem asks whether $C^*$-algebras with singleton spectra are necessarily elementary. The separable case was solved affirmatively in 1953 by Rosenberg. In 2004, Akemann and Weaver gave a counterexample to the Naimark problem for non-separable $C^*$-algebras in the setting of ZFC $+~\diamondsuit_{\aleph_1}$, where $\diamondsuit_{\aleph_1}$ is Jensen's diamond principle. From this, at least, the affirmative answer to the Naimark problem can no longer be expected although a counterexample is not constructed in ZFC alone yet. In this paper, we study the difference between elementary $C^*$-algebras and those with singleton spectra, and find a property $P$ written in the language of closure operators such that a $C^*$-algebra is elementary if and only if it has the singleton spectrum and the property $P$. Banach space theoretical construction of (primitive) spectra of $C^*$-algebras plays important roles in the theory. Characterizations of type I or CCR or (sub)homogeneous $C^*$-algebras are also given. These results are applied to a geometric nonlinear classification problem for $C^*$-algebras.
\end{abstract}
\maketitle

\section{Introduction}

The Naimark problem was a long-standing problem in the representation theory of $C^*$-algebras formulated in terms of spectra, where the spectrum $\hat{\A}$ of a $C^*$-algebra $\A$ is defined as the set of unitary equivalence classes of irreducible representations of $\A$.
\begin{NProb}

Let $\A$ be a $C^*$-algebra. Are the following statements equivalent?
\begin{itemize}
\item[{\rm (i)}] $\A$ is \emph{elementary}, that is, there exists a complex Hilbert space $\HH$ such that $\A$ is $*$-isomorphic to the algebra $K(\HH)$ of all compact operators on $\HH$.
\item[{\rm (ii)}] $\hat{\A}$ is a singleton.
\end{itemize}
\end{NProb}
The implication (i) $\Rightarrow$ (ii) was observed by Naimark~\cite{Nai48}, and the Naimark problem was asked in \cite{Nai51} based on his work. In 1953, Rosenberg~\cite{Ros53} quickly gave an affirmative answer to the separable (in fact, separably representable) case. However, the non-separable case had remained open for about fifty years. A great progress on the Naimark problem was relatively recently, in 2004, achieved by Akemann and Weaver~\cite{AW04}. They provide a counterexample to the non-separable case with the aid of Jensen's diamond principle $\diamondsuit_{\aleph_1}$, and noted that the existence of a counterexample to the Naimark problem which is generated by $\aleph_1$ elements is independent of ZFC (Zermelo-Fraenkel set theory with the axiom of choice); see \cite[Theorem 5 and Corollary 7]{AW04}. In this direction, very recently, Calder\'{o}n and Farah~\cite{CF23} also obtained a counterexample with the aid of $\diamondsuit^{\rm Cohen}~+$ CH, where $\diamondsuit^{\rm Cohen}$ is a weakening of $\diamondsuit_{\aleph_1}$ and CH is the continuum hypothesis. Moreover, $\diamondsuit^{\rm Cohen}~+$ CH is relatively consistent with the negation of $\diamondsuit_{\aleph_1}$; see~\cite[Appendix B]{CF23}. The both ideas of constructing counterexamples to the Naimark problem were based on the Kishimoto-Ozawa-Sakai theorem on the homogeneity of the pure state space of a separable simple $C^*$-algebra~\cite{KOS03}; see also the earlier work by Futamura, Kataoka and Kishimoto~\cite{FKK01}. Since $\diamondsuit_{\aleph_1}$ is relatively consistent with ZFC, the implication (ii) $\Rightarrow$ (i) cannot be proved in ZFC. However, this does not directly mean that the mentioned implication is false. It is still not known if there exists a counterexample to the Naimark problem in ZFC alone; see~\cite[Concluding remarks]{CF23}. It should be mentioned that partial positive results in the setting of graph $C^*$-algebras were given by Suri and Tomforde~\cite{ST17}. A modified version of the Naimark problem was asked in 2006 by Sakai~\cite{Sak06} and the long history of the Naimark problem toward its modification was summarized in his comprehensive survey~\cite{Sak08}.

The main objective of the present paper is to study the difference between (i) and (ii) in the Naimark problem for non-separable $C^*$-algebras, and to find a missing piece for characterizing non-separable elementary $C^*$-algebras related to spectra. More precisely, we investigate a property $P$ written in the language of closure operators such that (i) is equivalent to (ii) + $P$ as well as $P$ does not imply (ii). Naturally, if $\hat{A}$ is a singleton, it contains very little information as topological space. Hence, we have to consider another natural structure of a $C^*$-algebra including the information of its spectrum. A candidate for this is the concept of geometric structure spaces that was introduced in \cite{Tan22b} as closure spaces modeled after (primitive) spectra of $C^*$ algebras for classifying Banach spaces with respect to their structure determined by Birkhoff-James orthogonality. It was shown in~\cite[Theorem 5.2]{Tan22b} that the geometric structure space of an abelian $C^*$-algebra is homeomorphic to its spectrum. Meanwhile, the geometric structure space of a $C^*$-algebra forms a topological space if and only if the algebra is abelian; see~\cite[Theorem 3.5]{Tan23c}. Thus, geometric structure spaces of $C^*$-algebras are rarely topologizable and our work should be carried out in the class of closure spaces that is much wider than that of topological spaces.

The first milestone of our study is Theorem~\ref{main-theorem} which gives a Banach space theoretical way of constructing (primitive) spectra of $C^*$-algebras. This will be achieved by applying certain transformations (defined for general closure spaces) to geometric structure spaces of $C^*$-algebras. The necessary theories of closure spaces and Banach spaces for this purpose are developed in Sections~\ref{Sect:Cl} and \ref{Sect:BS}, respectively. It turns out in Corollary~\ref{spec-preserved} that (primitive) spectra of $C^*$-algebras are invariant under (closure space) homeomorphisms between geometric structure spaces of them. Hence, the geometric structure space of a $C^*$-algebra remenbers the information of its spectrum in this sense. Another main result of Section~\ref{Sect:CAlg}, Theorem~\ref{irr-rep}, is a useful tool for reducing our argument to each irreducible representation. Based on this, we develop nonlinear theory for irreducible $C^*$-algebras in Section~\ref{Sect:Irr}. The machinery is constructed by normal parts of geometric structure spaces and finite-homeomorphisms between them. In Theorem~\ref{continuity}, we characterize continuity of such mappings induced by pairs of bounded linear (or conjugate-linear) isomorphisms between complex Hilbert spaces on which given irreducible $C^*$-algebras act. Moreover, Theorem~\ref{weak-preserver}, the most technical part of the present paper, describes the possible forms of continuous finite-homeomorphisms between normal parts of geometric structure spaces of irreducible $C^*$-algebras acting on infinite-dimensional complex Hilbert spaces. We revisit the Naimark problem in Section~\ref{Sect:NP}. The main result is Theorem~\ref{elementary}, which characterizes non-separable elementary $C^*$-algebras in terms of geometric structure spaces. Since the geometric structure space and spectrum of a $C^*$-algebra can be recognized as Banach space structures, this gives a Banach space theoretical characterization of elementary $C^*$-algebras among all $C^*$-algebras. Similar characterizations of type I or CCR or (sub)homogeneous $C^*$-algebras are also given in Theorem~\ref{type-I-homogeneous}. The results in Sections~\ref{Sect:Irr} and \ref{Sect:NP} are applied to a nonlinear classification problem for $C^*$-algebras considered in Section~\ref{Sect:Cls} which asks whether $C^*$-algebras are Jordan $*$-isomorphic if and only if they have homeomorphic geometric structure spaces. Partial affirmative answers to this problem are given in Theorems~\ref{cpt-op-alg} and \ref{full-algebra}. The paper concludes with two geometric nonlinear classification problems for $C^*$-algebras related to the above-mentioned one.

\section{Preliminaries}\label{Sect:Pre}

\subsection{Categories}

A \emph{category} \textbf{C} consists of \emph{objects} $\mathrm{ob}(\textbf{C})$, \emph{hom-set} $\textbf{C}(X,Y)$ of objects (whose elements are called \emph{morphisms}), and \emph{composition} of morphisms $\textbf{C}(X,Y) \times \textbf{C}(Y,Z) \ni (f,g) \mapsto g \circ f \in \textbf{C}(X,Z)$. These data must satisfy the following conditions:
\begin{itemize}
\item[(i)] If $f:X \to Y$, $g:Y \to Z$ and $h:Z \to W$ are morphisms, then $h \circ (g \circ f) = (h \circ g) \circ f$.
\item[(ii)] For each object $X$, there exists a morphism $\mathrm{id}_X :X \to X$ such that $\mathrm{id}_X \circ f = f$ and $g \circ \mathrm{id}_X = g$ whenever $f:Y \to X$ and $g:X \to Z$ are morphisms.
\end{itemize}
A morphism $f:X \to Y$ is called an \emph{isomorphism} if there exists a morphism $g:Y \to X$ such that $g \circ f =\mathrm{id}_X$ and $f \circ g = \mathrm{id}_Y$, in which case, $g$ is the \emph{inverse} of $f$ and written as $f^{-1}$. If there exists an isomorphism $f:X \to Y$, then objects $X$ and $Y$ are \emph{isomorphic}.

Let \textbf{C} is a category. A \emph{subcategory} \textbf{D} of \textbf{C} is a category such that $\mathrm{ob}(\textbf{D}) \subset \mathrm{ob}(\textbf{C})$ and $\textbf{D}(X,Y) \subset \textbf{C}(X,Y)$ for each objects $X,Y$ of $\textbf{D}$. A subcategory \textbf{D} of \textbf{C} is said to be \emph{full} if $\textbf{D}(X,Y) = \textbf{C}(X,Y)$ for each objects $X,Y$ of $\textbf{D}$. We note that a full subcategory is uniquely determined by $\mathrm{ob}(\textbf{D})$.

Let \textbf{C} and \textbf{D} be categories. A \emph{functor} $F$ from \textbf{C} to \textbf{D} is mappings $X \mapsto F(X)$ from $\mathrm{ob}(\textbf{C})$ to $\mathrm{ob}(\textbf{D})$ and $f \mapsto F(f)$ from $\textbf{C}(X,Y)$ to $\textbf{D}(F(X),F(Y))$ such that $F(g \circ f) = F(g) \circ F(f)$ and $F(\mathrm{id}_X) = \mathrm{id}_{F(X)}$. The \emph{identity functor} $\mathrm{id}_\textbf{C}$ on \textbf{C} is the functor given by $F(X) = X$ and $F(f) = f$. Let $F$ be a functor from \textbf{C} to \textbf{D}. If there exists a functor $G$ from \textbf{D} to \textbf{C} such that $G \circ F = \mathrm{id}_\textbf{C}$ and $F \circ G = \mathrm{id}_\textbf{D}$, then $F$ is called an \emph{isomorphism of categories}. The categories \textbf{C} and \textbf{D} are \emph{isomorphic} if there exists an isomorphism of categories from \textbf{C} to \textbf{D}. See, for example, \cite{Gra18} for the general theory of categories.

\subsection{Banach spaces}

The scalar field of a Banach space is denoted by $\mathbb{K}$. A Banach space over $\mathbb{R}$ (or over $\mathbb{C}$) is called a \emph{real} (or \emph{complex}) Banach space. For a Banach space $X$, let $B_X = \{ x \in X :\|x\| \leq 1\}$ and $S_X = \{ x \in X :\|x\|=1\}$. The sets $B_X$ and $S_X$ are called the unit ball and unit sphere of $X$, respectively. If $X$ is a complex Banach space, then its complex conjugate is denoted by $\overline{X}$, that is, $\overline{X}$ is the Banach space having the same underlying vector space, addition and norm as $X$ and the conjugate scalar multiplication $(c,x) \mapsto \overline{c}x$. If $\HH$ is a Hilbert space, then so is $\overline{\HH}$. The inner product $\langle \cdot ,\cdot \rangle_a$ on $\overline{\HH}$ is given by $\langle x,y \rangle_a = \langle y,x \rangle$ for each $x,y \in \HH$. The (continuous) dual space of $X$ is denoted by $X^*$. Although the symbols $A^*$ and $\rho^*$ are also used for the adjoint of elements $A$ and $\rho$ of a $C^*$-algebra and its dual, the meanings will be clear from the context. A subset $C$ of a Banach space is said to be \emph{convex} if $(1-t)x+ty \in C$ whenever $x,y \in C$ and $t \in [0,1]$. A subset $D$ of a convex set $C$ is called a \emph{face} of $C$ if $x,y \in C$ and $(1-t)x+ty \in D$ for some $t\in (0,1)$ implies that $x,y \in D$. If further $D \neq C$, then $D$ is said to be \emph{proper}. A \emph{maximal face} of a convex set $C$ is a proper face of $C$ that is not contained in other proper faces of $C$. If a singleton $\{x\}$ is a face of $C$, then it is called an \emph{extreme point} of $C$. The set of all extreme points of $C$ is denoted by $\ext (C)$. An important example of a face of $B_X$ is an \emph{exposed face} which can be represented as $f^{-1}(\{1\}) \cap B_X$ for some $f \in S_{X^*}$. It is well-known that each maximal face of $B_X$ is exposed. Let $X$ and $Y$ are Banach spaces over the same scalar field. Then, a mapping $T:X \to Y$ is called an \emph{isometry} if $\|Tx\|=\|x\|$ for each $x \in X$. If an isometry $T$ maps $X$ onto $Y$, then it is called an \emph{isometric isomorphism}. If there exists an isometric isomorphism $T:X \to Y$, then $X$ is \emph{isometrically isomorphic} to $Y$. For basics of Banach spaces needed in the present paper, see~\cite{Meg98}.

\subsection{$C^*$-algebras}

A $C^*$-algebra is a complex Banach algebra with involution that satisfies the Naimark axiom $\|A^*A\| = \|A\|^2$ for each $A\in \A$. A $C^*$-algebra is said to be \emph{unital} if it contains the multiplicative unit $\mathbf{1}$. For a $C^*$-algebra $\A$, let $\uA =\A$ if $\A$ is unital, and let $\uA$ be its unitization of $\A$ if $\A$ is non-unital, where the unitization of a non-unital $C^*$-algebra $\A$ is the vector space $\A \oplus \mathbb{C}$ endowed with the $C^*$-norm 
\[
\|(A,c)\| = \sup \{ \|AB+cB\| : B \in B_\A\} .
\]
The multiplicative unit of $\uA$ is $(0,1)$ and an element $(A,c)$ of $\uA$ is denoted by $A+c\mathbf{1}$. If $\A$ is a $C^*$-subalgebra of a unital $C^*$-algebra, then it is natural to identify $\uA$ with the $C^*$-algebra generated by $\A$ and $\mathbf{1}$. The spectrum $\Sp (A)$ of $A \in \A$ is considered in $\uA$. Let $A$ be an element of a $C^*$-algebra. Then, $A$ is said to be \emph{self-adjoint} if $A^*=A$, and \emph{positive}, denoted by $A \geq 0$, if it is self-adjoint and $\Sp (A)$ consists of non-negative real numbers. A mapping $\varphi$ from a $C^*$-algebra $\A$ into another $\B$ is called a \emph{$*$-homomorphism} if it is linear, $\varphi (A^*) = \varphi (A)^*$ and $\varphi (AB) = \varphi (A) \varphi (B)$ for each $A,B \in \A$. It is known that $*$-homomorphisms between $C^*$-algebras are automatically bounded and norm-one. An injective $*$-homomorphism is called a \emph{$*$-isomorphism}. We note that each $*$-isomorphism is an isometry. If there exists a $*$-isomorphism from $\A$ onto $\B$, then  $\A$ is \emph{$*$-isomorphic} to $\B$. For basics of $C^*$-algebras, for example, the readers are referred to \cite{Arv76,Bla06,Dav96,KR97a,KR97b,Ped18}.

\subsection{States and representations}

Let $\A$ be a $C^*$-algebra. For each $\rho \in \A^*$, define another element $\rho^* \in \A^*$ by $\rho^* (A) = \overline{\rho (A^*)}$ for each $A \in \A$. A functional $\rho \in \A^*$ is said to be \emph{hermitian} if $\rho^* = \rho$, and \emph{positive} if $\rho (A) \geq 0$ whenever $A \in \A$ and $A \geq 0$. Each positive functional is hermitian. The \emph{state space} $\mathcal{S}(\A)$ of $\A$ is given by
\[
\mathcal{S}(\A) = \left\{ \rho \in \A^* : \|\rho\| = \lim_\lambda \rho (E_\lambda ) = 1\right\}
\]
for any fixed approximate identity $(E_\lambda )_\lambda$ in $\A$; see \cite[Lemmas I.9.5 and I.9.9]{Dav96}. An element of $\mathcal{S}(\A)$ is called a \emph{state} of $\A$. If $\rho \in \ext (\mathcal{S}(\A))$, then $\rho$ is said to be \emph{pure}. The set of all pure states of $\A$, that is, $\ext (\mathcal{S}(\A))$, is called the \emph{pure state space} of $\A$. We note that $\mathcal{S}(\A)$ is a (possibly not weakly$^*$ closed) face of $B_{\A^*}$ and each pure state of $\A$ belongs to $\ext (B_{\A^*})$. As is well-known, each (pure) state of a non-unital $C^*$-algebra $\A$ extends uniquely to a (pure) state of its unitization $\uA$; see~\cite[Proposition II.6.2.5 and Subsection II.6.3]{Bla06}.

Let $B(\HH)$ denote the $C^*$-algebra of all bounded linear operators on a complex Hilbert space $\HH$. A $*$-homomorphism from a $C^*$-algebra $\A$ into some $B(\HH)$ is called a \emph{representation} of $\A$. An injective representation of $\A$ on a Hilbert space $\HH$ (that is, a $*$-isomorphism from $\A$ into $B(\HH)$) is called a \emph{faithful} representation of $\A$. Let $\pi :\A \to B(\HH)$ is a representation of $\A$. Then, $\pi$ is said to be \emph{irreducible} if the $C^*$-subalgebra $\pi (\A)$ of $B(\HH)$ \emph{acts irreducibly} on $\HH$, that is, there is no nontrivial subspace of $\HH$ invariant under $\pi (\A)$. If $\pi_1:\A \to B(\HH_1)$ and $\pi_2 :\A \to B(\HH_2)$ are representations of $\A$, then $\pi_1$ and $\pi_2$ are \emph{unitarily equivalent} if exists an isometric isomorphism $U:\HH_1 \to \HH_2$ such that $U\pi_1 (A)U^* = \pi_2 (A)$ for each $A \in \A$.

By the Gelfand-Naimark-Segal (GNS) construction, each state $\rho$ of a $C^*$-algebra $\A$ generates a representation $\pi_\rho :\A \to B(\HH_\rho)$ and a unit cyclic vector $x_\rho \in \HH_\rho$ for $\pi_\rho (\A)$ such that $\rho (A) = \langle \pi_\rho (A)x_\rho ,x_\rho \rangle$ for each $A \in \A$. As an application of GNS construction, we can always find a faithful representation $\pi :\A \to B(\HH)$. In this sense, a $C^*$-algebra is considered as a norm-closed $*$-subalgebra of some $B(\HH)$. In the rest of this paper, ``Hilbert spaces'' always mean \emph{complex} Hilbert spaces if we work with $C^*$-algebras. It is known that the representation $\pi_\rho$ is irreducible if and only if $\rho$ is pure.

\subsection{Spectra and Primitive spectra}

Let $\A$ be a $C^*$-algebra. Then, an ideal $\mathcal{I}$ of $\A$ is said to be \emph{primitive} if $\mathcal{I}=\ker \pi$ for some irreducible representation $\pi$ of $\A$. The \emph{primitive spectrum} $\Prim (\A)$ of $\A$ is the set of all primitive ideals of $\A$ equipped with the hull-kernel topology, where the hull-kernel closure of $S \subset \Prim (\A)$ is given by $\{ \mathcal{I} \in \Prim (\A) : \bigcap \{ \mathcal{I}': \mathcal{I}' \in S\} \subset \mathcal{I}\}$. The \emph{spectrum} $\hat{\A}$ is defined as the set of all unitary equivalence classes of irreducible representations of $\A$ that equips with the pullback topology with respect to the mapping $\pi \mapsto \ker \pi$ from $\hat{\A}$ onto $\Prim (\A)$. As was noted in \cite[Theorem 3.13.2]{Ped18}, each irreducible representation is unitarily equivalent to $\pi_\rho$ for some pure state $\rho$ of $\A$. Spectra and primitive spectra of $C^*$-algebras are basic tools in the representation theory of $C^*$-algebras. One of the most successful applications of (primitive) spectra to the theory of operator algebras is characterizations of Type I $C^*$-algebras which is known as the Glimm-Sakai theorem.

Let $\A$ be a $C^*$-algebra. Then, $A \in \A$ is said to be \emph{abelian} if the norm-closure of $A\A A$ (hereditary $C^*$-algebra) is abelian. If every quotient of $\A$ contains a nonzero abelian element, then $\A$ is \emph{of internally type I}. Meanwhile, $\A$ is \emph{of bidual type I} if $\A^{**}$ is a type I von Neumann algebra. By a deep result by Glimm~\cite{Gli61} and Sakai~\cite{Sak66,Sak67}, these two definitions of type I $C^*$-algebras are mutually equivalent and they are characterized in terms of irreducible representations. Let $\pi:\A \to B(\HH)$ be an irreducible representation. Then, $\pi$ is a \emph{CCR representation} if $K(\HH) = \pi (\A)$, and \emph{GCR representation} if $K(\HH) \subset \pi (\A)$. If every irreducible representation of $\A$ is CCR (or GCR), then $\A$ is called a \emph{CCR} (or \emph{GCR}) $C^*$-algebra. The following is (a part of) the Glimm-Sakai theorem.
\begin{GST}

Let $\A$ be a $C^*$-algebra. Then, the following are equivalent:
\begin{itemize}
\item[{\rm (i)}] $\A$ is of internally type I.
\item[{\rm (ii)}] $\A$ is of bidual type I.
\item[{\rm (iii)}] $\A$ is GCR.
\end{itemize}
Moreover, (i) to (iii) implies the following:
\begin{itemize}
\item[{\rm (iv)}] The mapping $\pi \mapsto \ker \pi$ from $\hat{\A}$ onto $\Prim (\A)$ is injective.
\end{itemize}
If $\A$ is separable, then (i) to (iv) are equivalent.
\end{GST}
The separable case was due to Glimm, and the general case was completed by Sakai. We note that the separability is required to obtain (i) to (iii) from (iv). Indeed, if $\A$ is a non-separable simple $C^*$-algebra, then the implication (iv) to (iii) is (essentially the same as) the Naimark problem.


\section{Closure spaces}\label{Sect:Cl}

We begin with drawing the framework of our study. For a set $K$, the symbol $2^K$ stands for the power set of $K$. First, we introduce the concept of (Fr\'{e}chet (V)) closure operators, which is a basis for the theory developed in the present paper; see~\cite{Rio59,Sie52} for backgrounds. 
\begin{definition}\label{F-cl}

Let $K$ be a set. Then, $c:2^K \to 2^K$ is called a \emph{(Fr\'{e}chet (V)) closure operator} if it satisfies the following conditions:
\begin{itemize}
\item[(i)] $c(\emptyset ) = \emptyset$.
\item[(ii)] $S \subset c(S)$ for each $S\subset K$.
\item[(iii)] $c(c(S)) = c(S)$ for each $S \subset K$.
\item[(iv)] $c(S_1) \subset c(S_2)$ whenever $S_1 \subset S_2 \subset K$.
\end{itemize}
The ordered pair $(K,c)$ is called a \emph{(Fr\'{e}chet (V)) closure space}.
\end{definition}
\begin{remark}

There is another well-known concept of closure operators. Let $K$ be a set. Then, $c:2^K \to 2^K$ is called a \emph{preclosure operator} or \emph{\v{C}ech closure operator} if it satisfies the following conditions:
\begin{itemize}
\item[(i)] $c(\emptyset ) = \emptyset$.
\item[(ii)] $S \subset c(S)$ for each $S\subset K$.
\item[(iii)] $c(S_1 \cup S_2) = c(S_1)\cup C(S_2)$ for each subsets $S_1,S_2$ of $K$.
\end{itemize}
Although preclosure operators are also called closure operators in the literature, in this paper, the term ``a closure opeartor'' always means a Fr\'{e}chet (V) closure operator.
\end{remark}
In a closure space, we can consider closed sets in a natural manner.

\begin{definition}

Let $(K,c)$ be a closure space, and let $S \subset K$. Then, $S$ is said to be \emph{closed} if $c(S) = S$. Let $\mathfrak{C}(K,c)$ denote the family of all closed subsets of $K$, that is, $\mathfrak{C}(K,c) = \{ S \subset K : c(S) = S\}$.
\end{definition}
\begin{remark}

Let $(K,c)$ be a closure space. Then, the identity
\[
c(S) = \bigcap \{ C \in \mathfrak{C}(K,c) : S \subset C\}
\]
holds for each $S \subset K$. Indeed, if $C \in \mathfrak{C}(K,c)$ and $S \subset C$, then $c(S) \subset c(C) = C$. The converse follows from the fact that $S \subset c(S) \in \mathfrak{C}(K,c)$.
\end{remark}
\begin{remark}

Let $(K,c)$ be a closure space. Then, $\emptyset ,K \in \mathfrak{C}(K,c)$ by $c (\emptyset ) = \emptyset$ and $K \subset c(K) \subset K$.
\end{remark}
Let $K$ be a set, and let $c:2^K \to 2^K$. Then, as is well-known, the family $\{ S \subset K : c(S) = S\}$ satisfies the axioms of closed sets if and only if it satisfies the following conditions:
\begin{itemize}
\item[(i)] $c(\emptyset ) = \emptyset$.
\item[(ii)] $S \subset c(S)$ for each $S\subset K$.
\item[(iii)] $c(c(S)) = c(S)$ for each $S \subset K$.
\item[(iv)] $c(S_1 \cup S_2) = c(S_1) \cup c(S_2)$ for each subsets $S_1,S_2$ of $K$.
\end{itemize}
This set of conditions is known as the Kuratowski closure axioms.
\begin{definition}

Let $(K,c)$ be a closure space. Then, $(K,c)$ is said to be \emph{topologizable} if the closure operator $c$ satisfies the Kuratowski closure axioms. If $(K,c)$ is topologizable, then it is called a \emph{topological (closure) space} and $K$ is equipped with the topology induced by $\mathfrak{C}(K,c)$.
\end{definition}
\begin{remark}

If $(K,c)$ is a closure space, and if $S_1,S_2$ are subsets of $K$, then the monotonicity of $c$ implies that $c(S_1 \cup S_2) \supset c(S_1) \cup c(S_2)$. Hence, $(K,c)$ is topologizable if and only if $c(S_1 \cup S_2) \subset c(S_1) \cup c(S_2)$ for each subsets $S_1,S_2$ of $K$. It should be mentioned that this was already noted in the restricted case of geometric structure spaces; see~\cite[Lemma 3.6]{Tan22b}.
\end{remark}
The concept of closure spaces generalizes that of topological spaces. Hence, it is natural to consider appropriate morphisms between them.
\begin{definition}

Let $(K,c)$ and $(L,d)$ be closure spaces, and let $f:K \to L$. Then, $f$ is said to be \emph{continuous} if $f(c(S))\subset d(f(S))$ for each $S \subset K$. If $f$ is bijecitive, and if both $f$ and $f^{-1}$ are continuous, then $f$ is called a \emph{homeomorphism}. The closure spaces $(K,c)$ and $(L,d)$ are said to be \emph{homeomorphic}, denoted by $(K,c) \simeq (L,d)$, if there exists a homeomorphism between $(K,c)$ and $(L,d)$.
\end{definition}
\begin{remark}

Let $(K,c)$ and $(L,d)$ be closure spaces, and let $f :K \to L$. Then, $f$ is continuous if and only if $f^{-1}(S') \in \mathfrak{C}(K,c)$ whenever $S'\in \mathfrak{C}(L,d)$, and $f$ is a homeomorphism if and only if $f(c(S)) = d(f(S))$ for each $S \subset K$. Although these facts were essentially already noted in \cite[Remarks 3.10 and 3.11]{Tan23a}, we realize them here for the reader's convenience.

Suppose that $f$ is continuous and $S' \in \mathfrak{C}(Y,d)$. Then, it follows from 
\[
f(c(f^{-1}(S'))) \subset d(f(f^{-1}(S'))) \subset d(S') = S'
\]
that $c(f^{-1}(S')) \subset f^{-1}(S')$. Hence, $f^{-1}(S') \in \mathfrak{C}(K,c)$. Conversely, if $f^{-1}(S') \in \mathfrak{C}(K,c)$ whenever $S'\in \mathfrak{C}(L,d)$, and if $S \subset K$, then we infer that $d(f(S)) \in \mathfrak{C}(L,d)$ and $f^{-1}(d(f(S)) \in \mathfrak{C}(K,c)$. Since $S \subset f^{-1}(f(S)) \subset f^{-1}(d(f(S)))$, we obtain
\[
c(S) \subset c(f^{-1}(d(f(S)))) = f^{-1}(d(f(S))) ,
\]
which implies that $f(c(S)) \subset d(f(S))$. This shows that $f$ is continuous.

Next, we assume that $f$ is a homeomorphism. Then, it follows from $f(c(S)) \subset d(f(S))$ and $f^{-1}(d(f(S))) \subset c(f^{-1}(f(S)) = c(S)$ that $f(c(S)) = d(f(S))$ for each $S \subset K$. Conversely, suppose that $f(c(S)) = d(f(S))$ for each $S \subset K$. If $S' \in \mathfrak{C}(L,d)$, then
\[
f(c(f^{-1}(S'))) = d(f(f^{-1}(S'))) = d(S') = S'
\]
and $c(f^{-1}(S')) = f^{-1}(S')$. Therefore, $f$ is continuous. Meanwhile, if $S \in \mathfrak{C}(K,c)$, then $d(f(S)) = f(c(S)) = f(S)$, that is, $(f^{-1})^{-1}(S) = f(S) \in \mathfrak{C}(L,d)$. Thus, $f^{-1}$ is also continuous.
\end{remark}
\begin{remark}

Let $(K,c)$ be a closure space. Then, the identity mapping $\id_K$ on $K$ is a homeomorphism on $K$. This is realized by the equation $\id_K (c(S)) = c(S) = c(\id_K (S))$ for each $S \subset K$.
\end{remark}
The following proposition is what is needed for defining the category of closure spaces.
\begin{proposition}\label{category-lemma}

Let $(K,c)$, $(L,d)$ and $(M,e)$ be closure spaces, and let $f:K \to L$ and $g:L \to M$. If $f$ and $g$ are continuous, then $g \circ f : K \to M$ is continuous.
\begin{proof}
Let $S \subset K$. Then, it follows that
\[
(g\circ f)(c(S)) = g(f(c(S))) \subset g(d(f(S))) \subset e(g(f(S))) = e((g\circ f)(S)) .
\]
Hence, $g \circ f$ is continuous.
\end{proof}
\end{proposition}
Now, we introduce the category of closure spaces.

\begin{definition}

The category of (Fr\'{e}chet (V)) closure spaces, denoted as $\FCl$, is the category whose objects are (Fr\'{e}chet (V)) closure spaces. The morphisms between closure spaces $(K,c)$ and $(L,d)$ are the continous mappings, and the composition of morphisms is the composition of mappings. The full subcategory of $\FCl$ whose objects are topological closure spaces is denoted by $\FTop$.
\end{definition}
\begin{remark}

The isomorphisms in $\FCl$ are just the homeomorphisms.
\end{remark}
As the following proposition says, $\FTop$ has some desirable properties. Recall that a full subcategory of a category is said to be \emph{strictly full} if it is isomorphism-closed. Let $\Top$ denote the category of topological spaces.

\begin{proposition}\label{strictly-full}

$\FTop$ is strictly full subcategory of $\FCl$. Moreover, $\FTop$ and $\Top$ are isomorphic.
\begin{proof}
The proof of the former statement is essentially the same as that of \cite[Theorem 3.11]{Tan22b}. Let $(K,c),(L,d)$ are closure spaces. Suppose that $(K,c)$ is topologizable, and that $f:K \to L$ be a homeomorphism. Since $\FTop$ is full, it is sufficient to show that $L$ is topologizable. Let $S'_1,S'_2 \subset L$. Then, it follows from
\begin{align*}
f^{-1}(d(S'_1\cup S'_2)) = c(f^{-1}(S'_1 \cup S'_2)) 
&= c(f^{-1}(S'_1) \cup f^{-1}(S'_2)) \\
&= c(f^{-1}(S'_1)) \cup c(f^{-1}(S'_2)) \\
&= f^{-1}(d(S'_1)) \cup f^{-1}(d(S'_2))
\end{align*}
that
\begin{align*}
d(S'_1 \cup S'_2) = f(f^{-1}(d(S'_1)) \cup f^{-1}(d(S'_2))) 
&= f(f^{-1}(d(S'_1))) \cup f(f^{-1}(d(S'_2))) \\
&= d(S'_1) \cup d(S'_2) .
\end{align*}
This shows that $L$ is topologizable.

Now, we note that each topological space $(K,\mathfrak{T})$ can be viewed as a closure space by setting $c_\mathfrak{T}(S) = \overline{S}^\mathfrak{T}$ for each $S \subset K$. Moreover, $(K,c_\mathfrak{T})$ is topologizable and $\mathfrak{T} = \{ K \setminus S : S \in \mathfrak{C}(K,c_\mathfrak{T})\}$. Meanwhile, if $(K,c)$ is a topological closure space, then $\mathcal{I}_c = \{ K \setminus S: S \in \mathfrak{C}(K,c)\}$ is a topology on $K$ and $c = c_{\mathfrak{T}_c}$. Remark also that $f:(K,\mathfrak{T}) \to (L,\mathfrak{U})$ is continuous if and only if $f:(K,c_\mathfrak{T}) \to (L,c_\mathfrak{U})$ is continuous. 

Define functors $F$ from $\Top$ to $\FTop$ and $G$ from $\FTop$ to $\Top$ by $F((K,\mathfrak{T})) = (K,c_\mathfrak{T})$ for each topological space $(K,\mathfrak{T})$, $F(f)=f$ for each continuous mapping $f$ between topological spaces, $G((K,c)) = (K,\mathfrak{T}_c)$ for each topological closure space $(K,c)$, and $G(f)=f$ for each continuous mapping $f$ between topological closure spaces. Then, $FG$ and $GF$ are the identity functors on $\FTop$ and $\Top$, respectively. Thus, $\Top$ and $\FTop$ are (covariantly) isomorphic.
\end{proof}
\end{proposition}
As was shown in Proposition~\ref{strictly-full}, each topological space $(K,\mathfrak{T})$ naturally corresponds to the topological closure space $(K,c_\mathfrak{T})$. In the rest of this paper, we identify $(K,\mathfrak{T})$ with $(K,c_\mathfrak{T})$ unless otherwise stated.

Next, we consider subspaces of a closure space.
\begin{definition}

Let $(K,c)$ be a closure space, and let $A \subset K$. Define a mapping $c_A : 2^A \to 2^A$ by $c_A(S) = c(S) \cap A$ for each $S \subset A$. Then, the ordered pair $(A,c_A)$ is called a \emph{subspace} of $(K,c)$.
\end{definition}
It will turn out that subspaces of closure spaces are indeed closure spaces.

\begin{proposition}\label{sub-closure-space}

Let $(K,c)$ be a closure space, and let $A \subset K$. Then, $(A,c_A)$ is a closure space. Moreover, the natural injection $\iota_A : A \to B$ is continuous whenever $A \subset B \subset K$, where $\iota_A$ is defined by $\iota_A (t) = t$ for each $t \in A$.
\begin{proof}
Since $c(\emptyset )=\emptyset$, we have $c_A(\emptyset ) = \emptyset$. If $S \subset A$, then it follows from $S \subset c(S)$ that
\[
S = S \cap A \subset c(S) \cap A = c_A (S) .
\]
We also note that $c(c_A(S)) \subset c(S) \cap c(A)$ since $c_A(S) = c(S)\cap A \subset c(S) \cap c(A)$ and
\[
c(c(S) \cap c(A)) \subset c(c(S)) \cap c(c(A)) = c(S) \cap c(A),
\]
which implies that
\[
c_A (c_A (S)) = c(c_A(S)) \cap A \subset c(S) \cap c(A) \cap A =  c(S) \cap A = c_A(S) .
\]
Finally, if $S_1 \subset S_2 \subset A$, then
\[
c_A(S_1) = c(S_1) \cap A \subset c(S_2) \cap A = c_A(S_2)
\]
by $c(S_1) \subset c(S_2)$. Therefore, $c_A$ is a closure operator on $A$.

The continuity of $\iota_A$ follows from the fact that
\[
\iota_A (c_A(S)) = c_A(S) = c(S) \cap A \subset c(S) \cap B = c_B(S) = c_B (\iota_A (S))
\]
for each $S \subset A$.
\end{proof}
\end{proposition}
\begin{remark}

Let $(K,c)$ be a closure space, and let $A \subset K$. If $A$ is closed, then $c(S) \subset c(A) = A$ whenever $S \subset A$. Hence,
\[
c_A (S) = c(S) \cap A = c(S)
\]
for each $S \subset A$, that is, $c_A = c|2^A$.
\end{remark}
As in the case of topological spaces, continuous functions between closure spaces restrict to continuous functions between subspaces.
\begin{proposition}\label{subspace-homeo}

Let $(K,c)$ and $(L,d)$ be closure spaces, and let $A \subset K$. If $f:K \to L$ is continuous (or a homeomorphism), then $f|A:A \to f(A)$ is continuous (or a homeomorphism).
\begin{proof}
Let $S \subset A$. If $f$ is continuous, then
\[
f(c_A(S)) = f(c(S) \cap A) \subset f(c(S)) \cap f(A) \subset d(f(S)) \cap f(A) = c_{f(A)}(f(S)) .
\]
Hence, $f|A:A \to f(A)$ is continuous. In the case that $f$ is homeomorphism, we obtain
\[
f(c_A(S)) = f(c(S) \cap A) = f(c(S)) \cap f(A) = d(f(S)) \cap f(A) = c_{f(A)}(f(S)) ,
\]
which implies that $f|A:A \to f(A)$ is a homeomorphism.
\end{proof}
\end{proposition}
The following lemma gives a way of constructing the \emph{pullback} closure operator.
\begin{lemma}\label{pullback}

Let $K$ be a set, let $(L,d)$ be a closure space, and let $f:K \to L$. Define a mapping $c:2^K \to 2^K$ by $c(S) = f^{-1}(d(f(S)))$ for each $S \subset K$. Then, the following hold:
\begin{itemize}
\item[{\rm (i)}] $c$ is a closure operator on $K$ and $f$ is continuous with respect to the closure operators $c$ and $d$.
\item[{\rm (ii)}] If $(L,d)$ is topologizable, then so is $(K,c)$.
\item[{\rm (iii)}] If $f$ is surjective, then $f(c(S)) = d(f(S))$ for each $S\subset K$. In particular, $f$ is a closed mapping with respect to the closure operators $c$ and $d$.
\item[{\rm (iv)}] If $f$ is bijective, then it is a homeomorphism with respect to the closure operators $c$ and $d$.
\end{itemize}
\begin{proof}
(i) First, we note that $f(c(S)) = f(f^{-1}(d(f(S)))) \subset d(f(S))$ for each $S \subset K$. From this, it is sufficient to prove that $c$ is a closure operator on $K$. Since $d(\emptyset ) = \emptyset$, we have $c(\emptyset ) = f^{-1}(d(f(\emptyset ))) = \emptyset$. If $S \subset K$, then $S \subset f^{-1}(f(S)) \subset f^{-1}(d(f(S)) = c(S)$ and
\[
c(c(S)) = f^{-1}(d(f(c(S)))) \subset f^{-1}(d(d(f(S)))) = f^{-1}(d(f(S))) = c(S) .
\]
Moreover,
\[
c(S_1) = f^{-1}(d(f(S_1))) \subset f^{-1}(d(f(S_2))) = c(S_2)
\]
whenever $S_1 \subset S_2 \subset K$. This shows that $c$ is a closure operator on $K$.

(ii) Suppose that $(L,d)$ is topologizable, and that $S_1,S_2$ are subsets of $K$. Then, it follows that
\begin{align*}
c(S_1 \cup S_2) = f^{-1}(d(f(S_1 \cup S_2))) 
&= f^{-1}(d(f(S_1))) \cup f^{-1}(d(f(S_2))) \\
&= c(S_1) \cup c(S_2) .
\end{align*}
Hence, $(K,c)$ is topologizable.

(iii) If $f$ is surjective, then $f(f^{-1}(S')) = S'$ for each $S' \subset L$. Therefore, $f(c(S)) = d(f(S))$ for each $S \subset K$. In particular, if $S$ is closed, then $f(S) = f(c(S)) = d(f(S))$, that is, $f(S)$ is closed.

(iv) If $f$ is bijective, then it is continuous by (i) and closed by (iii). Thus, $f$ is a homeomorphism.
\end{proof}
\end{lemma}
The preceding lemma is useful for constructing new closure spaces from a given closure space. The following is the basis for our construction.
\begin{definition}\label{2p-eq}

Let $(K,c)$ be a closure space. Define binary relations on $K$ by declaring that $t_1 \rightleftarrows t_2$ if $t_1=t_2$ or $\{t_1,t_2\} \not \in \mathfrak{C}(K,c)$, and that $t_1 \sim t_2$ if there exists a finite subset $\{ s_1,\ldots ,s_n\}$ of $K$ such that $t_1 \rightleftarrows s_1 \rightleftarrows \cdots \rightleftarrows s_n \rightleftarrows t_2$. 
\end{definition}
\begin{remark}\label{arrows}

Let $(K,c)$ be a closure space, and let $t,t_1,t_2 \in K$. Then, the following hold:
\begin{itemize}
\item[(i)] $t \rightleftarrows t$.
\item[(ii)] If $t_1 \rightleftarrows t_2$, then $t_2 \rightleftarrows t_1$.
\item[(iii)] If $t_1 \rightleftarrows t_2$, then $t_1 \sim t_2$.
\end{itemize}
The statements (i) and (ii) are obvious from the definition. Moreover, if $t_1 \rightleftarrows t_2$, then $\{t_1\}$ is a finite subset of $\mathfrak{S}(X)$ and $t_1 \rightleftarrows t_1 \rightleftarrows t_2$, which implies that $t_1 \sim t_2$. This proves (iii).

In particular, if ``$\rightleftarrows$'' is transitive, then it becomes an equivalence relation on $K$, in which case, ``$\rightleftarrows$'' is equivalent to ``$\sim$.''
\end{remark}
\begin{remark}\label{equal}

Let $(K,c)$ be a closure space. If $\{t_1,t_2\} \in \mathfrak{C}(K,c)$ for each $t_1,t_2 \in K$, then ``$\sim$'' is equivalent to ``$=$.'' The converse is also true if $(K,c)$ is $T_1$. Indeed, if $\{ t_1,t_2\} \in \mathfrak{C}(K,c)$ for each $t_1,t_2 \in K$, and if $t_1 \rightleftarrows t_2$, then $t_1 = t_2$. This shows that $t_1 \sim t_2$ if and only if $t_1=t_2$. Conversely, if $(K,c)$ is $T_1$, and if $\{t_1,t_2\} \not \in \mathfrak{C}(K,c)$ for some $t_1,t_2 \in K$, then $t_1 \sim t_2$. We note that $t_1 \neq t_2$ since each singleton is closed.
\end{remark}
It will turn out that the relation ``$\sim$'' is an equivalence relation.

\begin{proposition}\label{equiv}

Let $(K,c)$ be a closure space. Then, ``$\sim$'' is an equivalence relation on $K$.
\begin{proof}
It follows from Remark~\ref{arrows} that $t \sim t$ for each $t \in K$ and $t_1 \sim t_2$ implies that $t_2 \sim t_1$. Now, suppose that $t_1 \sim t_2$ and $t_2 \sim t_3$. Then, there exist finite subsets $\{s_1,\ldots ,s_n\}$ and $\{s'_1,\ldots ,s'_m\}$ of $K$ such that $t_1 \rightleftarrows s_1 \rightleftarrows \cdots \rightleftarrows s_n \rightleftarrows t_2$ and $t_2 \rightleftarrows s'_1 \rightleftarrows \cdots \rightleftarrows s'_m \rightleftarrows t_3$. In particular, $\{ s_1,\ldots ,s_n ,t_2,s'_1,\ldots ,s'_m \}$ is a finite subset of $K$ and
\[
t_1 \rightleftarrows s_1 \rightleftarrows \cdots \rightleftarrows s_n \rightleftarrows t_2 \rightleftarrows s'_1 \rightleftarrows \cdots \rightleftarrows s'_m \rightleftarrows t_3 .
\]
Therefore, $t_1 \sim t_3$. This completes the proof.
\end{proof}
\end{proposition}
Using the equivalence relation ``$\sim$'' on a given closure space, we obtain two natural way of constructing new closure spaces. In the case of topological spaces, the most natural candidate of topologies on a quotient set is the quotient topology. However, we will go to another way to obtain closure spaces which are more compatible with nonlinear geometry.
\begin{definition}\label{trans-def}

Let $(K,c)$ be a closure space, let $K^\mathfrak{G}=K/\!\sim ~ = \{ C(t) : t\in K\}$ and $K^\mathfrak{P} = \{ c(C(t)) : t\in K \}$. Define mappings $c^\mathfrak{G} :2^{K^\mathfrak{G}} \to 2^{K^\mathfrak{G}}$ and $c^\mathfrak{P} : 2^{K^\mathfrak{P}} \to 2^{K^\mathfrak{P}}$ by
\[
c^\mathfrak{P}(T) = \left\{ c(C(t)) : c(C(t)) \subset c(\{ s \in K : c(C(s)) \in T\} ) \right\}
\]
for each $T \subset K^\mathfrak{P}$, and
\[
c^\mathfrak{G}(S) = \left\{ C(t) : c(C(t)) \in c^\mathfrak{P} ( \{ c(C(s)) : C(s) \in S \}) \right\}
\]
for each $S \subset K^\mathfrak{G}$. The ordered pairs $(K^\mathfrak{G},c^\mathfrak{G})$ and $(K^\mathfrak{P},c^\mathfrak{P})$ are called \emph{$\mathfrak{G}$-transform} and \emph{$\mathfrak{P}$-transform} of $(K,c)$, respectively. Define a mapping $\delta_K : K^\g \to K^\p$ by $\delta_K  (C(t)) = c(C(t))$ for each $t \in K$. 
\end{definition}
We present the first properties of $\g$- and $\p$-transforms of closure spaces.
\begin{proposition}\label{two-closures}

Let $(K,c)$ be a closure space. Then, $(K^\g,c^\g)$ and $(K^\p,c^\p)$ are closure spaces. Moreover, $c^\g (S) = \delta_K^{-1}(c^\p (\delta_K (S)))$ for each $S\subset K^\g$.
\begin{proof}
The identity $c^\g (S) = \delta_K^{-1}(c^\p (\delta_K (S)))$ for each $S\subset K^\g$ follows from the definitions of $c^\g$ and $c^\p$. Combining this with Lemma~\ref{pullback}, it is sufficient to show that $c^\p$ is a closure operator on $K^\p$. First, we note that $\{ s \in K : c(C(s)) \in \emptyset \} = \emptyset$, which together with the fact that $t \in C(t) \subset c(C(t))$ for each $t$ implies $c^\p (\emptyset ) = \emptyset$. Moreover, if $T_1 \subset T_2 \subset K^\p$, then
\[
\{ s \in K : c(C(s)) \in T_1 \} \subset \{ s \in K : c(C(s)) \in T_2\} 
\]
and $c^\p (T_1) \subset c^\p (T_2)$ by the monotonicity of $c$. Now, let $T \subset K^\p$. If $c(C(t)) \in T$, and if $s \in C(t)$, then $C(s) = C(t)$ and $c(C(s))=c(C(t)) \in T$. This shows that $C(t) \subset \{ s \in K : c(C(s)) \in T\}$. Hence, $c(C(t)) \subset c(\{ s \in K : c(C(s)) \in T\})$ and $c(C(t)) \in c^\mathfrak{P}(T)$, that is, $T \subset c^\p (T)$, which in turn implies that
\[
c(\{ s \in K : c(C(s)) \in T\} ) \subset c(\{ s \in K : c(C(s)) \in c^\p (T)\}) .
\]
Meanwhile, if $c(C(t)) \in c^\p (T)$, then
\[
t \in C(t) \subset c(C(t)) \subset c(\{ s \in K : c(C(s)) \in T\} )
\]
by the definition of $c^\p$. Therefore,
\[
\{ s \in K : c(C(s)) \in c^\p (T)\} \subset c(\{ s \in K : c(C(s)) \in T\} )
\]
and $c(\{ s \in K : c(C(s)) \in c^\p (T)\}) \subset c(\{ s \in K : c(C(s)) \in T\} )$, that is,
\[
c(\{ s \in K : c(C(s)) \in T\} ) = c(\{ s \in K : c(C(s)) \in c^\p (T)\}) .
\]
This proves that $c^\p (c^\p (T)) = c^\p (T)$.
\end{proof}
\end{proposition}
\begin{remark}\label{T0}

Let $(K,c)$ be a closure space. Then, $(K^\g,c^\g)$ is always $T_0$, that is, whenever $c(C(t_1)) \neq c(C(t_2))$, there exists an $F \in \mathfrak{C}(K^\p,c^\p)$ such that $c(C(t_1)) \in F$ and $c(C(t_2)) \not \in F$, or $c(C(t_2)) \in F$ and $c(C(t_2)) \not \in F$. Indeed, if $s,t \in K$ and $c(C(s)) = c(C(t))$, then $s \in C(s) \subset c(C(s)) = c(C(t))$, which implies that $\{ s \in K : c(C(s)) = c(C(t)) \} \subset c(C(t))$. Hence,
\[
c^\p (\{ c(C(t)) \}) \subset \{ c(C(t')) : c(C(t')) \subset c(C(t)) \}
\]
for each $t \in K$. In particular, if $t_1,t_2 \in K$, $c(C(t_1)) \in c^\p (\{ c(C(t_2))\})$, and $c(C(t_2)) \in c^\p (\{ c(C(t_1))\})$, then $c(C(t_1))=c(C(t_2))$. This shows that $(K^\p,c^\p)$ is $T_0$.

Meanwhile, $(K^\g,c^\g)$ is $T_0$ if and only if it is homeomorphic to $(K^\p,c^\p)$, in which case, $\delta_K :K^\g \to K^\p$ is a homeomorphism. Indeed, if $(K^\g,c^\g)$ is $T_0$, and if $C(t_1) \neq C(t_2)$, then there exists an $F \in \mathfrak{C}(K^\g,c^\g)$ such that $C(t_1) \in F$ and $C(t_2) \not \in F$, or $C(t_1) \not \in F$ and $C(t_2) \in F$. Replacing $C(t_1)$ with $C(t_2)$ if necessary, we may assume that $C(t_1) \in F$ and $C(t_2) \not \in F$. Since $c^\g (F)=F$, it follows that $c^\g (\{C(t_1)\}) \subset F$, which implies that
\begin{align*}
C(t_2) \not \in c^\g (\{ C(t_1)\}) 
&= \{ C(t) : c(C(t)) \in c^\p (\{ c(C(t_1))\}) \} .
\end{align*}
Combining this with $c(C(t_1)) \in c^\p (\{ c(C(t_1))\})$, we have $c(C(t_1)) \neq c(C(t_2))$, that is, $\delta_K$ is injective. Therefore, $\delta_K :K^\g \to K^\p$ is a homeomorphism by Lemma~\ref{pullback} (iv) and Proposition~\ref{two-closures}. The converse immediately follows from the first paragraph.
\end{remark}
Althought the morphisms in $\FCl$ are the continuous mappings, in the context of $\g$- and $\p$-transforms of closure spaces, a weaker concept of continuity plays important roles.
\begin{definition}

Let $(K,c)$ and $(L,d)$ be closure spaces, and let $f:K \to L$. Then, $f$ is said to be \emph{finitely continuous} if $f(c(S)) \subset d(f(S))$ for each finite subset $S$ of $K$. A finitely continuous bijection with finitely continuous inverse is called a \emph{finite-homeomorphism}.
\end{definition}
We naturally expect that homeomorphisms between closure spaces preserve the equivalence of elements. In fact, as the following lemma says, finite-homeomorphisms already have this property.
\begin{lemma}\label{finite-conti}

Let $(K,c)$ and $(L,d)$ be closure spaces, and let $f: K \to L$ be finitely continuous injection. Then, $t_1 \rightleftarrows t_2$ implies $f(t_1) \rightleftarrows f(t_2)$. Consequently, $t_1 \sim t_2$ implies $f(t_1) \sim f(t_2)$ .
\begin{proof}
Suppose that $t_1 \rightleftarrows t_2$ and $t_1 \neq t_2$. Then, $\{t_1,t_2\} \not \in \mathfrak{C}(K,c)$. Since $f$ is finitely continuous, it follows that
\[
f(c(\{t_1,t_2\})) \subset d(f(\{t_1,t_2\})) .
\]
If $\{f(t_1),f(t_2)\} \in \mathfrak{C}(L,d)$, then $f(c(\{t_1,t_2\})) \subset f(\{t_1,t_2\})$, which together with the injectivity of $f$ implies that
\[
c(\{t_1,t_2\}) = f^{-1}(f(c(\{t_1,t_2\}))) \subset f^{-1}(f(\{t_1,t_2\})) = \{t_1,t_2\} .
\]
However, this contradicts $\{t_1,t_2\} \not \in \mathfrak{C}(K,c)$. Hence, $\{f(t_1),f(t_2)\} \not \in \mathfrak{C}(L,d)$ and $f(t_1) \rightleftarrows f(t_2)$.
\end{proof}
\end{lemma}
We conclude this section with the following result. It will turn out that the $\g$- and $\p$-transforms of a closure space are invariant under homeomorphisms.
\begin{theorem}\label{two-homeo}

Let $(K,c)$ and $(L,d)$ be closure spaces, and let $f: K \to L$ be a homeomorphism. Define mappings $f^\g :K^\g \to L^\g$ and $f^\p :K^\p \to L^\p$ by $f^\g (C(t)) = C(f(t))$ and $f^\p (c(C(t))) = d(C(f(t)))$ for each $t \in K$. Then, $f^\g$ and $f^\p$ are homeomorphisms with inverses $(f^{-1})^\g$ and $(f^{-1})^\p$. Moreover, $\delta_L \circ f^\g = f^\p \circ \delta_K$.
\begin{proof}
By Lemma~\ref{finite-conti}, we have $f(C(t)) = C(f(t))$ for each $t \in K$. From this, the mapping $f^\g$ is well-defined and bijective. Moreover since $f(c(S)) = d(f(S))$ for each $S \subset K$, it follows that $f(c(C(t))) = d(f(C(t))) = d(C(f(t)))$ for each $t \in K$, which implies that $f^\p$ is well-defined and bijective. By the definitions of $f^\g$ and $f^\p$, we obtain $(f^\g)^{-1}=(f^{-1})^\g$, $(f^\p)^{-1}=(f^{-1})^\p$ and $\delta_L \circ f^\g = f^\p \circ \delta_K$.

Suppose that $T \subset K^\p$. If $c(C(t)) \in c^\p (T)$, then $c(C(t)) \subset c(\{ s\in K :c(C(s)) \in T\})$. Since $f$ is a homeomorphism, we derive that
\begin{align*}
d(C(f(t))) = d(f(C(t))) = f(c(C(t))) 
&\subset f(c(\{ s\in K :c(C(s)) \in T\})) \\
&= d(f(\{s \in K :c(C(s)) \in T\})) .
\end{align*}
Meanwhile, if $s \in K$ and $c(C(s)) \in T$, then $d(C(f(s))) = f^\p (c(C(s))) \in f^\p (T)$, that is, $f(s) \in \{ s' \in L : d(C(s')) \in f^\p (T)\}$. Hence,
\begin{align*}
f^\p (c(C(t))) = d(C(f(t))) 
&\subset d(f(\{s \in K :c(C(s)) \in T\})) \\
&\subset d(\{ s' \in L : d(C(s')) \in f^\p (T)\} .
\end{align*}
This shows that $f^\p (c(C(t))) \in d^\p (f^\p (T))$ whenenver $c(C(t)) \in c^\p (T)$, and that $f^\p (c^\p (T)) \subset d^\p (f^\p (T))$. Therefore, $f^\p$ is continuous. Applying the same argument as above to $f^{-1}$, we can conclude that $f^\p$ is a homeomorphism.

Next, suppose that $S \subset K^\g$. By the preceding paragraph, we have
\begin{align*}
f^\g (c^\g (S)) = f^\g (\delta_K^{-1}(c^\p (\delta_K(S))))
&= f^\g (\delta_K^{-1}(c^\p (((f^\p)^{-1}\circ f^\p \circ \delta_K) (S)))) \\
&= f^\g (\delta_K^{-1}(c^\p (((f^\p)^{-1}\circ \delta_L \circ f^\g ) (S)))) \\
&= f^\g (\delta_K^{-1}((f^\p)^{-1} (d^\p ((\delta_L \circ f^\g ) (S))))) \\
&= f^\g ((f^\p \circ \delta_K )^{-1}(d^\p ((\delta_L \circ f^\g )(S)))) \\
&= f^\g ((\delta_L \circ f^\g)^{-1}(d^\p ((\delta_L \circ f^\g )(S)))) \\
&= f^\g ((f^\g)^{-1}(\delta_L^{-1}(d^\p ((\delta_L \circ f^\g )(S))))) \\
&= \delta_L^{-1}(d^\p (\delta_L (f^\g (S)))) \\
&= d^\g (f^\g (S)) .
\end{align*}
Thus, $f^\g$ is also a homeomorphism.
\end{proof}
\end{theorem}


\section{Geometric spectra and geometric primitive spectra}\label{Sect:BS}

In this section, we focus on special closure spaces that are naturally induced by facial structure of unit balls of Banach spaces.
\begin{definition}

Let $X$ be a (real or complex) nontrivial Banach space, and let $I_F = \bigcup \{ \ker f : f \in \Phi^*(F)\}$ for each maximal face $F$ of $B_X$, where
\[
\Phi^*(F) = \{ f \in B_{X^*} : f(F) = \{1\} \} .
\]
Then, the set
\[
\mathfrak{S}(X) = \{ I_F : \text{$F$ is a maximal face of $B_X$} \}
\]
equipped with the mapping $S \mapsto S^=$ given by
\[
S^= = \left\{ I \in \s (X) : \bigcap \{ I' : I' \in S\} \subset I \right\}
\]
for each $S \subset \s (X)$ is called the \emph{geometric structure space} of $X$. Let $\mathfrak{C}(X) = \{ S \subset \s (X) : S^= = S \}$.
\end{definition}
See~\cite[Definition 3.4]{Tan22b} for the original definition of $\s (X)$, and \cite[Theorem 4.15]{Tan22b} for the equivalence of two definitions of geometric structure spaces. It was shown in \cite[Proposition 3.5]{Tan22b} that $\s (X)$ is a closure space whenever $X$ is a Banach space. Although the scalar field should be specified in categories whose objects are Banach spaces, in the context of closure spaces, it is not necessary to distinguish the real and complex cases. More precisely, it is appropriate to categorize the scalar-specified versions as a subcategory of (homeomorphs of) geometric structure spaces of \emph{all} (real or complex) Banach spaces.
\begin{definition}

The full subcategory of $\FCl$ whose objects are homeomorphs of geometric structure spaces of nontrivial Banach spaces is denoted by $\Gss$. The full subcategory of $\Gss$ whose objects are homeomorphs of geometric structure spaces of nontrivial Banach spaces over $\mathbb{K}$ is denoted by $\mathbb{K}$-$\Gss$.
\end{definition}
\begin{remark}

From the definition, $\Gss$ is a strictly full subcategory of $\FCl$, and $\mathbb{K}$-$\Gss$ is a strictly full subcategory of $\Gss$.
\end{remark}
\begin{remark}

It was shown in \cite[Theorem 4.11]{Tan22b} that the geometric structure space $\mathfrak{S}(X)$ of a Banach space $X$ is always $T_1$. Namely, each singleton in $\mathfrak{S}(X)$ is closed. Hence, in particular, not $T_1$ topological spaces cannot be objects in $\Gss$.
\end{remark}
The main purpose of this section is to obtain more visible forms of $\g$- and $\p$-transforms of geometric structure spaces of Banach spaces. As in the preceding section, the symbol ``$\sim$'' stands for the equivalence relation on a closure space defined in Definition~\ref{2p-eq}.

Before developing the theory, we give an example which explains why ``$\sim$'' was adopted nevertheless ``$\rightleftarrows$'' is much simpler.
\begin{example}\label{not-transitive}

There exists a Banach space $X$ such that ``$\rightleftarrows$'' is not transitive on $\mathfrak{S}(X)$. Let $X$ be the space $\mathbb{R}^3$ endowed with the norm $\|\cdot\|$ defined by
\[
\|(a,b,c)\| = \|(\|(a,b)\|_\infty ,c)\|_2 = (\max \{ |a|^2,|b|^2 \} + |c|^2)^{1/2}
\]
for each $(a,b,c) \in \mathbb{R}^3$. Then, the dual norm of $\|\cdot\|$ is given by
\[
\|(a,b,c)\|^* = \|(\|(a,b)\|_1,c)\|_2 = ((|a|+|b|)^2+|c|^2)^{1/2}
\]
for each $(a,b,c) \in \mathbb{R}^3$. To be precise, $X^*$ is isometrically isomorphic to $(\mathbb{R}^3,\|\cdot\|^*)$ under the isometric isomorphism
\[
X^* \ni f \mapsto (f(1,0,0),f(0,1,0),f(0,0,1)) \in \mathbb{R}^3.
\]
Let $F_1=\{ (a,b,c) \in B_X : a=1 \}$ and $F_2 = \{ (a,b,c) \in B_X : b=1\}$. Then, $F_1$ and $F_2$ are exposed faces of $B_X$ with respect to the norm-one functionals $f_1=(1,0,0)$ and $f_2=(0,1,0)$. Moreover, if $f=(a,b,c) \in \Phi^* (F_1)$, then it follows from $(1,0,0) \in F_1$ that $a=f_1(1,0,0)=1$. Combining this with $\|(a,b,c)\|^* =1$, we obtain $f=(1,0,0)$. Hence, $\Phi^*(F_1) = \{ f_1 \}$. From this, in particular, $F_1$ is a maximal face of $B_X$. Indeed, if $F$ is a maximal face of $B_X$ containing $F_1$, then $\Phi^*(F) \subset \Phi^*(F_1) = \{ f_1 \}$, which together with the maximality of $F$ implies that $F = f_1^{-1}(1) \cap B_X = F_1$. Similarly, it turns out that $\Phi^* (F_2) = \{ f_2 \}$ and $F_2$ is a maximal faces of $B_X$.

Next, let $F_3=\{ (a,b,c) \in B_X : a=c=1/\sqrt{2}\}$, $F_4=\{ (a,b,c) \in B_X : b=c=1/\sqrt{2}\}$, $f_3 = (1/\sqrt{2},0,1/\sqrt{2})$, and $f_4 = (0,1/\sqrt{2},1/\sqrt{2})$. Then, $\|f_3\|=\|f_4\|=1$. Moreover, by the equality condition for the Cauchy-Schwarz inequality, we have $F_j = f_j^{-1}(\{1\}) \cap B_X$ for $j=3,4$. If $f=(a,b,c) \in \Phi^* (F_3)$, then $(1/\sqrt{2},0,1/\sqrt{2}) \in F_3$ implies that $a/\sqrt{2}+c/\sqrt{2}=1$. Since $(|a|^2+|c|^2)^{1/2} \leq \|(a,b,c)\|^*=1$, again by the equality condition for the Cauchy-Schwarz inequality, we have $a=c=1/\sqrt{2}$ and $b=0$. Therefore, $f=f_3$, that is, $\Phi^*(F_3) = \{f_3\}$. An argument similar to above also shows that $\Phi^*(F_4) = \{f_4\}$. In particular, $F_3$ and $F_4$ are maximal faces of $B_X$.

We note that $e = (0,0,1) \in S_X$ is an extreme point of $B_X$ which is also a point smoothness. Indeed, if $f=(a,b,c) \in S_{X^*}$, and if $f(0,0,1) = 1$, then $c=1$ and $a=b=0$, that is $f = e$. Meanwhile, if $(a,b,c) \in B_X$ satisfies $e(a,b,c)=1$, then $c=1$ and $a=b=0$, that is, $(a,b,c)=e$. Hence, $\{e\}=e^{-1}(\{1\}) \cap B_X$. This means that $\{e\}$ is a one-point face (that is, an extreme point) of $B_X$. In particular, $\{e\}$ is a maximal face of $B_X$ with $\Phi^* (\{e\}) = \{e\}$.

From the above, it follows that $\ker f_j ,\ker e \in \mathfrak{S}(X)$ for $j=1,2,3,4$. Now, we note that $\ker f_1 \cap \ker e \subset \ker f_3$ and $\ker f_2 \cap \ker e \subset \ker f_4$. Hence,
\[
\ker f_3 \in \{ \ker f_1,\ker e\}^= \setminus \{ \ker f_1,\ker e\}
\]
and
\[
\ker f_4 \in \{ \ker f_2,\ker e\}^= \setminus \{ \ker f_2,\ker e\} ,
\]
that is, $\{ \ker f_1,\ker e \} \not \in \mathfrak{C}(X)$ and $\{ \ker f_2,\ker e\} \not \in \mathfrak{C}(X)$. Thus, $\ker f_1 \rightleftarrows \ker e$ and $\ker e \rightleftarrows \ker f_2$. Meanwhile, if $F$ is a maximal face of $B_X$ that contains $x=(a,b,c) \in S_X$ with $c \neq 0$, then $I_F \not \in \{ \ker f_1,\ker f_2\}^=$. Indeed, such an $x$ satisfies $\max \{ |a|^2,|b|^2\} = 1-|c|^2 < 1$, which implies that $c' \neq 0$ whenever $f=(a',b',c') \in \Phi^*(F)$; for otherwise, $|f(x)| = |aa' + bb'| \leq \max \{|a|,|b|\}(|a'|+|b'|) <1$ by $f \in S_{X^*}$. Therefore,
\[
(0,0,1) \in (\ker f_1 \cap \ker f_2) \setminus \bigcup \{ \ker f : f\in \Phi^*(F)\} ,
\]
which proves that $I_F \not \in \{ \ker f_1,\ker f_2\}^=$. Now, let $I_F \in \{ \ker f_1,\ker f_2\}^=$, where $F$ is a maximal face of $B_X$. Then, $F$ is contained in the subspace $M=\{ (a,b,c) \in X : c = 0\}$ which is naturally isometrically isomorphic to $\ell_\infty^2$. Since $F (\subset M)$ is also a maximal face of $B_M$ and $S_M = F_1 \cup (-F_1) \cup F_2 \cup (-F_2)$, it follows that $F \in \{ F_1,-F_1,F_2,-F_2\}$. Hence, $I_F = \ker f_1$ if $F \in \{ F_1,-F_1\}$, while $I_F = \ker f_2$ if $F \in \{ F_2,-F_2\}$. This proves that $\{\ker f_1,\ker f_2 \} \in \mathfrak{C}(X)$, that is, $\ker f_1 \not \rightleftarrows \ker f_2$.
\end{example}
Now, we proceed to develop the theory. First, we give two ways of constructing new closure spaces from the geometric structure space of a given Banach space.
\begin{definition}\label{quotient-def}

Let $X$ be a nontrivial Banach space, and let
\[
J_{C(I)} = \bigcap \{ I' : I' \in C(I) \}
\]
for each $I \in \s (X)$. The \emph{geometric spectrum} $\gs (X)$ and \emph{geometric primitive spectrum} $\ps (X)$ of $X$ are defined as the sets $\gs (X) = \{ C(I) : I \in \s (X)\}$ and
\[
\ps (X) = \{ J_{C(I)} : I \in \s(X) \} .
\]
endowed with the mappings $S \mapsto S^{\equiv}$ on $2^{\gs (X)}$ and $S \mapsto S^{\cong}$ on $2^{\ps (X)}$ given by
\begin{align*}
S^{\equiv} &= \left\{ C(I) \in \gs (X) : \bigcap \{ J_{C(I')}: C(I') \in S \} \subset J_{C(I)} \right\} \\
T^{\cong} &= \left\{ J \in \ps (X) : \bigcap \{ J': J'\in T\} \subset J \right\}
\end{align*}
for each $S \subset \gs (X)$ and each $T \subset \ps (X)$. Let $\mathfrak{GC}(X) = \{ S \subset \gs (X) : S^{\equiv}=S\}$ and $\mathfrak{PC}(X) = \{ T \subset \ps (X) : T^{\cong}=T\}$. Define a mapping $\gamma_X : \gs (X) \to \ps (X)$ by $\gamma_X (C(I)) = J_{C(I)}$ for each $C(I) \in \mathfrak{GS}(X)$.
\end{definition}
It will turn out that geometric spectra and geometric primitive spectra of Banach spaces are closure spaces.
\begin{theorem}\label{quotient}

Let $X$ be a nontrivial Banach space. Then, $\gs (X)$ and $\ps (X)$ are closure spaces. Moreover, $S^{\equiv} = \gamma_X^{-1} (\gamma_X (S)^{\cong})$ for each $S \subset \gs (X)$.
\begin{proof}
Let $S \subset \mathfrak{GS}(X)$ and $C(I) \in \mathfrak{GS}(X)$. Then, we have
\[
\bigcap \{ J' : J' \in \gamma_X (S) \} = \bigcap \{ J_{C(I')} : C(I') \in S\} \subset J_{C(I)} = \gamma_X (C(I))
\]
if and only if $\gamma_X (C(I)) \in \gamma_X (S)^{\cong}$, which implies that
\[
S^{\equiv} = \{ C(I) \in \mathfrak{GS}(X) : \gamma_X (C(I)) \in \gamma_X (S)^{\cong} \} = \gamma_X^{-1}(\gamma_X (S)^{\cong}) .
\]
Hence, in the light of Lemma~\ref{pullback}, it is sufficient to show that $\ps (X)$ is a closure space. First, we note that $\bigcap \{J: J \in \emptyset \} = X$, and $J_{C(I)} \subset I \neq X$ for each $I \in \mathfrak{S}(X)$ since $I \cap F =\emptyset$ provided that $I=I_F$ for some maximal face $F$ of $B_X$. Therefore, $\emptyset^{\cong} = \emptyset$. If $J \in T$, then $\bigcap \{ J' : J' \in T\} \subset J$, which implies that $J \in T^{\cong}$. This shows that $T \subset T^{\cong}$. Since $\bigcap \{ J' : J' \in T\} \subset J$ whenever $J \in T^{\cong}$, it follows that $\bigcap \{ J' : J' \in T\} \subset \bigcap \{ J' : J' \in T^{\cong}\}$. Combining this with (ii), we obtain $\bigcap \{ J' : J' \in T\} = \bigcap \{ J' : J' \in T^{\cong}\}$, that is, $(T^{\cong})^{\cong}=T^{\cong}$. Finally, if $T_1 \subset T_2 \subset \mathfrak{PS}(X)$, then $\bigcap \{ J' : J' \in T_2\} \subset \bigcap \{ J' : J' \in T_1\}$. Therefore, $T_1^{\cong} \subset T_2^{\cong}$. This completes the proof.
\end{proof}
\end{theorem}
\begin{remark}\label{coincide}

Let $X$ be a Banach space such that ``$\sim$'' and ``$=$'' are equivalent on $\s (X)$. Then, $C(I) = \{I\}$ and $J_{C(I)} = I$ for each $I \in \s (X)$. Hence, $\s (X)$ coincides with $\ps (X)$, and $\gamma_X : \gs (X) \to \ps (X)$ is bijective. In particular, $\gs (X)$ is homeomorphic to $\mathfrak{S}(X)$ by Lemma~\ref{pullback} (iv).
\end{remark}
We give some examples of geometric spectra and geometric primitive spectra of Banach spaces.
\begin{example}\label{topologizable-ex}

Let $X$ be a Banach space. Suppose that $\s (X)$ is topologizable. Then, each finite subset of $\mathfrak{S}(X)$ is closed, which together with Remark~\ref{equal} implies that ``$\sim$'' and ``$=$'' are equivalent on $\s (X)$. Therefore, $\s (X),\gs (X)$ and $\ps (X)$ are homeomorphic to each other by Remark~\ref{coincide}.
\end{example}
\begin{example}\label{refsm-ex}

Let $X$ be a reflexive smooth Banach space with $\dim X \geq 2$. As in the proof of \cite[Theorem 3.8]{Tan22b}, each two-point set in $\mathfrak{S}(X)$ is not closed. Hence, we have $I_1 \sim I_2$ for any choice of $I_1,I_2 \in \s (X)$, that is, $\gs (X) = \{ \s (X)\}$. Moreover since $J_{\s (X)} = \bigcap \{ I : I \in \s (X)\} = \{0\}$ by \cite[Lemma 4.1]{Tan23a}, we obtain $\ps (X) = \{ \{0\}\}$. Consequently, $\gs (X)$ and $\ps (X)$ are topologizable, and they are both homeomorphic to a one-point topological space endowed with the trivial topology.
\end{example}
\begin{example}\label{1-ex}

The geometric spectrum and geometric primitive spectrum of $\ell_1$ depend on the scalar field $\mathbb{K}$. Recall that $\s (\ell_1) = \{ \ker f_u : u \in \mathcal{U} \}$ by \cite[Theorem 5.8]{Tan22b}, where
\[
\mathcal{U} = \{ (u_n)_n \in \ell_\infty : \text{$|u_n|=1$ for each $n \in \mathbb{N}$}\}
\]
and $f_u \in \ell_1^*$ is defined by $f_u ((a_n)_n) = \sum_n a_n u_n$ for each $u=(u_n)_n \in \mathcal{U}$ and each $(a_n)_n \in \ell_1$.

(I) Let $\mathbb{K}=\mathbb{R}$. Suppose that $u=(u_n)_n$ and $v=(v_n)_n$ are elements of $\mathcal{U}$ such that $\ker f_u \neq \ker f_v$. Replacing $u$ (or $v$) with $-u$ (or $-v$) if necessary, we may assume that $u_1=v_1=1$. Since $\ker f_u \neq \ker f_v$, there exists an $n \in \mathbb{N}$ such that $u_n \neq v_n$. Let $\ker f_w \in \{ \ker f_u,\ker f_v\}^=$. It follows from $\ker f_u \cap \ker f_v \subset \ker f_w$ that $f_w = \alpha f_u + \beta f_v$ for some $\alpha ,\beta \in \mathbb{R}$, in which case, $w = \alpha u + \beta v$. In particular, we have $|\alpha + \beta |= |w_1| = 1$ and $|\alpha -\beta |=|w_n|=1$. Since $\alpha ,\beta \in \mathbb{R}$, it turns out that $\alpha \beta =0$. Therefore, $f_w \in \{ f_u,-f_u,f_v,-f_v\}$, that is, $\ker f_w \in \{ \ker f_u,\ker f_v\}$. This shows that $\{ \ker f_u,\ker f_v\} \in \mathfrak{C}(\ell_1)$, which together with Remark~\ref{equal} implies that ``$\sim$'' and ``$=$'' are equivalent on $\mathfrak{S}(\ell_1)$. Thus, in the real case, $\mathfrak{S}(\ell_1)$ coincides with $\mathfrak{PS}(\ell_1)$, and $\mathfrak{GS}(\ell_1)$ is homeomorphic to $\mathfrak{S}(\ell_1)$. It is worth mentioning here that $\mathfrak{S}(\ell_1)$ is not topologizable; see \cite[Theorem 5.8]{Tan22b}.

(II) Let $\mathbb{K}=\mathbb{C}$. Suppose that $u=(u_n)_n$ and $v=(v_n)_n$ are elements of $\mathcal{U}$. We claim that $\ker f_u \rightleftarrows \ker f_v$ if and only if the set $\{ u_n\overline{v_n} : n \in \mathbb{N}\}$ contains at most two elements. To show this, suppose that $\{ u_n\overline{v_n} : n \in \mathbb{N}\}$ contains at most two elements. We note that if $\{ u_n\overline{v_n} : n \in \mathbb{N}\}$ is a singleton, then $v$ is a scalar multiple of $u$, in which case, $\ker u = \ker v$. Next, let $\{ u_n\overline{v_n} : n \in \mathbb{N}\} = \{ e^{i\theta_1 },e^{i\theta_2}\}$. Replacing $u$ with $e^{i\theta_1}u$, we may assume that $\{ u_n\overline{v_n} : n \in \mathbb{N}\} = \{ 1,e^{i\theta}\}$ and $1 \neq e^{i\theta}$. Set $N_1=\{ n \in \mathbb{N} : u_n = v_n \}$ and $N_2 = \mathbb{N} \setminus N_1$. If $(a_n)_n \in \ker f_u \cap \ker f_v$, then $0 = f_u ((a_n)_n)= \sum_n a_nu_n$ and
\[
0 = f_v ((a_n)_n) = \sum_n a_nv_n = \sum_{n \in N_1} a_nu_n + e^{-i\theta}\sum_{n \in N_2} a_nu_n ,
\]
which implies that $\sum_{n \in N_1}a_nu_n = \sum_{n \in N_2}a_nu_n = \sum_{n \in N_2}a_nv_n= 0$. Now, let
\[
w_n = \left\{ \begin{array}{ll}
u_n & (n \in N_1) \\
e^{i\theta /2}v_n & (n \in N_2)
\end{array}
\right. .
\]
Then, $w=(w_n)_n \in \mathcal{U}$ and
\[
\ker f_w \in \{ \ker f_u,\ker f_v \}^= \setminus \{ \ker f_u,\ker f_v\} .
\]
Hence, $\{ \ker f_u,\ker f_v\} \not \in \mathfrak{C}(\ell_1)$ and $\ker f_u \rightleftarrows \ker f_v$.

Conversely, suppose that $\{ u_n\overline{v_n} : n \in \mathbb{N}\}$ contains three elements. In this case, we have $\ker f_u \neq \ker f_v$; for otherwise, $\{u,v\}$ is linearly dependent, and $\{ u_n\overline{v_n} : n \in \mathbb{N}\}$ becomes a singleton. As in the preceding paragraph, it may be assumed that $\{ u_n\overline{v_n} : n \in \mathbb{N}\} = \{ 1,e^{i\theta_1},e^{i\theta_2}\}$ with $0<\theta_1<\theta_2<2\pi$. Let $w \in \mathcal{U}$ be such that $\ker f_w \in \{ \ker f_u,\ker f_v\}^=$. Then, $w = \alpha u+\beta v$ for some $\alpha ,\beta \in \mathbb{C}$. Choose $n,n_1,n_2 \in \mathbb{N}$ such that $u_n=v_n$, and $u_{n_j}=e^{i\theta_j}v_{n_j}$ for $j=1,2$. It follows that $|\alpha + \beta| = |w_n|=1$ and $|e^{i\theta_j}\alpha +\beta|=|w_{n_j}| =1$ for $j=1,2$, which implies that $\re \alpha \overline{\beta} = \re e^{i\theta_1}\alpha \overline{\beta} = \re e^{i\theta_2}\alpha \overline{\beta}$. If $\alpha \overline{\beta} = re^{i\theta} \neq 0$, then we obtain $\theta <\theta +\theta_1<\theta +\theta_2<\theta +2\pi$. Meanwhile, $\re e^{i\varphi_1} = \re e^{i\varphi_2}$ implies that $e^{i\varphi_1} = \overline{e^{i\varphi_2}}$. Therefore, $\alpha \overline{\beta}=0$. This proves that $\{ \ker f_u,\ker f_v\} \in \mathfrak{C}(X)$. Combining this with $\ker f_u \neq \ker f_v$, we derive that $\ker f_u \not \rightleftarrows \ker f_v$.

Next, we claim that $\ker f_u \sim \ker f_v$ if and only if $\{ u_n\overline{v_n} : n \in \mathbb{N}\}$ is a finite set. Indeed, if $\ker f_u \sim \ker f_v$, then there exists a finite subset $\{w^{(1)},\ldots ,w^{(m)}\} \subset \mathcal{U}$ such that $\ker f_u \rightleftarrows \ker f_{w^{(1)}} \rightleftarrows \cdots \rightleftarrows \ker f_{w^{(m)}} \rightleftarrows \ker v$. By what we have proved in above, it follows that each of
\begin{align*}
&\left\{ u_n\overline{w^{(1)}_n} : n \in \mathbb{N} \right\}, \\
&\left\{ w^{(1)}_n\overline{w^{(2)}_n} : n \in \mathbb{N} \right\},\ldots ,\left\{ w^{(m-1)}_n\overline{w^{(m)}_n} : n \in \mathbb{N} \right\} ,\\
&\left\{ w^{(m)}_n\overline{v_n} : n \in \mathbb{N} \right\}
\end{align*}
contains at most two elements. Moreover, we have
\[
u_n\overline{v_n} = \left(u_n\overline{w^{(1)}_n}\right)\left(w^{(1)}_n\overline{w^{(2)}_n}\right) \cdots \left(w^{(m-1)}_n\overline{w^{(m)}_n}\right)\left(w^{(m)}_n\overline{v_n}\right)
\]
for each $n \in \mathbb{N}$. Thus, $\{ u_n\overline{v_n} : n \in \mathbb{N}\}$ contains at most $2^{m+1}$ elements.

For the converse, let $\{ u_n\overline{v_n} : n \in \mathbb{N}\} = \{ e^{i\theta_1},\ldots ,e^{i\theta_m}\}$ with $0 \leq \theta_1<\cdots <\theta_m <2\pi$, and let $N_j = \{ n \in \mathbb{N} : u_n\overline{v_n} = e^{i\theta_j} \}$. Set
\[
w^{(j)}_n = \left\{ \begin{array}{ll}
u_n & (n \not \in N_1 \cup \cdots \cup N_j) \\
v_n & (n \in N_1 \cup \cdots \cup N_j)
\end{array}
\right.
\]
for each $j \in \{ 1,\ldots ,m\}$ and each $n \in \mathbb{N}$. Since $N_1 \cup \cdots \cup N_m = \mathbb{N}$, we have $w^{(m)} = v$. Moreover, it follows that
\[
\left\{ u_n\overline{w^{(1)}_n} : n \in \mathbb{N} \right\} = \{ 1,e^{i\theta_1}\}
\]
and
\[
\left\{ w^{(j)}_n\overline{w^{(j+1)}_n} : n \in \mathbb{N} \right\} = \{ 1,e^{i\theta_{j+1}}\}
\]
for each $j \in \{1,\ldots ,m-1\}$. This means that
\[
\ker f_u \rightleftarrows \ker f_{w^{(1)}} \rightleftarrows \cdots \rightleftarrows \ker f_{w^{(m)}} = \ker v,
\]
that is, $\ker f_u \sim \ker f_v$.

Now, let $u \in \mathcal{U}$. Then, we know that
\[
C(\ker f_u) = \{ \ker f_v : \text{$v \in \mathcal{U}$ and $\{ u_n\overline{v_n} : n \in \mathbb{N}\}$ is a finite set}\} .
\]
Suppose that $(a_n)_n \in J_{C(\ker f_u)}$. For each $m \in \mathbb{N}$, let
\[
v^{(m)}_n = \left\{ \begin{array}{ll}
u_n & (n \neq m) \\
-u_n & (n=m)
\end{array}
\right. .
\]
Since
\[
\left\{ u_n\overline{v^{(m)}_n} : n \in \mathbb{N} \right\} = \{ -1,1\},
\]
it follows that $\ker u \sim \ker v^{(m)}$. From this, we obtain $0=f_u ((a_n)_n) = \sum_n a_n u_n$ and
\[
0 = f_{v^{(m)}} ((a_n)_n) = \sum_n a_nv_n = -a_mu_m + \sum_{n \neq m}a_nu_n ,
\]
which together with $|u_m|=1$ implies that $a_m=0$. Thus, $(a_n)_n=0$. This shows that $J_{C(\ker u)} = \{0\}$ for each $\ker u \in \mathcal{U}$, that is, $\ps (\ell_1) = \{ \{0\}\}$ and $\mathfrak{GC}(\ell_1) = \{ \emptyset ,\gs (\ell_1)\}$ in the complex case. In particular, if $\mathbb{K}=\mathbb{C}$, then $\ps (\ell_1)$ and $\gs (\ell_1)$ are both topologizable and have the trivial topology.
\end{example}
\begin{remark}

The complex $\ell_1$-space gives another example of a (complex) Banach space $X$ such that ``$\rightleftarrows$'' is not transitive on $\mathfrak{S}(X)$.
\end{remark}
Next, we clarify the positions of geometric (primitive) spectra of Banach spaces as closure spaces. For this purpose, three auxiliary results will be needed.
\begin{lemma}\label{closed}

Let $X$ be a Banach space, and let $I \in \mathfrak{S}(X)$. Then,
\[
C(I)^= = \{ I' \in \s (X) : J_{C(I)} \subset I' \}
\]
and $J_{C(I)} = \bigcap \{ I' : I' \in C(I)^= \}$. In particular, $C(I_1)^= = C(I_2)^=$ if and only if $J_{C(I_1)}=J_{C(I_2)}$ for $I_1,I_2 \in \s (X)$.
\begin{proof}
From the definitions of $C(I)^=$ and $J_{C(I)}$, we have
\[
C(I)^= = \{ I' \in \s (X) : J_{C(I)} \subset I' \} ,
\]
which implies that $J_{C(I)} \subset \bigcap \{ I' : I' \in C(I)^=\}$. Meanwhile, since $C(I) \subset C(I)^=$, it follows that
\[
J_{C(I)} = \bigcap \{ I' : I' \in C(I) \} \supset \bigcap \{ I': I' \in C(I)^=\} .
\]
Therefore, $J_{C(I)} = \bigcap \{ I' : I' \in C(I)^= \}$.
\end{proof}
\end{lemma}
\begin{lemma}\label{lifting}

Let $X$ be a Banach space. Then,
\[
\bigcap \{ J : J \in T\} = \bigcap \{ I \in \s (X) : J_{C(I)} \in T\}
\]
for each $T \subset \ps (X)$.
\begin{proof}
Suppose that $x \in \bigcap \{J :J \in T\}$, and that $J_{C(I)} \in T$. Then, $x \in J_{C(I)} \subset I$, which shows that
\[
\bigcap \{ J : J \in T\} \subset \bigcap \{ I \in \s (X) : J_{C(I)} \in T\} .
\]
Conversely, suppose that $x \in \bigcap \{ I \in \s (X) : J_{C(I)} \in T\}$. If $J \in T$, then $J=J_{C(I)}$ for some $I \in \s (X)$, in which case, $J_{C(I')} = J_{C(I)} \in T$ whenever $I' \in C(I)$. Hence, $x \in \bigcap \{ I' : I' \in C(I)\} = J_{C(I)} =J$, which implies that $x \in \bigcap \{J : J \in T\}$. Therefore,
\[
\bigcap \{ J : J \in T\} \supset \bigcap \{ I \in \mathfrak{S}(X) : J_{C(I)} \in T\} .
\]
This completes the proof.
\end{proof}
\end{lemma}
\begin{lemma}\label{closure-point}

Let $X$ be a Banach space, let $I \in \mathfrak{S}(X)$, and let $T \subset \mathfrak{PS}(X)$. Then, 
$J_{C(I)} \in T^{\cong}$ if and only if $C(I)^= \subset \{ I' \in \mathfrak{S}(X) : J_{C(I')} \in T\}^=$.
\begin{proof}
Suppose that $J_{C(I)} \in T^{\cong}$. If $I_0 \in C(I)$. then we have
\[
\bigcap \{ I' \in \s (X) : J_{C(I')} \in T\}=\bigcap \{ J : J \in T\} \subset J_{C(I)} \subset I_0
\]
by Lemma~\ref{lifting}, that is, $I_0 \in \{ I' \in \s (X) : J_{C(I')} \in T\}^=$. Hence,
\[
C(I)^= \subset \{ I' \in \s (X) : J_{C(I')} \in T\}^=.
\]

For the converse, we assume that $C(I)^= \subset \{ I' \in \s (X) : J_{C(I')} \in T\}^=$. If $I_0 \in C(I)$, then $I_0 \in \{ I' \in \s (X) : J_{C(I')} \in T\}^=$, which together with Lemma~\ref{lifting} implies that
\[
\bigcap \{ J : J \in T\} = \bigcap \{ I' \in \s (X) : J_{C(I')} \in T\} \subset I_0 .
\]
Thus,
\[
\bigcap \{ J : J \in T\} \subset \bigcap \{ I_0 : I_0 \in C(I) \} = J_{C(I)} ,
\]
that is, $J_{C(I)} \in T^{\cong}$.
\end{proof}
\end{lemma}
The following is the main result in this section.

\begin{theorem}\label{spec-trans}

Let $X$ be a Banach space. Then, $\gs (X)$ coincides with the $\g$-transform of $\s (X)$ and the mapping $C(I)^= \mapsto J_{C(I)}$ is a homeomorphism from the $\p$-transform of $\s (X)$ onto $\ps (X)$.
\begin{proof}
Let $(\s (X)^\g,c^\g)$ and $(\s (X)^\p,c^\p)$ be the $\g$- and $\p$-transforms of $\s (X)$, respectively, and let $\delta (C(I)) = C(I)^=$ for each $I \in \s (X)$. Then, $\s (X)^\g = \{ C(I) : I \in \s (X)\} = \gs (X)$, $\s (X)^\p = \{ C(I)^= : I \in \s (X) \}$,
\[
c^\p (T) = \{ C(I)^= : C(I)^= \subset \{ I' \in \s (X) : C(I')^= \in T\}^= \}
\]
for each $T \subset \s (X)^\p$, and $c^\g (S) = \delta^{-1}(c^\p (\delta (S)))$ for each $S\subset \s (X)^\g$. Define a mapping $\Phi : \s (X)^\p \to \ps (X)$ by $\Phi (C(I)^=) = J_{C(I)}$ for each $I \in \s (X)$. Then, $\Phi$ is well-defined and bijective by Lemma~\ref{closed}.

Let $T\subset \s (X)^\p$ and $I \in \s (X)$. Then, $C(I)^= \in c^\p (T)$ if and only if
\[
C(I)^= \subset \{ I' \in \s (X) : C(I')^= \in T\}^= = \{ I' \in \s (X) : J_{C(I')} \in \Phi (T)\}^= ,
\]
which occurs if and only if $\Phi (C(I)^=) = J_{C(I)} \in \Phi (T)^{\cong}$ by Lemma~\ref{closure-point}. Therefore, $\Phi (c^\p (T)) = \Phi (T)^{\cong}$. This shows that $\Phi$ is a homeomorphism.

Finally, let $S \subset \s (X)^\g = \gs (X)$. Since $\Phi \circ \delta = \gamma_X$, it follows from the preceding paragraph that
\begin{align*}
S^{\equiv} = \gamma_X^{-1}(\gamma_X (S)^{\cong})
&= (\Phi \circ \delta)^{-1}((\Phi \circ \delta)(S)^{\cong}) \\
&= (\Phi \circ \delta)^{-1}(\Phi (c^\p (\delta (S)))) \\
&= \delta^{-1}(\Phi^{-1}(\Phi (c^\p (\delta (S))))) \\
&= \delta^{-1}(c^\p (\delta (S))) \\
&= c^\g (S) .
\end{align*}
Thus, the closure operators on $\gs (X)$ and $\s (X)^\g$ coincide.
\end{proof}
\end{theorem}
Combining the preceding theorem with Theorem~\ref{two-homeo}, we immediately have the following result.
\begin{corollary}\label{s-ps}

Let $X$ and $Y$ be Banach spaces, and let $\Phi :\s (X) \to \s (Y)$ be a homeomorphism. Define $\Phi^{\gs}(C(I))= C(\Phi (I))$ and $\Phi^{\ps} (J_{C(I)}) = J_{C(\Phi (I))}$ for each $I \in \s (X)$. Then, $\Phi^{\gs} :\mathfrak{GS}(X) \to \mathfrak{GS}(Y)$ and $\Phi^{\ps} : \mathfrak{PS}(X) \to \mathfrak{PS}(Y)$ are homeomorphisms. Moreover, $\gamma_Y \circ \Phi^{\gs} = \Phi^{\ps} \circ \gamma_X$ on $\mathfrak{GS}(X)$. In particular, $\gamma_X :\mathfrak{GS}(X) \to \mathfrak{PS}(X)$ is injective if and only if $\gamma_Y :\mathfrak{GS}(Y) \to \mathfrak{PS}(Y)$ is.
\end{corollary}
\begin{remark}\label{converse}

Let $X$ and $Y$ be Banach spaces. Then, by the preceding theorem, $\s (X) \simeq \s (Y)$ implies that $\gs (X) \simeq \gs (Y)$ and $\ps (X)\simeq \ps (Y)$. We remark that the converse to this is false in general since all reflexive smooth Banach spaces $X$ with $\dim X \geq 2$ have the same geometric (primitive) spectra. In particular, it should be emphasized that $\mathfrak{GS}(X)$ and $\mathfrak{PS}(X)$ no longer remember the dimension of $X$ in general.
\end{remark}
We conclude this section with a useful tool for obtaining homeomorphisms between geometric structure spaces from isometric isomorphisms between underlying Banach spaces. Recall that an element $x$ of a Banach space over $\mathbb{K}$ is said to be \emph{Birkhoff-James orthogonal} to another $y$, denoted by $x \perp_{BJ}y$, if $\|x+\lambda y\| \geq \|x\|$ for each $\lambda \in \mathbb{K}$. This is known as a generalization of the usual orthogonality in Hilbert spaces, and was named after Birkhoff and James since Birkhoff~\cite{Bir35} first introduced it and James~\cite{Jam47a,Jam47b} greatly developed its theory. For nice surveys on Birkhoff-James orthogonality, see \cite{AMW12,AMW22,AGKZ22}.

Let $X$ and $Y$ be Banach spaces, and let $T:X \to Y$. Then, $T$ is called a \emph{Birkhoff-James orthogonality preserver} if it is bijective and $x \perp_{BJ} y$ is equivalent to $Tx \perp_{BJ}Ty$. The following result comes from \cite[Theorem 3.10]{Tan22b}, which is useful for constructing homeomorphisms between geometric structure spaces of Banach spaces.
\begin{theorem}\label{Tan-310}

Let $X$ and $Y$ be Banach spaces, and let $T:X \to Y$ be a Birkhoff-James orthogonality preserver. Define a mapping $\Phi_T :\s (X) \to \s (Y)$ by $\Phi_T(I) = T(I)$ for each $I \in \s (X)$. Then, $\Phi_T$ is a homeomorphism with $\Phi_T^{-1} = \Phi_{T^{-1}}$.
\end{theorem}
Since isometric isomorphism between Banach spaces are obviously Birkhoff-James orthogonality preservers, we can apply the preceding theorem to such mappings.

\section{The case of $C^*$-algebras}\label{Sect:CAlg}

This section is devoted to the case of $C^*$-algebras. Let $\A$ be a $C^*$-algebra. Then, as was shown in~\cite[Proposition 3.2]{Tan23c}, the geometric structure space of $\A$ is given by
\[
\s (\A) = \{ \ker \rho : \rho \in \ext (B_{\A^*})\} .
\]
Based on this information, we aim to identify the $\g$- and $\p$-transforms of $\s (\A)$ via $\gs (\A)$ and $\ps (\A)$.
\begin{remark}\label{polar-u}

Let $\A$ be a $C^*$-algebra, and let $\rho \in \ext (B_{\A^*})$. Then, as was noted (without a detailed proof) in \cite[5.4]{AP92}, there exist a pure state $\tau$ of $\A$ and a unitary element $U$ of $\uA$ such that $\rho U = \tau$, where $\rho U \in \A^*$ is defined by $(\rho U)(A) =\rho (UA)$ for each $A \in \A$. We realize this fact for the reader's convenience.

Let $\rho \in \ext (B_{\A^*})$, and let $F = \{ \tau \in B_{\uA^*} : \tau |\A = \rho\}$. Then, $F \neq \emptyset$ by the Hahn-Banach theorem, and $F$ is a weakly$^*$ closed face of $B_{A^*}$. Hence, by the Krein-Milman theorem, there exists a $\overline{\rho} \in \ext (F) \subset \ext (B_{\uA^*})$ such that $\overline{\rho}|\A = \rho$. If there exist a pure state $\tau$ of $\uA$ and a unitary $U \in \uA$ such that $\overline{\rho}U=\tau$, then $\rho U = \tau |\A$ and $\rho = (\tau |\A)U^*$ since $\A$ is an ideal of $\uA$, in which case, $\rho |\A$ is a pure state of $\A$. Indeed, we have $\|\tau |\A\| = \|\rho \|=1$ by
\[
\|\tau |\A\| = \|\rho U\| \leq \|\rho \| = \|(\tau |\A)U^*\| \leq \|\tau |\A\| .
\]
If $\tau |\A = (1-t)\tau_1 +t\tau_2$ for some state $\tau_1,\tau_2$ of $\A$ and some $t \in (0,1)$, then $\tau_j$ is extended uniquely to a state $\overline{\tau_j}$ of $\uA$ for $j=1,2$. To be precise, let $\overline{\tau_j}(A+c\mathbf{1}) = \tau_j (A)+c$ for each $A+c\mathbf{1} \in \uA$ and $j=1,2$. It follows that
\begin{align*}
((1-t)\overline{\tau_1}+t\overline{\tau_2})(A+c\mathbf{1}) 
&= (1-t)(\tau_1 (A)+c) + t(\tau_2 (A)+c) \\
&= ((1-t)\tau_1+t\tau_2)(A)+c \\
&= \tau (A)+ \tau (c\mathbf{1}) \\
&= \tau (A+c\mathbf{1})
\end{align*}
for each $A+c\mathbf{1} \in \uA$, that is, $\tau = (1-t)\overline{\tau_1}+t\overline{\tau_2}$. Since $\tau$ is a pure state of $\uA$, we have $\tau = \overline{\tau_1}=\overline{\tau_2}$ and $\tau |A = \tau_1 =\tau_2$. This shows that $\tau |\A$ is a pure state of $\A$. Hence, it is sufficient to prove the claim for unital $C^*$-algebras.

Suppose that $\A$ is unital, and that $\rho \in \ext (B_{\A^*})$. Let $\pi_u : \A \to \HH$ be the universal representation of $\A$, and let $\mathcal{R}$ be the weak-operator closure of $\pi_u (\A)$ in $B(\HH)$. We note that each element of $\pi_u (\A)^*$ is normal and extends, uniquely and without changing norm, to a normal element of $\mathcal{R}^*$. Let $\rho_0 \in \mathcal{R}^*$ denote the normal extension of $\rho \circ \pi_u^{-1} \in \ext (B_{\pi_u (\A)^*})$. By the polar decomposition, there exist a normal state $\tau_0$ of $\mathcal{R}$ and a partial isometry $V\in \mathcal{R}$ such that $\rho_0 V=\tau_0$ and $\rho_0 =\tau_0 V^*$. We note that $\tau_0 |\pi_u (\A)$ is a pure state of $\pi_u(\A)$. Indeed, we have $\|\tau_0 |\pi_u (\A)\| = \|\tau_0 \|= 1$ by the Kaplansky density theorem since $\tau$ is normal. Moreover, if $\tau_0 = (1-t)\tau_1 + t\tau_2$ on $\pi_u (\A)$ for some states $\tau_1,\tau_2$ of $\pi_u (\A)$ and some $t \in (0,1)$, then $\tau_0 = (1-t)\tau_1 + t\tau_2$ on $\mathcal{R}$ by the weak-operator density of $\pi_u (\A)$ in $\mathcal{R}$, where $\tau_j$ is extended uniquely to a normal linear functional for $j=1,2$. It follows from $\tau_0 (V^*V) = \rho_0 (V) = \tau_0 (\mathbf{1}) = 1$ that $\tau_1 (V^*V) = \tau_2 (V^*V) = 1$, which implies that $\tau_j = \tau_j V^*V$ for $j=1,2$. Now, we note that $\rho_0 = \tau_1 V^* = \tau_2 V^*$ since $\rho_0 = \tau_0 V^* = (1-t)(\tau_1 V^*) + t(\tau_2 V^*)$ and $\rho \in \ext (B_{\pi_u(\A)^*})$. Therefore, $\tau_0 = \rho V = \tau_1 = \tau_2$ by $\tau_j V^*V=\tau_j$ for $j=1,2$. This shows that $\tau_0 |\pi_u (\A)$ is a pure state of $\pi_u (\A)$, and hence, $\tau = \tau_0 \circ \pi_u$ is a pure state of $\A$. Since $\pi_u$ is a direct sum of all representations obtained from states of $\A$ by the GNS construction, the representation $(\pi_\tau ,\HH_\tau )$ is a summand of $(\pi_u,\HH )$. In particular, we can regard $\HH_\tau$ as a subspace of $\HH$ in a natural manner. Let $x \in S_{\HH_\tau}$ be such that $\tau (A) = \langle \pi_\tau (A)x,x \rangle$ for each $A \in \A$, and let $P$ be the projection from $\HH$ onto $\HH_\tau$. Then, we obtain $Px=x$ and
\[
\tau_0 (\pi_u (A)) = \tau (A) = \langle \pi_\tau (A)x,x \rangle = \langle \pi_u (A)x,x \rangle 
\]
for each $A \in \A$. Combining this with the weak-operator density of $\pi_u (A)$ in $\mathcal{R}$, we derive that $\tau_0 (A) = \langle Ax,x \rangle$ for each $A \in \mathcal{R}$. Moreover since $P \in \pi_u (A)'$, it follows that
\begin{align*}
\rho (A) = \rho_0 (\pi_u (A)) = \tau_0 (V^*\pi_u (A)) 
&= \langle V^*\pi_u (A)x,x \rangle \\
&= \langle \pi_u (A)x,PVPx \rangle \\
&= \langle \pi_\tau (A)x,PVPx \rangle
\end{align*}
for each $A \in \A$. We note that $\|PVPx\|=1$ by $\|\rho\|=1$. Since $\tau$ is a pure state, the representation $\pi_\tau$ is irreducible. Thus, there exists a unitary element $U \in \A$ such that $\pi_\tau (U)PVPx = x$, in which case,
\[
\rho (U^*A) = \langle \pi_\tau (U^*A)x,PVPx \rangle = \langle \pi_\tau (A)x,\pi_\tau (U)PVPx \rangle = \tau (A)
\]
for each $A \in \A$, that is, $\rho U^* = \tau$. This proves the claim.
\end{remark}
\begin{remark}\label{polar-uu}

The converse to the claim in the preceding remark is also true. Namely, if $\A$ is a $C^*$-algebra, then $\rho U \in \ext (B_{\A^*})$ whenever $\rho$ is a pure state of $\A$ and $U$ is a unitary element of $\uA$. Indeed since $\A$ is an ideal of $\uA$, the mapping $\tau \mapsto \tau U$ is an isometric isomorphism on $\A^*$. Combining this with the fact that $\rho \in \ext (B_{\A^*})$, we obtain $\rho U \in \ext(B_{\A^*})$. In summary, the identity
\[
\ext (B_{\A^*}) = \{ \rho U : \text{$\rho$ is a pure state of $\A$ and $U \in \mathcal{U}$}\}
\]
holds, where $\mathcal{U}$ is the unitary group of $\uA$.
\end{remark}
The following lemma readily follows from Remark~\ref{polar-uu}.

\begin{lemma}\label{gss-alg}

Let $\A$ be a $C^*$-algebra. Then,
\[
\mathfrak{S}(\A) = \{ \ker \rho U : \text{$\rho$ is a pure state of $\A$ and $U \in \mathcal{U}$} \} ,
\]
where $\mathcal{U}$ is the unitary group of $\uA$.
\end{lemma}
To identify $\gs (\A)$ and $\ps (\A)$ for $C^*$-algebras $\A$, we give a characterization of the equivalence relation on $\s (\A)$ in terms of representations. To this end, some auxiliary results will be needed. Let $\A$ be a $C^*$-algebra, and let $\rho$ be a state of $\A$. Then, $\pi_\rho : \A \to B(\HH_\rho)$ denotes the $*$-representation of $\A$ obtained from $\rho$ by the GNS construction. If $x_\rho \in S_{\HH_\rho}$ satisfies $\rho (A) = \langle \pi_\rho (A)x_\rho,x_\rho \rangle$ for each $A \in \A$, then
\[
\langle (\pi_\rho (A)+c\mathbf{1})x_\rho,x_\rho \rangle = \langle \pi_\rho (A)x_\rho,x_\rho \rangle + c = \rho (A)+c = \overline{\rho}(A+c\mathbf{1})
\]
for each $A \in \A$ and each $c \in \mathbb{C}$, where $\overline{\rho}$ is a unique extension of $\rho$ to $\uA$. Hence, we can extend $\pi_\rho$ to $\uA$ by setting $\pi_\rho (A+c\mathbf{1}) = \pi_\rho (A)+c\mathbf{1}$ for each $A+c\mathbf{1} \in \uA$, in which case, $\overline{\rho}(A+c\mathbf{1}) = \langle \pi_\rho (A+c\mathbf{1})x_\rho ,x_\rho \rangle$ for each $A+c\mathbf{1} \in \uA$. In the following, if $\A$ is non-unital, the symbol $\pi_\rho$ stands for this extended representation of $\uA$. We note that $\pi_\rho$ is an irreducible representation of $\uA$ whenever $\rho$ is a pure state of $\A$ since $\pi_\rho (\A) \subset \pi_\rho (\uA)$.
\begin{lemma}\label{ext-xy}

Let $\A$ be a $C^*$-algebra, and let $\rho$ be a pure state of $\A$. Define an element $\rho_{x,y}$ of $\A^*$ by
\[
\rho_{x,y} (A) = \langle \pi_\rho (A)x,y \rangle
\]
for each $A \in \A$ and each $x,y \in \HH_\rho$. Then, $\rho_{x,y} \in \ext (B_{\A^*})$ for each $x,y \in S_{\HH_\rho}$.
\begin{proof}
Let $x,y \in S_{\HH_\rho}$, and let $x_\rho \in S_{\HH_\rho}$ be such that $\rho (A) = \langle \pi_\rho (A)x_\rho ,x_\rho \rangle$ for each $A \in \A$. Since the representation $\pi_\rho$ of $\uA$ is irreducible, there exist unitary elements $U,V$ of $\uA$ such that $\pi (U)x_\rho = x$ and $\pi (V)x_\rho = y$. It follows that
\[
(U\rho V^*)(A) = \rho (V^*AU) = \langle \pi_\rho (V^*AU)x_\rho ,x_\rho \rangle = \langle \pi_\rho (A)x,y\rangle =\rho_{x,y}(A)
\]
for each $A \in \A$, that is, $U\rho V^* =\rho_{x,y}$. Combining this with the fact that the mapping $\tau \mapsto U\tau V^*$ is an isometric isomorphism on $\A^*$, we infer that $\rho_{x,y} \in \ext (B_{\A^*})$.
\end{proof}
\end{lemma}
Let $\rho$ be a pure state of a $C^*$-algebra $\A$. Then, it follows from
\[
\rho_{x,y}(A) = \langle \pi_\rho (A)x,y \rangle = \overline{\langle \pi_\rho (A^*)y,x \rangle} = \overline{\rho_{y,x}(A^*)}
\]
for each $A \in \A$ that $\rho_{x,y}(A)=0$ if and only if $\rho_{y,x}(A^*)=0$.
\begin{lemma}\label{ker-equal}

Let $\A$ be a $C^*$-algebra, and let $\rho$ be a pure state of $\A$. Then, the following hold:
\begin{itemize}
\item[{\rm (i)}] Suppose that $x,y_1,y_2 \in \HH_\rho \setminus \{0\}$. Then, $\ker \rho_{x,y_1} = \ker \rho_{x,y_2}$ if and only if $\{y_1,y_2\}$ is linearly dependent.
\item[{\rm (ii)}] Suppose that $x_1,x_2,y \in \HH_\rho \setminus \{0\}$. Then, $\ker \rho_{x_1,y} = \ker \rho_{x_2,y}$ if and only if $\{x_1,x_2\}$ is linearly dependent.
\end{itemize}
\begin{proof}
(i) Suppose that $\ker \rho_{x,y_1} = \ker \rho_{x,y_2}$. Then, $\rho_{x,y_2}=c\rho_{x,y_1} = \rho_{x,\overline{c}y_1}$ for some $c \in \mathbb{C} \setminus \{0\}$. Since $\pi_\rho$ is irreducible, we have $\pi_\rho (A)x = y_2-cy_1$ for some $A \in \A$, which implies that
\[
0 = \rho_{x,y_2}(A)-\rho_{x,cy_1}(A) = \langle \pi_\rho (A)x,y_2-cy_1 \rangle = \|y_2-cy_1\|^2 .
\]
Hence, $y_2=cy_1$, that is, $\{y_1,y_2\}$ is linearly dependent. Conversely, if $\{y_1,y_2\}$ is linearly dependent, then $y_2 = cy_1$ for some $c \in \mathbb{C} \setminus \{0\}$. It follows from $\rho_{x,y_2} = \rho_{x,cy_1} = \overline{c}\rho_{x,y_1}$ that $\ker \rho_{x,y_1} = \ker \rho_{x,y_2}$. This proves (i).

(ii) We note that $\ker \rho_{x_j,y} = \{ A^* : A \in \ker \rho_{y,x_j}\}$ and $\ker \rho_{y,x_j} = \{ A^* : A \in \ker \rho_{x_j,y} \}$ for $j=1,2$. Hence, we have $\ker \rho_{x_1,y} = \ker \rho_{x_2,y}$ if and only if $\ker \rho_{y,x_1} = \ker \rho_{y,x_2}$, which occurs if and only if $\{x_1,x_2\}$ is linearly dependent by (i).
\end{proof}
\end{lemma}
\begin{lemma}\label{two-points}

Let $\A$ be a $C^*$-algebra, and let $\rho$ be a pure state of $\A$. Then, the following hold:
\begin{itemize}
\item[{\rm (i)}] Suppose that $x,y_1,y_2 \in S_{\HH_\rho}$. If $\{y_1,y_2\}$ is linearly independent, then
\[
\{\ker \rho_{x,y_1}, \ker \rho_{x,y_2}\} \not \in \mathfrak{C}(\A).
\]
\item[{\rm (ii)}] Suppose that $x_1,x_2,y \in S_{\HH_\rho}$. If $\{x_1,x_2\}$ is linearly independent, then
\[
\{\ker \rho_{x_1,y}, \ker \rho_{x_2,y}\} \not \in \mathfrak{C}(\A).
\]
\end{itemize}
\begin{proof}
Suppose that $\{y_1,y_2\}$ is linearly independent. Set $M=\langle \{y_1,y_2\} \rangle$, that is, $M$ is the linear span of $\{y_1,y_2\}$. If $y \in S_M$, then $y=\lambda y_1+\mu y_2$, which implies that $\rho_{x,y} = \overline{\lambda}\rho_{x,y_1}+\overline{\mu}\rho_{x,y_2}$. Hence, $\ker \rho_{x,y} \in \{ \ker \rho_{x,y_1},\ker \rho_{x,y_2} \}^=$. Meanwhile, $\ker \rho_{x,y} \not \in \{ \ker \rho_{x,y_1},\ker \rho_{x,y_2}\}$ if $\{y,y_j\}$ is linearly independent for $j=1,2$. Therefore, $\{\ker \rho_{x,y_1}, \ker \rho_{x,y_2}\} \not \in \mathfrak{C}(\A)$. This proves (i). An argument similar to above also shows (ii).
\end{proof}
\end{lemma}
Before characterizing the equivalence relation on $\s (\A)$, we give a technical remark on elements of $\ext (B_{\A^*})$.
\begin{remark}\label{abs-unique}

Let $\rho$ be a pure state of a $C*$-algebra $\A$. If $\rho U$ is also a pure state for some unitary element $U$ in $\uA$, then $\rho = \rho U$. Indeed since $\rho (A) = \langle \pi_\rho (A)x_\rho,x_\rho \rangle$ for each $A \in \A$, it follows that $\rho U = \rho_{x_\rho,y}$ for $y = \pi_\rho (U)^* x_\rho$. If $\pi_\rho (A) \geq 0$, then $\pi_\rho (A)^{1/2} \in \pi_\rho (\A)$ since $\pi_\rho (\A)$ is a $C^*$-subalgebra of $B(\HH_\rho)$. Let $\pi_\rho (B) = \pi_\rho (A)^{1/2}$. It may be assumed that $B^*=B$ since $\pi_\rho (2^{-1}(B+B^*)) = \pi_\rho (B)$, in which case, $B^2 \geq 0$ and $\pi_\rho (B^2) = \pi_\rho (B)^2 = \pi_\rho (A)$. Since $\rho U$ is positive, we obtain
\[
\langle \pi_\rho (A)x_\rho,y \rangle = \langle \pi_\rho (B^2)x_\rho,y \rangle = (\rho U)(B^2) \geq 0 .
\]
Meanwhile, $\pi_\rho (\A)$ is strong-operator dense in $B(\HH_\rho)$ since it acts irreducibly on $\HH_\rho$. By the Kadison transitivity theorem, each positive element of $B(\HH_\rho )$ is the strong-operator limit of a net of positive elements of $\pi_\rho (\A)$. Therefore, $\langle Ax_\rho,y \rangle \geq 0$ for each positive element $A$ of $B(\HH_\rho)$. This means that the functional $A \mapsto \langle Ax_\rho,y \rangle$ is a state of $B(\HH_\rho)$, which implies that $x_\rho=y$. Thus, $\rho = \rho U$.

As a consequence of the preceding paragraph, if $\rho \in \ext (B_{\A^*})$, then there exists a unique pure state $\tau$ of $\A$ such that $\rho = \tau U$ for some unitary element $U$ of $\uA$. Indeed, if $\tau_1,\tau_2$ are pure state of $\A$, and if $\tau_1 U_1 = \rho = \tau_2 U_2$ for some unitary elements $U_1,U_2$ of $\uA$, then $\tau_1 U_1U_2^* = \tau_2$, that is, $\tau_1 U_1U_2^*$ is a pure state of $\A$, which implies that $\tau_1 = \tau_1 U_1U_2^* = \tau_2$.
\end{remark}
Let $\A$ be a $C^*$-algebra, and let $\rho \in \ext (B_{\A^*})$. Then, $|\rho|$ denotes a unique pure state of $\A$ such that $\rho = |\rho |U$ for some unitary element $U$ of $\uA$.

If $\pi_1:\A \to B(\HH_1)$ and $\pi_2:\A \to B(\HH_2)$ are irreducible representations, then $(\pi_1\oplus \pi_2)(\A)$ acts non-degenerately on $ \HH_1 \oplus \HH_2$. Indeed, if $x \oplus y$ is a nonzero element of $\HH_1 \oplus \HH_2$, then $x \neq 0$ or $y \neq 0$. In the case that $x \neq 0$, there exists an $A \in \A$ such that $\pi_1 (A)x = x$. It follows that $((\pi_1 \oplus \pi_2 )(A))(x,y) = (\pi_1 (A)x , \pi_2(A)y) \neq 0$. We have the same consequence as above in the other case.

Recall that if $\pi :\A \to B(\HH)$ is a representation, then $\ker \pi$ is a norm-closed ideal of $\A$ and $\A /\ker \pi$ is $*$-isomorphic to $\pi (\A)$ under the $*$-isomorphism $A+\ker \pi \mapsto \pi (A)$. If $A \in \A$, then \cite[Proposition II.5.1.5]{Bla06} generates an element $A' \in \A$ such that $\pi (A) = \pi (A')$ and $\|A'\| = \|A+\ker \pi \|=\|\pi (A)\|$. This shows that $\pi (B_\A) = B_{\pi (\A)}$.

Now, the equivalence relation on $\s (\A)$ is characterized as follows.
\begin{theorem}\label{sim-equiv}

Let $\A$ be a $C^*$-algebra, and let $\rho_1 ,\rho_2 \in \ext (B_{\A^*})$. Then, $\ker \rho_1 \sim \ker \rho_2$ if and only if $\pi_{|\rho_1|}$ and $\pi_{|\rho_2|}$ are unitarily equivalent.
\begin{proof}
Let $U_1$ and $U_2$ be unitary elements of $\uA$ such that $\rho_j = |\rho_j|U_j$ for $j=1,2$. Suppose that $\pi_{|\rho_1|}$ and $\pi_{|\rho_2|}$ are not unitarily equivalent. In this case, the set $\{ \rho_1,\rho_2\}$ is linearly independent; in particular, $\ker \rho_1 \neq \ker \rho_2$. Indeed, if $\{ \rho_1,\rho_2\}$ is linearly dependent, then $\rho_1=c\rho_2 = c|\rho_2|U_2 = |\rho_2|(cU)$ for some $c \in\mathbb{C}$ with $|c|=1$, which implies that $|\rho_1|=|\rho_2|$, a contradiction. We note that $|\rho_j| = |\rho_j|_{x_{|\rho_j|},x_{|\rho_j|}}$ and $\rho_j = |\rho_j|_{x_{|\rho_j|},y_j}$ for $j=1,2$, where $y_j = \pi_{|\rho_j|}(U_j)^*x_{|\rho_j|}$. 

Since $\pi_{|\rho_1|}$ and $\pi_{|\rho_2|}$ are irreducible and not unitarily equivalent, they are disjoint from each other, which implies that the strong-operator closure of $(\pi_{|\rho_1|} \oplus \pi_{|\rho_2|})(\A)$ in $B(\HH_{|\rho_1|} \oplus \HH_{|\rho_2|})$ coincides with the direct sum $B(\HH_{|\rho_1|}) \oplus B(\HH_{|\rho_2|})$; see \cite[Corollary 2.2.5, Theorems 3.8.11 and 3.13.2, and 3.8.13]{Ped18}. Now, let $\ker \rho \in \{ \ker \rho_1,\ker \rho_2 \}^=$. Since $\ker \rho_1 \cap \ker \rho_2 \subset \ker \rho$, we have $\rho = c_1 \rho_1 + c_2 \rho_2$. Let $V_j$ be a unitary operator in $B(\HH_{|\rho_j|})$ such that $V_j(c_j x_{|\rho_j|}) = |c_j|y_j$. Then, $(V_1,V_2)$ is a unitary element of $B(\HH_{|\rho_1|}) \oplus B(\HH_{|\rho_2|})$. Meanwhile, since
\[
\rho_1 (A) = \langle \pi_{|\rho_1|}(A)x_{|\rho_1|},y_1 \rangle = \langle (\pi_{|\rho_1|}\oplus \pi_{|\rho_2|})(A)(x_{|\rho_1|},0),(y_1,0) \rangle
\]
and
\[
\rho_2 (A) = \langle \pi_{|\rho_2|}(A)x_{|\rho_2|},y_2 \rangle = \langle (\pi_{|\rho_1|}\oplus \pi_{|\rho_2|})(A)(0,x_{|\rho_2|}),(0,y_2) \rangle
\]
for each $A \in \A$, it follows that
\begin{align*}
\lefteqn{(c_1 \rho_1 + c_2 \rho_2)(A)}\\
&= \langle (\pi_{|\rho_1|}\oplus \pi_{|\rho_2|})(A)(c_1 x_{|\rho_1|},0),(y_1,0) \rangle + \langle (\pi_{|\rho_1|}\oplus \pi_{|\rho_2|})(A)(0,c_2 x_{|\rho_2|}),(0,y_2) \rangle
\end{align*}
for each $A \in \A$. Since $(V_1,V_2) \in B(\HH_{|\rho_1|}) \oplus B(\HH_{|\rho_2|})$, there exists a net $(B_\lambda )_\lambda$ in the unit ball of $(\pi_{|\rho_1|} \oplus \pi_{|\rho_2|})(\A)$ that converges to $(V_1,V_2)$ with respect to the strong-operator topology. Moreover, the unit ball of $(\pi_{|\rho_1|} \oplus \pi_{|\rho_2|})(\A)$ coincides with $(\pi_{|\rho_1|} \oplus \pi_{|\rho_2|})(B_\A)$. Hence, $B_\lambda = (\pi_{|\rho_1|}\oplus \pi_{|\rho_2|})(A_\lambda)$ for some $A_\lambda \in B_\A$ and
\begin{align*}
|c_1|+|c_2|
&= |\langle (V_1,V_2)(c_1 x_{|\rho_1|},0),(y_1,0) \rangle + \langle (V_1,V_2)(0,c_2 x_{|\rho_2|}),(0,y_2) \rangle| \\
&= \lim_\lambda |(c_1\rho_1+c_2\rho_2)(A_\lambda )| \\
&\leq \|c_1\rho_1+c_2\rho_2\| \\
&\leq |c_1|+|c_2| .
\end{align*}
This together with $\rho \in \ext (B_{\A^*})$ shows that $|c_1|+|c_2|=\|c_1\rho_1+c_2\rho_2\|=1$. If $c_1c_2 \neq 0$, we derive that
\[
\rho = |c_1|(|c_1|^{-1}c_1\rho_1)+|c_2|(|c_2|^{-1}c_2\rho_2) ,
\]
and $|c_1|^{-1}c_1\rho_1 \neq |c_2|^{-1}c_2\rho_2$ since $\{\rho_1,\rho_2\}$ is linearly independent. However, this contradicts $\rho \in ext (B_{\A^*})$. Thus, $c_1=0$ or $c_2=0$, that is, $\rho = c_1\rho_1$ or $\rho = c_2 \rho_2$, which shows that $\ker \rho \in \{\ker \rho_1,\ker \rho_2\}$. Hence, $\{\ker \rho_1,\ker \rho_2\} \in \mathfrak{C}(\A)$, and $\ker \rho_1 \not \rightleftarrows \ker \rho_2$.

Next, suppose that $\pi_{|\rho_1|}$ and $\pi_{|\rho_2|}$ are unitarily equivalent, and that $\ker \rho_1 \neq \ker \rho_2$. Let $U:\HH_1 \to \HH_2$ be an isometric isomorphism such that $\pi_{|\rho_2|}(A)= U\pi_{|\rho_1|}(A)U^*$ for each $A \in \A$. Then, it follows that
\[
\rho_2 (A) = \langle \pi_{|\rho_2|}(A)x_{|\rho_2|},y_2 \rangle = \langle U\pi_{|\rho_1|}(A)U^*x_{|\rho_2|},y_2 \rangle = \langle \pi_{|\rho_1|}(A)U^*x_{|\rho_2|},U^*y_2 \rangle
\]
for each $A \in \A$. Here, we note that either $\{x_{|\rho_1|},U^*x_{|\rho_2|}\}$ or $\{y_1,U^*y_2\}$ is linearly independent by $\ker \rho_1 \neq \ker \rho_2$. Indeed, if $\{x_{|\rho_1|},U^*x_{|\rho_2|}\}$ and $\{y_1,U^*y_2\}$ are both linearly dependent, then $U^*x_{|\rho_2|}=cx_{|\rho_1|}$ and $U^*y_2=dy_1$ for some $c,d \in \mathbb{C}$ with $|c|=|d|=1$, in which case, $\rho_2 = c\overline{d}\rho_1$. This contradicts $\ker \rho_1 \neq \ker \rho_2$. The argument will be divided into three parts.

(i) Suppose that $\{x_{|\rho_1|},U^*x_{|\rho_2|}\}$ is linearly dependent. Then, $\{y_1,U^*y_2\}$ is linearly independent. Let $c \in \mathbb{C}$ be such that $|c|=1$ and $U^*x_{|\rho_2|} = cx_{|\rho_1|}$. Then, we have
\[
\rho_2 (A) = \langle \pi_{|\rho_1|}(A)cx_{|\rho_1|},U^*y_2 \rangle = c|\rho_1|_{x_{|\rho_1|},U^*y_2}(A)
\]
for each $A \in \A$, which implies that $\ker \rho_2 = \ker |\rho_1|_{x_{|\rho_1|},U^*y_2}$. Meanwhile, Lemma~\ref{two-points} ensures that
\[
\{\ker |\rho_1|_{x_{|\rho_1|},y} ,\ker \rho_1 \} = \{\ker |\rho_1|_{x_{|\rho_1|},y} ,\ker |\rho_1|_{x_{|\rho_1|},y_1} \} \not \in \mathfrak{C}(\A)
\]
and
\[
\{\ker |\rho_1|_{x_{|\rho_1|},y} ,\ker \rho_2 \} = \{\ker |\rho_1|_{x_{|\rho_1|},y} ,\ker |\rho_1|_{x_{|\rho_1|},U^*y_2} \} \not \in \mathfrak{C}(\A)
\]
whenever $\{ y,y_1\}$ and $\{y,U^*y_2\}$ are linearly independent. We note that such an element $y \in S_{\HH_{|\rho_1|}}$ exists since $\{y_1,U^*y_2\}$ is linearly independent. This shows that $\ker \rho_1 \sim \ker \rho_2$.

(ii) Suppose that $\{y_1,U^*y_2\}$ is linearly dependent. Then, $\{x_{|\rho_1|},U^*x_{|\rho_2|}\}$ is linearly independent, and $U^*y_2 = dy_1$ for some $d \in \mathbb{C}$ with $|d|=1$. In particular, we obtain
\[
\rho_2 (A) = \langle \pi_{|\rho_1|}(A)U^*x_{|\rho_2|},U^*y_2 \rangle = \overline{d}|\rho_1|_{U^*x_{|\rho_2|},y_1}(A)
\]
for each $A \in \A$, that is, $\ker \rho_2 = \ker |\rho_1|_{U^*x_{|\rho_2|},y_1}$. Now, let $x \in S_{\HH_{|\rho_1|}}$ be such that $\{x,x_{|\rho_1|}\}$ and $\{x,U^*x_{|\rho_2|}\}$ are linearly independent. Then, it follows from Lemma~\ref{two-points} that
\[
\{\ker |\rho_1|_{x,y_1} ,\ker \rho_1 \} = \{\ker |\rho_1|_{x,y_1} ,\ker |\rho_1|_{x_{|\rho_1|},y_1} \} \not \in \mathfrak{C}(\A)
\]
and
\[
\{\ker |\rho_1|_{x,y_1} ,\ker \rho_2 \} = \{\ker |\rho_1|_{x,y_1} ,\ker |\rho_1|_{U^*x_{|\rho_2|},y_1} \} \not \in \mathfrak{C}(\A) .
\]
Therefore, $\ker \rho_1 \sim \ker \rho_2$.

(iii) Suppose that $\{x_{|\rho_1|},U^*x_{|\rho_2|}\}$ and $\{y_1,U^*y_2\}$ are linearly independent. In this case, we derive that
\[
\{ \ker |\rho_1|_{x_{|\rho_1|},U^*y_2},\ker \rho_1 \} = \{ \ker |\rho_1|_{x_{|\rho_1|},U^*y_2},\ker |\rho_1|_{x_{|\rho_1|},y_1} \} \not \in \mathfrak{C}(\A)
\]
and
\[
\{ \ker |\rho_1|_{x_{|\rho_1|},U^*y_2},\ker \rho_2 \} = \{ \ker |\rho_1|_{x_{|\rho_1|},U^*y_2},\ker |\rho_1|_{U^*x_{|\rho_2|},U^*y_2} \} \not \in \mathfrak{C}(\A)
\]
by Lemma~\ref{two-points}. This proves that $\ker \rho_1 \sim \ker \rho_2$. Thus, in either case, the unitary equivalence of $\pi_{|\rho_1|}$ and $\pi_{|\rho_2|}$ implies $\ker \rho_1 \sim \ker \rho_2$.

Now, we know the following two implications: $\ker \rho_1 \rightleftarrows \ker \rho_2$ implies that $\pi_{|\rho_1|}$ and $\pi_{|\rho_2|}$ are unitarily equivalent, and the unitary equivalence of $\pi_{|\rho_1|}$ and $\pi_{|\rho_2|}$ implies that $\ker \rho_1 \sim \ker \rho_2$. Since the unitary equivalence of representations of $\A$ is transitive, we can conclude that $\ker \rho_1 \sim \ker \rho_2$ if and only if $\pi_{|\rho_1|}$ and $\pi_{|\rho_2|}$ are unitarily equivalent.
\end{proof}
\end{theorem}
\begin{corollary}\label{polar-equi}

Let $\A$ be a $C^*$-algebra, and let $\rho \in \ext (B_{\A^*})$. Then, $\ker \rho \sim \ker |\rho|$. In particular, each element of $\mathfrak{GS}(\A)$ contains a pure state of $\A$.
\end{corollary}
Based on Theorem~\ref{sim-equiv}, we can find a natural relationship between elements of geometric primitive spectra and primitive ideals of $C^*$-algebras.
\begin{lemma}\label{primitive}

Let $\A$ be a $C^*$-algebra, and let $\rho$ be a pure state of $\A$. Then, $C(\ker \rho) = \{ \ker \rho_{x,y} : x,y \in S_{\HH_\rho}\}$ and $J_{C(\ker \rho)} = \ker \pi_\rho$.
\begin{proof}
Suppose that $\tau \in \ext (B_{\A^*})$ and $\ker \tau \in C(\ker \rho)$. Let $U$ be a unitary element of $\uA$ such that $\tau = |\tau |U$. Since $\ker \rho \sim \ker \tau$, the representations $\pi_\rho$ and $\pi_{|\tau|}$ are unitarily equivalent by Theorem~\ref{sim-equiv}. Let $U_0 :\HH_\rho \to \HH_{|\tau|}$ be an isometric isomorphism such that $\pi_{|\tau|}(A) = U_0\pi_\rho (A) U_0^*$ for each $A \in \A$. Then, it follows that
\begin{align*}
\tau (A) = |\tau |(UA) = \langle \pi_{|\tau|}(UA)x_{|\tau|},x_{|\tau|}\rangle
&= \langle U_0\pi_\rho (UA)U_0^*x_{|\tau|},x_{|\tau|}\rangle \\
&= \langle \pi_\rho (A)U_0^*x_{|\tau|},\pi_\rho (U)^*U_0^*x_{|\tau|} \rangle
\end{align*}
for each $A \in \A$. Therefore, setting $x=U_0^*x_{|\tau|}$ and $y = \pi_\rho (U)^*U_0^*x_{|\tau|}$, we have $\tau = \rho_{x,y}$ and $\ker \tau = \ker \rho_{x,y}$.

Conversely, suppose that $x,y \in S_{\HH_\rho}$. Since $\pi_\rho$ is an irreducible representation of $\uA$, there exists a unitary element $U$ of $\uA$ such that $\pi_\rho (U)x = y$, which implies that $\rho_{x,y}=\rho_{x,x}U^*$. In particular, we obtain $|\rho_{x,y}| = \rho_{x,x}$. Moreover, since each unit vector in $\HH_\rho$ is cyclic for $\pi_\rho (\A)$, the representations $\pi_\rho$ and $\pi_{\rho_{x,x}}$ are unitarily equivalent by the essential uniqueness of the GNS construction. Therefore, $\ker \rho \sim \ker \rho_{x,x} \sim \ker \rho_{x,y}$ by Theorem~\ref{sim-equiv} and Corollary~\ref{polar-equi}. This shows that $\ker \rho_{x,y}\in C(\ker \rho)$ for each $x,y \in S_{\HH_\rho}$. Thus, $C(\ker \rho) = \{ \ker \rho_{x,y} : x,y \in S_{\HH_\rho}\}$.

Finally, we note that $A \in J_{C(\ker \rho)}$ in and only if $\langle \pi_\rho (A)x,y \rangle = \rho_{x,y}(A) = 0$ for each $x,y \in S_{\HH_\rho}$, which occurs if and only if $\pi_\rho (A) =0$. Hence, $J_{C(\ker \rho)} = \ker \pi_\rho$.
\end{proof}
\end{lemma}
Now, we give the main theorem in this section.

\begin{theorem}\label{main-theorem}

Let $\A$ be a $C^*$-algebra. Then, the spectrum $\hat{\A}$ of $\A$ is homeomorphic to $\mathfrak{GS}(\A)$, and the primitive spectrum $\Prim (A)$ of $\A$ coincides with $\mathfrak{PS}(\A)$. Moreover, a natural homeomorphism $\Phi :\hat{\A} \to \mathfrak{GS}(\A)$ is given by $\Phi (\overline{\pi_\rho}) = C(\ker \rho)$ for each pure state $\rho$ of $\A$, where $\overline{\pi_\rho}$ is the unitary equivalence class of $\pi_\rho$.
\begin{proof}
Recall that a subset $\mathcal{I}$ of $\A$ is a primitive ideal of $\A$ if and only if $\mathcal{I} = \ker \pi_\rho$ for some pure state $\rho$ of $\A$. Moreover, we have $\ker \pi_\rho = J_{C(\ker \rho)}$ by Lemma~\ref{primitive}. Hence, it follows from Corollary~\ref{polar-equi} that $\Prim (\A ) = \mathfrak{PS}(\A)$ as sets. Then, it is obvious that $S^{\cong}$ coincides with the hull-kernel closure of $S$ in $\Prim (\A)$.

Next, suppose that $\pi$ is an irreducible representation. Then, there exists a pure state $\rho$ of $\A$ such that $\pi$ is unitarily equivalent to $\pi_\rho$. Let $\Phi (\overline{\pi_\rho}) = C(\ker \rho)$ for each pure state $\rho$ of $\A$. Since $\overline{\pi_{\rho_1}}=\overline{\pi_{\rho_2}}$ if and only if $C(\ker \rho_1) = C(\ker \rho_2)$ by Theorem~\ref{sim-equiv}, the mapping $\Phi : \hat{\A} \to \mathfrak{GS}(\A)$ is well-defined and injective. Moreover, $\Phi$ is surjective by Corollary~\ref{polar-equi}.

Now, let $\kappa (\overline{\pi}) = \ker \pi$ for each $\overline{\pi} \in \hat{\A}$. Then, the topology on $\hat{\A}$ is defined as follows: A subset $F$ of $\hat{\A}$ is closed if and only if $C= \kappa^{-1}(F)$ for some hull-kernel closed subset $F$ of $\Prim (\A) = \mathfrak{PS}(\A)$. Meanwhile, we have $C' \in \mathfrak{GC}(\A)$ if and only if $C'=\gamma_\A^{-1}(F')$ for some $F' \in \mathfrak{PC}(\A)$ by Lemma~\ref{pullback} (i) and Theorem~\ref{quotient}, for $C' = C'^{\equiv} = \gamma_\A^{-1}(\gamma_\A (C')^{\cong})$ provided that $C' \in \mathfrak{GC}(\A)$. We claim that $\Phi (\kappa^{-1}(S)) = \gamma_\A^{-1}(S)$ for each $S \subset \mathfrak{PS}(\A)$. Indeed, if $C(\ker \rho) \in \Phi (\kappa^{-1}(S))$, then $C(\ker \rho ) = \Phi (\overline{\pi_\rho})$ for some $\overline{\pi_\rho} \in \kappa^{-1}(S)$, in which case, $\kappa (\overline{\pi_\rho}) \in S$. Since $\kappa (\overline{\pi_\rho}) = \ker \pi_\rho = J_{C(\ker \rho)} = \gamma_\A (C(\ker \rho ))$, it follows that $C(\ker \rho) \in \gamma_\A^{-1}(S)$. Therefore, $\Phi (\kappa^{-1}(S)) \subset \gamma_\A^{-1}(S)$. Conversely, if $C(\ker \rho) \in \gamma_\A^{-1}(S)$, then $\kappa (\overline{\pi_\rho}) = \gamma_\A (C(\ker \rho )) \in S$, which implies that $C(\ker \rho) = \Phi (\overline{\pi_\rho}) \in \Phi (\kappa^{-1}(S))$. This proves the claim.

Finally, if $C$ is closed in $\hat{\A}$, then $C = \kappa^{-1}(F)$ for some $F \in \mathfrak{PC}(\A)$ by the preceding paragraph. Since $\Phi (C) = \Phi(\kappa^{-1}(F)) = \gamma_\A^{-1}(F)$, it follows that $\Phi (C) \in \mathfrak{GC}(\A)$. Conversely, if $\Phi (C) \in \mathfrak{GC}(\A)$, then $\Phi (C) = \gamma_\A^{-1}(F)$ for some $F \in \mathfrak{PC}(\A)$, which together with $\gamma_\A^{-1}(F) = \Phi (\kappa^{-1}(F))$ implies that $C = \kappa^{-1}(F)$, that is, $C$ is closed in $\hat{\A}$. This shows that $\Phi$ preserves closed sets in both directions. Thus, $\Phi$ is a homeomorphism.
\end{proof}
\end{theorem}
Since $\s (\A),\gs (\A)$ and $\ps (\A)$ are constructed from the Banach space structure of $\A$, the preceding theorem gives a Banach space theoretical construction of the (primitive) spectrum of $\A$. Moreover, in the light of Theorem~\ref{spec-trans}, the spectrum $\hat{\A}$ and primitive spectrum $\Prim (\A)$ of $\A$ are homeomorphic to $\g$- and $\p$-transforms of $\s (\A)$, respectively.

Here, we recall that (primitive) spectra of $C^*$-algebras are always topologizable. Combining this with Theorem~\ref{main-theorem}, we immediately have the following result.
\begin{corollary}\label{top-able}

Let $\A$ be a $C^*$-algebra. Then, $\mathfrak{GS}(\A)$ and $\mathfrak{PS}(\A)$ are topologizable.
\end{corollary}
The following is a consequence of Corollary~\ref{s-ps} and Theorem~\ref{main-theorem}.
\begin{corollary}\label{spec-preserved}

Let $\A,\B$ be $C^*$-algebras. If $\s (\A)$ and $\s (\B)$ are homeomorphic, then so are the (primitive) spectra of $\A$ and $\B$.
\end{corollary}
This will be a useful tool for classifying $C^*$-algebras in terms of geometric structure spaces.

The rest of this section focuses on geometric structure spaces of images of irreducible representations of $C^*$-algebras. We begin with a simple lemma.
\begin{lemma}\label{pure-unitary}

Let $\A$ be a $C^*$-algebra, and let $\rho$ be a pure state of $\A$. Suppose that $B \in B_{\uA}$ and $\|\rho B\|=1$. Then, $\rho B \in \ext (B_{\A^*})$.
\begin{proof}
We note that $\rho B = \rho_{x_\rho,y}$ for $y = \pi_\rho (B)^*x_\rho$. Moreover, it follows from $\|B\| \leq 1$ and $\|\rho B\|=1$ that $\|y\|=1$. By the irreducibility of $\pi_\rho$, there exists a unitary element $U \in \uA$ such that $\pi_\rho (U)x_\rho = y$, which implies that $\rho B = \rho_{x_\rho,y} = \rho U^* \in \ext (B_{\A^*})$.
\end{proof}
\end{lemma}
We also need the following basic lemma.

\begin{lemma}\label{quotient-unitary}

Let $\A$ be a $C^*$-algebra, and let $\mathcal{I}$ be a closed ideal of $\A$. Then, the unitization of $\A/\mathcal{I}$ is $*$-isomorphic to $\uA /\mathcal{I}$ under the $*$-isomorphism $(A+\mathcal{I})+c\mathbf{1}' \mapsto (A+c\mathbf{1})+\mathcal{I}$ from the unitization of $\A/\mathcal{I}$ onto $\uA /\mathcal{I}$, where $\mathbf{1}$ and $\mathbf{1}'$ are the units in $\uA$ and the unitization of $\A/\mathcal{I}$ respectively. Consequently, $(A+\mathcal{I})+c\mathbf{1}'$ is unitary in the unitization of $\A/\mathcal{I}$ if and only if $(A+c\mathbf{1})+\mathcal{I}$ is unitary in $\uA /\mathcal{I}$.
\begin{proof}
If $\A$ is unital, then $\A /\mathcal{I}$ is also unital. Hence, we consider the non-unital case. First, we note that $\mathcal{I}$ is also an ideal of $\uA$. To show this, suppose that $A \in \mathcal{I}$ and $A \geq 0$. Since $\mathcal{I}$ is a $C^*$-algebra and $\A$ is an ideal of $\uA$, we have $A^{1/2} \in \mathcal{I}$, $BA = (BA^{1/2})A^{1/2} \in \mathcal{I}$ and $AB = A^{1/2}(A^{1/2}B) \in \mathcal{I}$ for each $B \in \uA$. Combining this with the fact that $\mathcal{I}$ coincides with the linear span of its positive elements, we have $BA \in \mathcal{I}$ and $AB \in \mathcal{I}$ whenever $A \in I$ and $B \in \uA$. We note that $A \in \A$ and $A+c\mathbf{1} \in \mathcal{I}$ implies that $c=0$; for otherwise, $I = c^{-1}(B-A) \in \A$ for some $B \in \mathcal{I} \subset \A$.

Next, let $A,B \in \A$ and $c,d \in \mathbb{C}$. Then, $(A+\mathcal{I})+c\mathbf{1}'=(B+\mathcal{I})+d\mathbf{1}'$ if and only if $A-B \in \mathcal{I}$ and $c=d$, which occurs if and only if $(A+c\mathbf{1})-(B+d\mathbf{1}) \in \mathcal{I}$, that is, $(A+c\mathbf{1})+\mathcal{I} = (B+d\mathbf{1})+\mathcal{I}$. Hence, setting $\varphi ((A+\mathcal{I})+c\mathbf{1}') = (A+c\mathbf{1})+\mathcal{I}$ for each $A \in \A$ and each $c \in \mathbb{C}$, we have a bijection $\varphi$ from the unitization of $\A /\mathcal{I}$ onto $\uA /\mathcal{I}$. Moreover, it follows that
\begin{align*}
\varphi ((A+\mathcal{I})+c\mathbf{1}')^*) = \varphi ((A^*+\mathcal{I})+\overline{c}\mathbf{1}') = (A^*+\overline{c}\mathbf{1})+\mathcal{I} 
&= (A+c\mathbf{1})^*+\mathcal{I} \\
&= \varphi ((A+\mathcal{I})+c\mathbf{1}')^* ,
\end{align*}
and 
\begin{align*}
\varphi (((A+\mathcal{I})+c\mathbf{1}')((B+\mathcal{I})+d\mathbf{1}')) 
&= \varphi ((AB+cB+dA+\mathcal{I})+cd\mathbf{1}') \\
&= (AB+cB+dA+cd\mathbf{1})+\mathcal{I} \\
&= (A+c\mathbf{1})(B+d\mathbf{1})+\mathcal{I} \\
&= ((A+c\mathbf{1})+\mathcal{I})((B+d\mathbf{1})+\mathcal{I}) \\
&= \varphi ((A+\mathcal{I})+c\mathbf{1}')\varphi ((B+\mathcal{I})+d\mathbf{1}') .
\end{align*}
Therefore, $\varphi$ is a $*$-isomorphism.
\end{proof}
\end{lemma}
Applying the preceding lemma to primitive ideals, we have the following identity. Recall that, for a subset $M$ of $\A$, the symbol $M^\perp$ stands for its annihilator, that is, $M^\perp = \{ \rho \in \A^* : M \subset \ker \rho\}$.
\begin{lemma}\label{quotient-ext}

Let $\A$ be a $C^*$-algebra, and let $\pi : \A \to B(\HH)$ be an irreducible representation. Then, $\ext (B_{(\ker \pi )^\perp}) = \ext (B_{\A^*}) \cap (\ker \pi )^\perp$.
\begin{proof}
Suppose that $\rho \in \ext (B_{(\ker \pi)^\perp})$. First, we assume that $\rho$ is a state of $\A$. If $\rho_1,\rho_2$ are states of $\A$ such that $\rho = (1-t)\rho_1 +t\rho_2$ for some $t \in (0,1)$, and if $A$ is a positive element of $\ker \pi$, then $\rho_1 (A) = \rho_2(A) = 0$ by $\rho_1(A) \geq 0$, $\rho_2 (A) \geq 0$ and $(1-t)\rho_1(A)+t\rho_2(A) = \rho (A) = 0$. Since $\ker \pi$ is a $C^*$-subalgebra of $\A$, every element of $\ker \pi$ is a linear combination of (at most four) positive elements of $\ker \pi$. Therefore, $\rho_j \in (\ker \pi )^\perp$ for $j=1,2$, which together with $\rho \in \ext (B_{(\ker \pi)^\perp})$ shows that $\rho =\rho_1 =\rho_2$. This proves that $\rho$ is a pure state of $\A$, and $\rho \in \ext (B_{\A^*}) \cap (\ker \pi )^\perp$.

Next, let $(\varphi (\rho ))(A+\ker \pi ) = \rho (A)$ for each $\rho \in (\ker \pi)^\perp$ and each $A\in \A$. Then, $\varphi : (\ker \pi )^\perp \to (\A /\ker \pi )^*$ is an isometric isomorphism; see, for example, \cite[Theorem 1.10.17]{Meg98}. Suppose that $\rho \in \ext (B_{(\ker \pi)^\perp})$. We note that $\varphi (\rho ) \in \ext (B_{(\A /\ker \pi)^*})$, and hence, $\varphi (\rho )V$ is a pure state of $\A /\ker \pi$ for some unitary element $V$ of the unitization of $\A /\ker \pi$. Set $V = (A_0+\ker \pi )+c\mathbf{1}'$, where $A \in \A$, $c \in \mathbb{C}$ and $\mathbf{1}'$ is the unit of the unitization of $\A /\ker \pi$. By Lemma~\ref{quotient-unitary}, $(A_0+c\mathbf{1})+\ker \pi$ is a unitary element of $\uA /\ker \pi$. In particular, we may assume that $\|A_0+c\mathbf{1}\|=1$ by \cite[Proposition II.5.1.5]{Bla06}. Set $U=A_0+c\mathbf{1} \in S_{\uA}$. Since $U+\ker \pi$ is unitary in $\uA /\ker \pi$, we have $\mathbf{1}-U^*U \in \ker \pi$ and $\mathbf{1}-UU^* \in \ker \pi$. Moreover, it turns out that
\begin{align*}
V(A+\ker \pi ) 
&= ((A_0+\ker \pi )+c\mathbf{1}')(A+\ker \pi ) \\
&= (A_0+\ker \pi)(A+\ker \pi) + c(A+\ker \pi) \\
&= (A_0A+cA) +\ker \pi \\
&= UA + \ker \pi
\end{align*}
for each $A \in \A$, which implies that
\[
(\varphi (\rho )V)(A + \ker \pi ) = (\varphi (\rho ))(UA+\ker \pi) = \rho (UA) = (\rho U)(A)
\]
for each $A \in \A$. Hence, $\rho U$ is positive. We claim that the mapping $\tau \mapsto \tau U$ is an isometric isomorphism on $(\ker \pi )^\perp$. Since $\ker \pi$ is an ideal of $\uA$, we obtain $\tau U \in (\ker \pi )^\perp$ whenever $\tau \in (\ker \pi)^\perp$. We note that $\|\tau U\| \leq \|\tau \|$ since $\|U\|=1$. Similarly, $\tau U^* \in (\ker \pi )^\perp$ whenever $\tau \in (\ker \pi )^\perp$, and $\|\tau U^*\| \leq \|\tau\|$ by $\|U^*\|=\|U\|=1$. Meanwhile since $\mathbf{1}-U^*U \in \ker \pi$ and $\mathbf{1}-UU^* \in \ker \pi$, it follows that $A-U^*UA \in \ker \pi$ and $A-UU^*A \in \ker \pi$ for each $A \in \A$, that is, $\tau = (\tau U^*)U = (\tau U)U^*$ whenever $\tau \in (\ker \pi )^\perp$. Thus, the mapping $\tau \mapsto \tau U$ is an isometric isomorphism on $(\ker \pi )^\perp$. Consequently, $\rho U \in \ext (B_{(\ker \pi)^\perp})$ as $\rho \in \ext (B_{(\ker \pi)^\perp})$. Combining this with the the first paragraph, it turns out that $\rho U$ is a pure state of $\A$. Now, Lemma~\ref{pure-unitary} applies, and $\rho = (\rho U)U^* \in \ext (B_{\A^*})$.

The converse is obvious since $(\ker \pi)^\perp$ is a subspace of $\A^*$.
\end{proof}
\end{lemma}
Based on the preceding lemma, we can identify the geometric structure space of the image of an irreducible representation of a $C^*$-algebra. Recall that, for a bounded linear operator $T$ from a Banach space $X$ onto another $Y$, the (Banach space) adjoint $T^*$ of $T$ is defined by $(T^*f)(x) = f(Tx) = (f\circ T)(x)$ for each $f \in X^*$ and each $x \in X$.
\begin{theorem}\label{irr-rep}

Let $\A$ be a $C^*$-algebra, let $\pi$ be an irreducible representation of $\A$, and let $I \in \mathfrak{S}(\A)$ be such that $\ker \pi = J_{C(I)}$. Then, $\pi^* : \pi (\A)^* \to (\ker \pi )^\perp$ is an isometric weak$^*$-to-(relative) weak$^*$ homeomorphism, and 
\[
\mathfrak{S}(\pi (\A)) = \{ \ker \rho : \pi^* (\rho ) \in \ext (B_{\A^*})  \}
\]
Moreover, the mapping $\ker \rho \mapsto \ker \pi^* (\rho )$ is a homeomorphism from $\mathfrak{S}(\pi (\A))$ onto the subspace $C(I)^=$ of $\mathfrak{S}(\A)$.
\begin{proof}
We first note that $\pi^* (\rho ) = \rho \circ \pi \in (\ker \pi)^\perp$ whenever $\rho \in \pi (\A)^*$. Let $\psi (A + \ker \pi ) = \pi (A)$ and $(\varphi (\rho ))(A+\ker \pi )=\rho (A)$ for each $A \in \A$ and each $\rho \in (\ker \pi )^\perp$. Then, $\psi :\A /\ker \pi \to \pi (\A)$ is a $*$-isomorphism and $\varphi : (\ker \pi )^\perp \to (\A /\ker \pi)^*$ is an isometric isomorphism, which implies that $(\psi^{-1})^* \circ \varphi : (\ker \pi)^\perp \to \pi (\A)^*$ is also an isometric isomorphism. We note that since $\varphi$ is a (relative) weak$^*$-to-weak$^*$ homeomorphism, so is $(\psi^{-1})^* \circ \varphi$. Moreover, setting $\overline{\rho} = ((\psi^{-1})^* \circ \varphi)(\rho)$, we obtain
\begin{align*}
\overline{\rho}(\pi (A)) = ((\psi^{-1})^*\circ \varphi )(\rho ))(\pi (A)) &= ((\psi^{-1})^*(\varphi (\rho )))(\pi (A)) \\
&= (\varphi (\rho ))(\psi^{-1}(\pi (A))) \\
&= \rho (A)
\end{align*}
for each $A \in \A$ and each $\rho \in (\ker \pi )^\perp$, that is, $\pi^* (\overline{\rho}) = \rho$ for each $\rho \in (\ker \pi )^\perp$. This shows that $\pi^* = ((\psi^{-1})^* \circ \varphi)^{-1}$, and $\pi^* : \pi (\A)^* \to (\ker \pi )^\perp$ is an isometric weak$^*$-to-(relative) weak$^*$ homeomorphism, which together with Lemma~\ref{quotient-ext} implies that
\begin{align*}
\mathfrak{S}(\pi (\A)) = \{ \ker \rho : \rho \in \ext (B_{\pi (\A )^*}) \} &= \{ \ker \rho : \pi^*(\rho ) \in \ext (B_{(\ker \pi)^\perp}) \} \\
&= \{ \ker \rho : \pi^*(\rho ) \in \ext (B_{\A^*}) \}.
\end{align*}
Moreover, for an element $\tau \in \ext (B_{\A^*})$, it follows from $\ker \pi = J_{C(I)}$ and Lemma~\ref{closed} that $\tau \in (\ker \pi )^\perp$ if and only if $\ker \tau \in C(I)^=$. Meanwhile, if $\ker \tau \in C(I)^=$, then $\tau \in \ext (B_{\A^*})$ implies that $\tau = \pi^* (\rho )$ for some $\rho \in \ext (B_{\pi (\A)^*})$. Therefore, $C(I)^= = \{ \ker \pi^* (\rho ) : \ker \rho \in \mathfrak{S}(\pi (\A)) \}$. We note that, for elements $\ker \rho_1$ and $\ker \rho_2$ of $\mathfrak{S}(\pi (\A))$, the equality $\ker \rho_1 = \ker \rho_2$ holds if and only if $\rho_2 = c\rho_1$ for some nonzero $c \in \mathbb{C}$, which occurs if and only if $\pi^* (\rho_2)=c\pi^* (\rho_1)$ for some nonzero $c \in \mathbb{C}$, that is, the equality $\ker \pi^*(\rho_1) = \ker \pi^* (\rho_2 )$ holds. Thus, setting $\Phi (\ker \rho) = \ker \pi^* (\rho )$ for each $\ker \rho \in \mathfrak{S}(\pi (\A))$ yields a bijection from $\mathfrak{S}(\pi (\A))$ onto $C(I)^=$.

Finally, suppose that $S \subset \mathfrak{S}(\pi (\A))$ and $\ker \rho \in \mathfrak{S}(\pi (\A))$. We assume that $\ker \rho \in S^=$. If $A \in \bigcap \{ \ker \pi^* (\tau ) : \ker \tau \in S \}$, then $\tau (\pi (A)) = (\pi^*(\tau))(A) = 0$ whenever $\ker \tau \in S$, that is, $\pi (A) \in \bigcap \{ \ker \tau : \ker \tau \in S \} \subset \ker \rho$. This means that $(\pi^* (\rho ))(A) = \rho (\pi (A)) = 0$. Hence,
\[
\bigcap \{ \Phi (\ker \tau ) : \ker \tau \in S \} = \bigcap \{ \ker \pi^*(\tau ) : \ker \tau \in S \} \subset \ker \pi^* (\rho ) = \Phi (\ker \rho ) ,
\]
that is, $\Phi (\ker \rho) \in \Phi (S)^=$. Conversely, if $\ker \rho \not \in S^=$, then there exists a $\pi (A) \in \pi (\A)$ such that $\pi (A) \in \bigcap \{ \ker \tau : \ker \tau \in S\}$ and $\pi (A) \not \in \ker \rho$, in which case, we have $(\pi^*(\tau ))(A) = \tau (\pi (A)) = 0$ whenever $\ker \tau \in S$, and $(\pi^* (\rho))(A) = \rho (\pi (A)) \neq 0$. Therefore, $A \in \bigcap \{ \ker \pi^* (\tau ) : \ker \tau \in S\}$ and $A \not \in \ker \pi^* (\rho)$, that is, $\Phi (\ker \rho) = \ker \pi^* (\rho ) \not \in \Phi (S)^=$. Now, we know that $\ker \rho \in S^=$ if and only if $\Phi (\ker \rho) \in \Phi (S)^=$, which shows that $\Phi (S^=) = \Phi (S)^=$. Thus, $\Phi$ is a homeomorphism.
\end{proof}
\end{theorem}
Let $\HH$ be a Hilbert space, and let $x,y \in \HH$. In the rest of this paper, let $\omega_{x,y}$ denote the weak-operator continuous functional on $B(\HH)$ defined by $\omega_{x,y}(A) = \langle Ax,y \rangle$ for each $A \in B(\HH)$.
\begin{corollary}\label{vector-funct}

Let $\A$ be a $C^*$-algebra, let $\pi$ be an irreducible representation of $\A$, let $\HH$ be the Hilbert space on which $\pi (\A)$ acts, and let $I \in \mathfrak{S}(\A)$ be such that $\ker \pi = J_{C(I)}$. Then, $\mathfrak{S}(\pi (\A)) \supset \{ \ker \omega_{x,y} \cap \pi (\A): x,y \in S_{\HH} \}$.
\begin{proof}
Let $\rho$ be a pure state of $\A$ such that $\pi_\rho$ is unitarily equivalent to $\pi$, and let $U : \HH \to \HH_\rho$ be an isometric isomorphism such that $\pi (A) = U^*\pi_\rho (A)U$ for each $A \in \A$. Suppose that $x,y \in S_{\HH}$. Then, it follows from
\begin{align*}
\omega_{x,y}(\pi (A)) = \langle \pi (A)x,y \rangle = \langle U^*\pi_\rho (A)Ux,y \rangle = \rho_{Ux,Uy}(A)
\end{align*}
that $\pi^* (\omega_{x,y}|\pi (\A)) = \rho_{Ux,Uy} \in \ext (B_{\A^*})$ by Lemma~\ref{ext-xy}. Therefore,
\[
\ker \omega_{x,y} \cap \pi (\A) = \ker (\omega_{x,y}|\pi (\A)) \in \mathfrak{S}(\pi (A))
\]
by Theorem~\ref{irr-rep}.
\end{proof}
\end{corollary}
\begin{corollary}\label{irr-alg}

Let $\HH$ be a Hilbert space, and let $\A$ be a $C^*$-subalgebra of $B(\HH)$ that acts irreducibly on $\HH$. Then, $\mathfrak{S}(\A ) \supset \{ \ker \omega_{x,y} \cap \A : x,y \in S_{\HH} \}$.
\begin{proof}
It is sufficient to apply the preceding corollary to the identity representation of $\A$.
\end{proof}
\end{corollary}


\section{Nonlinear theory for irreducible $C^*$-algebras}\label{Sect:Irr}

An \emph{irreducible} $C^*$-algebra is a $C^*$-subalgebra of some $B(\HH)$ that acts irreducibly on $\HH$. This section is devoted to special mappings between geometric structure spaces of irreducible $C^*$-algebras that play crucial roles in the rest of this paper. We begin with modifications of Lemma~\ref{primitive}.
\begin{lemma}\label{pure-adjoint}

Let $\A$ be a $C^*$-algebra, and let $\rho$ be a pure state of $\A$. Then, $C(\ker \rho ) = \{ \ker \pi_\rho^* (\omega_{x,y}|\pi_\rho (\A)) : x,y \in S_{\HH_\rho} \}$.
\begin{proof}
Recall that $C(\ker \rho) = \{ \ker \rho_{x,y} : x,y \in S_{\HH_\rho}\}$ by Lemma~\ref{primitive}. Since $\rho_{x,y} = \omega_{x,y} \circ \pi_\rho = \pi_\rho^* (\omega_{x,y}|\pi_\rho (\A))$, it follows that $C(\ker \rho) = \{ \ker \pi_\rho^* (\omega_{x,y}|\pi_\rho (\A)) : x,y \in S_{\HH_\rho}\}$.
\end{proof}
\end{lemma}
\begin{remark}

Let $\HH$ be a Hilbert space. Then, $\omega_{x,y} = \|x\|\|y\|\omega_{\|x\|^{-1}x,\|y\|^{-1}y}$ and $\ker \omega_{x,y} = \ker \omega_{\|x\|^{-1}x,\|y\|^{-1}y}$ whenever $x,y \in \HH \setminus \{0\}$. Hence, the condition $x,y \in S_{\HH_\rho}$ is not essential in the preceding lemma, that is, the identity
\[
\{ \ker \pi_\rho^* (\omega_{x,y}|\pi_\rho (\A)) : x,y \in S_{\HH_\rho} \} = \{ \ker \pi_\rho^* (\omega_{x,y}|\pi_\rho (\A)) : x,y \in \HH_\rho \setminus \{0\} \}
\]
holds.
\end{remark}
\begin{lemma}\label{normal-class}

Let $\A$ be an irreducible $C^*$ algebra acting on a Hilbert space $\HH$. Suppose that $x_0,y_0 \in \HH \setminus \{0\}$. Then, $C(\ker \omega_{x_0,y_0} \cap \A) = \{ \ker \omega_{x,y} \cap \A : x,y \in \HH \setminus \{0\}\}$.
\begin{proof}
It may be assumed that $x_0,y_0 \in S_{\HH}$. Then, $\omega_{x_0,y_0},\omega_{x_0,x_0} \in \ext (B_{\A^*})$ by Corollary~\ref{irr-alg}. Moreover, since $\A$ acts irreducibly on $\HH$, we have a unitary element $U \in \uA$ such that $U^*y_0=x_0$, that is, $\omega_{x_0,y_0}U = \omega_{x_0,x_0}$. This means that $|\omega_{x_0,y_0} | \A |=\omega_{x_0,x_0}|\A$, and $C(\ker \omega_{x_0,y_0}\cap \A) =\ker (\omega_{x_0,x_0}\cap \A)$ by Theorem~\ref{sim-equiv}. Set $\rho = \omega_{x_0,x_0}|\A$ for short. Since $\rho$ is a pure state of $\A$, by Lemma~\ref{pure-adjoint}, we obtain
\[
C(\ker \omega_{x_0,x_0} \cap \A) = \{ \ker \pi_\rho^* (\omega_{z,w}|\pi_\rho (\A)) : z,w \in \HH_\rho \setminus \{0\} \} .
\]
Now, let $x,y \in \HH \setminus \{0\}$. By the irreducibility of $\A$, there exist $V_1,V_2 \in \A$ such that $V_1x_0=x$ and $V_2y_0 = y$. It follows from $\omega_{x,y} \neq 0$ and
\begin{align*}
\omega_{x,y}(A) = \omega_{x_0,x_0}(V_2^*AV_1) 
&= \langle \pi_\rho (A)\pi_\rho (V_1)x_\rho,\pi_\rho (V_2)x_\rho \rangle \\
&= (\pi_\rho^* (\omega_{\pi_\rho (V_1)x_\rho,\pi_\rho (V_2)x_\rho}|\pi_\rho (\A))(A)
\end{align*}
for each $A \in \A$ that $\ker \omega_{x,y} \cap \A \in C(\ker \omega_{x_0,x_0} \cap \A)$. Conversely, if $z,w \in \HH_\rho \setminus \{0\}$, then the irreducibility of $\pi_\rho$ ensures that $\pi_\rho (W_1)x_\rho = z$ and $\pi_\rho (W_2)x_\rho = w$, which implies that
\begin{align*}
(\pi_\rho^* (\omega_{z,w}|\pi_\rho (\A)))(A)
= \langle \pi_\rho (A)z,w \rangle
&= \langle \pi_\rho (W_2^*AW_1)z,w \rangle \\
&= \omega_{x_0,x_0}(W_2^*AW_1) \\
&= \omega_{W_1x_0,W_2x_0}(A)
\end{align*}
for each $A \in \A$. In particular, $W_1x_0 \neq 0$ and $W_2x_0 \neq 0$. Therefore, the identity
\[
C(\ker \omega_{x_0,x_0} \cap \A) = \{ \ker \omega_{x,y} \cap \A :x,y \in \HH \setminus \{0\} \}
\]
holds.
\end{proof}
\end{lemma}
The equivalence class of a pure vector state of an irreducible $C^*$-algebra plays important roles in the rest of this paper. The following definition assigns a special symbol to it.
\begin{definition}

Let $\A$ be an irreducible $C^*$-algebra acting on a Hilbert space $\HH$. Then, the equivalence class $\{ \ker \omega_{x,y} \cap \A : x,y \in \HH \setminus \{0\}\}$ in $\s (\A)$ is called the \emph{normal part} and denoted by $\mathfrak{S}_n(\A)$.
\end{definition}
In the rest of this section, we characterize (continuous) finite-homeomorphism between normal parts of geometric structure spaces of irreducible $C^*$-algebras. First, we give a wide class of such mappings. To this end, some preliminary works are needed.
\begin{remark}

Let $\HH$ be a Hilbert space, and let $H \in B(\HH)$ be such that $0 \leq H \leq I$. Then, there exists a sequence of projections $(R_n(H))_n$ with the following properties:
\begin{itemize}
\item[(i)] $HR_n(H) = R_n(H)H \geq n^{-1}R_n(H)$ for each $n$.
\item[(ii)] $R_n(H) \leq R(H)$ for each $n$, where $R(H)$ is the range projection of $H$.
\item[(iii)] The sequence $(R_n(H))_n$ converges to $R(H)$ with respect to the strong-operator topology.
\end{itemize}
This is an elementary application of functional calculus. For each $k \in (0,1)$, let $f_k$ be an element of $C(\Sp (H))$ defined by
\[
f_k(t) = \left\{ \begin{array}{ll}
0 & (t \in [0,k]) \\
(t-k)/(1-k) & (t \in [k,1])
\end{array}
\right. .
\]
Then, we obtain $tf_k(t)^{1/n} \geq kf_k(t)^{1/n}$,
\[
t(1-f_k(t)^{1/n}) \leq \max \{ (1-k)k+k, 1-f_k((1-k)k+k)^{1/n} \}
\]
and $(1-t)f(t)^{1/n} \leq 1-k$ for each $t \in [0,1]$ and each $n$, which implies that $Hf_k(H)^{1/n} \geq kf_k(H)^{1/n}$, $H(I-f_k(H)^{1/n}) \leq ((1-k)k+k)I$ and $(I-H)f_k(H)^{1/n} \leq (1-k)I$ for sufficiently large $n$. Since $(f_k(H)^{1/n})_n$ converges to the range projection $E_k$ of $f_k(H)$ with respect to the strong-operator topology, it follows that $HE_k = E_kH \geq kE_k$, $0 \leq H(I-E_k) \leq ((1-k)k+k)I$ and $0 \leq (I-H)E_k \leq (1-k)I$. In particular, we have $\|H-HE_k\| \leq (1-k)k+k$ and $\|E_k-E_kH\| \leq 1-k$. Here, we note that $E_k \leq R(H)$ for each $k \in (0,1)$ since $f_k(t) \leq t$ for each $t \in [0,1]$ and $(H^{1/n})_n$ converges to $R(H)$ with respect to the strong-operator topology. Moreover, if $x \in R(H)(\HH)$, and if $\varepsilon >0$, then $\|x-Hy\|<\varepsilon /3$ for some $y \in \HH$. It follows from
\begin{align*}
\|E_kx-R(H)x\|
&\leq \|E_kx-E_kHy\| + \|E_kHy-Hy\| + \|Hy-R(H)x\| \\
&\leq \|x-Hy\| +\|E_kH-H\|\|y\| + \|R(H)Hy-R(H)x\| \\
&\leq 2\|Hy-x\| + ((1-k)k+k)\|y\| \\
&< 2\varepsilon /3 + (2k-k^2)\|y\|
\end{align*}
that $\|E_kx-R(H)x\| < \varepsilon$ for sufficiently small $k$. This shows that $\lim_{k \to 0}\|E_kx-R(H)x\| = 0$. The desired sequence is then given by setting $R_n(H) = E_{n^{-1}}$ for each $n$. 
\end{remark}
\begin{lemma}\label{sublemma}

Let $\HH$ be a Hilbert space, and let $A,B \in B(\HH)$. If $\langle Bx,y \rangle = 0$ whenever $\langle Ax,y \rangle = 0$, then $B=cA$ for some $c \in \mathbb{C}$.
\begin{proof}
We only consider the case $A \neq 0$ since $A=0$ implies $B=0$ from the assumption. It may be assumed that $\|A\|=1$. Suppose that $\langle Bx,y \rangle = 0$ whenever $\langle Ax,y \rangle = 0$. Let $A=U|A|=|A^*|U$ be the polar decomposition of $A$, where $U$ is the partial isometry with $U^*U=R(|A|)$ and $UU^*=R(|A^*|)=R(A)$. Since $|A|=|A|R(|A|)$ and $|A^*|=|A^*|R(|A^*|)$, it follows that
\[
\langle A(I-R(|A|))x,y \rangle = \langle U|A|(I-R(|A|))x,y \rangle = 0
\]
and
\[
\langle Ax,(I-R(|A^*|))y \rangle = \langle |A^*|Ux,(I-R(|A^*|))y \rangle = 0
\]
for each $x,y \in \HH$, which together with the assumption implies that
\[
\langle B(I-R(|A|))x,y \rangle = \langle Bx,(I-R(|A^*|))y \rangle = 0 ,
\]
for each $x,y \in \HH$. Hence, $B=R(|A^*|)BR(|A|)$.

Next, recall that $|A|R_n (|A|) = R_n(|A|)|A| \geq n^{-1}R_n(|A|)$, which implies that $|A|R_n (|A|)$ is invertible in $R_n(|A|)B(\HH)R_n(|A|)$. Let $C_n \in R_n(|A|)B(\HH)R_n(|A|)$ be such that $C_n|A|R_n(|A|)=|A|R_n(|A|)C_n=R_n(|A|)$, and let $D_n=C_nU^* = R_n(|A|)C_nR_n(|A|)U^*$. We note that $D_nAR_n(|A|)=R_n(|A|)$ and $AD_n=UR_n(|A|)U^*$. In particular, $AD_n$ is the final projection of the partial isometry $UR_n(|A|)$. Now, suppose that $x,y \in (AD_n)(\HH)$ and $\langle x,y \rangle = 0$. Then, it turns out from
\[
0 = \langle x,y \rangle = \langle AD_nx,y \rangle 
\]
that $\langle BD_nx,y \rangle = 0$. In particular, if $x,y \in (AD_n)(\HH)$, $\|x\|=\|y\|=1$ and $\langle x,y \rangle =0$. Then,
\[
\langle x,y \rangle = \langle y,x \rangle = \langle x+y,x-y \rangle = 0 ,
\]
which implies that
\[
\langle BD_nx,y \rangle = \langle BD_ny,x \rangle = \langle BD_n(x+y),x-y \rangle = 0 .
\]
We obtain $\langle BD_nx,x \rangle = \langle BD_ny,y \rangle$ by
\[
0 = \langle BD_n(x+y),x-y \rangle = \langle BD_nx,x \rangle -\langle BD_ny,y \rangle .
\]
Now, let $(e_\lambda )_\lambda$ be an orthonormal basis for $(AD_n)(\HH)$. Then, we have
\[
\langle (BD_n-\langle BD_ne_{\lambda_0},e_{\lambda_0}\rangle AD_n)e_\lambda ,e_\lambda \rangle = \langle BDe_\lambda ,e_\lambda \rangle -\langle BD_ne_{\lambda_0},e_{\lambda_0} \rangle = 0
\]
for a fixed $\lambda_0$. Moreover, we obtain $\langle BD_ne_\lambda ,e_\mu \rangle = 0$ and
\[
\langle (BD_n-\langle BD_ne_{\lambda_0},e_{\lambda_0}\rangle AD_n)e_\lambda ,e_\mu \rangle = \langle BD_ne_\lambda ,e_\mu \rangle -\langle BD_ne_{\lambda_0},e_{\lambda_0} \rangle \langle e_\lambda,e_\mu \rangle = 0
\]
whenever $\lambda \neq \mu$. Therefore, $AD_nBD_nAD_n = \langle BD_ne_{\lambda_0},e_{\lambda_0}\rangle AD_n$ in $B(\HH)$ and
\begin{align*}
AD_nBR_n(|A|) = (AD_nBD_nAD_n)AR_n(|A|) 
&=  (\langle BD_ne_{\lambda_0},e_{\lambda_0}\rangle AD_n)AR_n(|A|) \\
&= \langle BD_ne_{\lambda_0},e_{\lambda_0}\rangle AD_nAR_n(|A|) .
\end{align*}
Set $c_n = \langle BD_ne_{\lambda_0},e_{\lambda_0}\rangle$ for each $n$. Since $(R_n(|A|))_n$ converges to $R(|A|) = U^*U$ with respect to the strong-operator topology, it follows that $\lim_n \|R_n(|A|)x-R(|A|)x\|=0$ and $\lim_n \|AD_ny-R(|A^*|)y\|=0$ for each $x,y \in \HH$. In particular, if $Ax_0 \neq 0$, then
\begin{align*}
\lim_n \langle AD_nAR_n(|A|)x_0,Ax_0 \rangle 
&= \lim_n \langle AR_n(|A|)x_0,AD_nAx_0 \rangle \\
&= \langle AR(|A|)x_0,R(|A^*|)Ax_0 \rangle \\
&= \|Ax_0\|^2 \\
&> 0
\end{align*}
which together with $\|AD_n\| \leq 1$ and $R_n (|A|) \leq 1$ implies that
\begin{align*}
|c_n | 
&= \left| \frac{\langle AD_nBR_n(|A|)x_0,Ax_0 \rangle}{\langle AD_nAR_n(|A|)x_0,Ax_0 \rangle}\right|
\leq \frac{2\|B\|\|x_0\|}{\|Ax_0\|}
\end{align*}
for sufficiently large $n$. Hence, the sequence $(c_n)_n$ is bounded. Taking a subsequence if necessary, we may assume that $(c_n )_n$ converges to some $c \in \mathbb{C}$. It follows that
\begin{align*}
\langle BR(|A|)x,R(|A^*|)y \rangle 
&= \lim_n \langle BR_n(|A|)x,AD_ny \rangle \\
&= \lim_n c_n \langle AR_n(|A|)x,AD_ny \rangle \\
&= c \langle AR(|A|)x,R(|A^*|)y \rangle
\end{align*}
for each $x,y \in \HH$, that is,
\[
B = R(|A^*|)BR(|A|) = c R(|A^*|)AR(|A|) = c A .
\]
This completes the proof.
\end{proof}
\end{lemma}
\begin{lemma}\label{like-tensor}

Let $\A$ be an irreducible $C^*$-algebra acting on a Hilbert space $\HH$. Suppose that $x_1,\ldots ,x_n ,y_1,\ldots ,y_n \in \HH$. Then, $\sum_{j=1}^n \omega_{x_j,y_j} = 0$ (on $\A$) if and only if there exists an $n\times n$ matrix $C=[c_{jk}]$ such that $\sum_{j=1}^n \overline{c_{jk}}x_j = 0$ for each $k \in \{1,\ldots ,n\}$ and $\sum_{k=1}^n c_{jk}y_k = y_j$ for each $j \in \{1,\ldots ,n\}$.
\begin{proof}
The proof is almost the same as that of \cite[Proposition 2.6.6]{KR97a}. If the matrix $[c_{jk}]$ with the stated properties exists, then it follows that
\begin{align*}
\sum_{j=1}^n \omega_{x_j,y_j} = \sum_{j=1}^n \sum_{k=1}^n \overline{c_{jk}}\omega_{x_j,y_k}
&= \sum_{k=1}^n \sum_{j=1}^n \overline{c_{jk}}\omega_{x_j,y_k} = \sum_{k=1}^n \omega_{0,y_k} = 0 .
\end{align*}

For the converse, suppose that $\sum_{j=1}^n \omega_{x_j,y_j} = 0$ on $\A$. If $y_1=\cdots =y_n=0$, then $C=0$ has the desired properties. We assume that $\langle \{ y_1,\ldots ,y_n \} \rangle \neq \{0\}$. Let $\{e_1,\ldots ,e_r\}$ be an orthonormal basis for $\langle \{ y_1,\ldots ,y_n\} \rangle$. Then, there exist $n \times r$ matrix $A = [a_{jl}]$ and $r \times n$ matrix $B=[b_{lk}]$ such that $y_j = \sum_{i=1}^r a_{ji}e_i$ for each $j \in \{ 1,\ldots ,n\}$ and $e_l = \sum_{k=1}^n b_{lk}y_k$ for each $l \in \{ 1,\ldots ,r\}$. Set $C=AB$. Then, $C$ is an $n \times n$ matrix and
\[
y_j = \sum_{l=1}^r a_{jl} \sum_{k=1}^n b_{lk}y_k = \sum_{k=1}^n \left( \sum_{l=1}^r a_{jl}b_{lk} \right) y_k = \sum_{k=1}^n c_{jk}y_k
\]
for each $j \in \{1,\ldots ,n\}$. Moreover, setting $z_l = \sum_{j=1}^n\overline{a_{jl}}x_j$ for each $l \in \{1,\ldots ,r\}$, we have
\begin{align*}
0=\sum_{j=1}^n \omega_{x_j,y_j} = \sum_{j=1}^n \sum_{l=1}^r \overline{a_{jl}}\omega_{x_j,e_l}
&= \sum_{l=1}^r \sum_{j=1}^n \overline{a_{jl}}\omega_{x_j,e_l}
= \sum_{l=1}^r \omega_{z_l,e_l} .
\end{align*}
Since $\A$ acts irreducibly on $\HH$, there exists an element $U$ of $\A$ such that $Ue_l = z_l$ for each $l \in \{ 1,\ldots ,r\}$. It follows that
\[
0 =  \sum_{l=1}^r \omega_{z_l,e_l}(U^*) = \sum_{l=1}^r \|z_l\|^2 .
\]
Therefore, $z_1=\cdots = z_r = 0$ and
\[
\sum_{j=1}^n \overline{c_{jk}}x_j
= \sum_{j=1}^n \sum_{l=1}^r \overline{a_{jl}b_{lk}}x_j
= \sum_{l=1}^r \overline{b_{lk}} \sum_{j=1}^n \overline{a_{jl}}x_j \\
= \sum_{l=1}^r \overline{b_{lk}}z_l = 0
\]
for each $k \in \{1,\ldots ,n\}$.
\end{proof}
\end{lemma}
\begin{remark}\label{semilinear}

Let $X$ and $Y$ be vector spaces over the scalar field $\mathbb{K}$, respectively. Suppose that $\sigma$ is a field automorphism on $\mathbb{K}$. Then, a mapping $A:X \to Y$ is said to be \emph{$\sigma$-semilinear} if it satisfies $T(x+y)=Tx+Ty$ and $T(cx) = \sigma (c)Tx$ for each $x,y \in X$ and each $c \in \mathbb{K}$. If $\mathbb{K}=\mathbb{R}$, then the identity mapping on $\mathbb{R}$ is a unique field automorphism on $\mathbb{R}$. However, in the complex case, there are many other situations. If $\sigma$ is the complex conjugation on $\mathbb{C}$, then a $\sigma$-similinear mapping are called \emph{conjugate-linear operator}. The identity mapping and complex conjugation on $\mathbb{C}$ are characterized as continuous field automorphism on $\mathbb{C}$.

Let $\HH$ and $\KK$ be Hilbert spaces, and let $A:\HH \to \KK$ be a continuous conjugate-linear. Then, $A:\HH \to \overline{\KK}$ and $A:\overline{\HH} \to \KK$ can be viewed as bounded linear operators, and hence, the adjoint $A^*$ of $A$ exists as a bounded conjugate-linear operator from $\HH \to \KK$.
\end{remark}
Now, we give a wide class of (continuous) finite-homeomorphisms between normal parts of geometric structure spaces of irreducible $C^*$-algebras.
\begin{theorem}\label{continuity}

Let $\A$ and $\B$ be irreducible $C^*$-algebras acting on Hilbert spaces $\HH$ and $\KK$. Suppose that either of the following holds:
\begin{itemize}
\item[{\rm (i)}] $U$ and $V$ are linear bijection from $\HH$ onto $\KK$.
\item[{\rm (ii)}] $U$ and $V$ are conjugate-linear bijection from $\HH$ onto $\KK$.
\end{itemize}
Define a mapping $\Phi:\mathfrak{S}_n(\A) \to \mathfrak{S}_n(\B)$ by
\[
\Phi (\ker \omega_{x,y} \cap \A) = \ker \omega_{Ux,Vy} \cap \B
\]
for each $x,y \in \HH \setminus \{0\}$. Then, $\Phi$ is a finite-homeomorphism. Moreover, provided that $U$ and $V$ are bounded, $\Phi$ is continuous if and only if $V^*\B U \subset \A$, and $\Phi$ is a homeomorphism if and only if $V^*\B U=\A$.
\begin{proof}
Since $U$ and $V$ are semilinear, it is straightforward to check that $\Phi$ is well-defined. Moreover, setting $\Psi (\ker \omega_{z,w} \cap \B) = \ker \omega_{U^{-1}z,V^{-1}w} \cap \A$ for each $z,w \in \KK$ yields the inverse mapping of $\Phi$. In particular, $\Phi$ is bijective. Suppose that $x,x_1,\ldots ,x_n,y,y_1,\ldots ,y_n \in \HH \setminus \{0\}$, and that
\[
\ker \omega_{x,y} \cap \A \in \{ \ker \omega_{x_j,y_j} \cap \A : j \in \{1,\ldots ,n\}\}^= \cap \s_n (\A).
\]
Then, $\bigcap_{j=1}^n (\ker \omega_{x_j,y_j} \cap \A) \subset \ker \omega_{x,y} \cap \A$ and $\omega_{x,y} = \sum_{j=1}^n c_j\omega_{x_j,y_j}$ for some $c_1,\ldots ,c_n \in \mathbb{C}$. Hence, setting $z_j = c_jx_j$ for each $j \in \{1,\ldots ,n\}$, $z_{n+1}=-x$ and $y_{n+1}=y$, we have $\sum_{j=1}^{n+1} \omega_{z_j,y_j}=0$. It follows from Lemma~\ref{like-tensor} that there exists an $(n+1) \times (n+1)$ matrix $[c_{jk}]$ such that $\sum_{j=1}^{n+1} \overline{c_{jk}}z_j = 0$ for each $k \in \{1,\ldots ,n+1\}$ and $\sum_{k=1}^{n+1} c_{jk}y_k = y_j$ for each $j \in \{1,\ldots ,n+1\}$. If $U$ and $V$ are linear, then $\sum_{j=1}^{n+1} \overline{c_{jk}}Uz_j = 0$ and $\sum_{k=1}^{n+1} c_{jk}Vy_k = Vy_j$, while if $U$ and $V$ are conjugate-linear, then $\sum_{j=1}^{n+1} c_{jk}Uz_j = 0$ and $\sum_{k=1}^{n+1} \overline{c_{jk}}Vy_k = Vy_j$. Therefore, in either case, $\sum_{j=1}^{n+1} \omega_{Uz_j,Vy_j} = 0$ by Lemma~\ref{like-tensor}, which implies that $\omega_{Ux,Vy} = \sum_{j=1}^n \omega_{U(c_jx_j),Vy_j}$. This shows that $\bigcap_{j=1}^n (\ker \omega_{Ux_j,Vy_j} \cap \B) \subset \ker \omega_{Ux,Vy} \cap \B$ and
\[
\ker \omega_{Ux,Vy} \cap \B \in \{ \ker \omega_{Ux_j,Vy_j} \cap \B : j \in \{1,\ldots ,n\} \}^= \cap \s_n (\B),
\]
that is, $\Psi$ is finitely continuous. Applying the same argument as above to $\Psi^{-1}$, we can conclude that $\Psi$ is a finite-homeomorphism.

In the rest of this proof, we assume that $U$ and $V$ are bounded. We note that $U^{-1}$ and $V^{-1}$ are also bounded by the open mapping theorem. Suppose that $\Phi$ is continuous, and that $B \in \B \setminus \{0\}$. Set
\[
S(B) = \{ \ker \omega_{z,w} \cap \B \in \mathfrak{S}(\B) : B \in \ker \omega_{z,w} \cap B \} \subset \mathfrak{S}_n(\B).
\]
Then, $S(B)^= \cap \mathfrak{S}_n(\B) = S(B)$ and $S(B) \neq \mathfrak{S}_n(\B)$. Indeed, if $\ker \omega_{z_0,w_0} \cap \B \in S(B)^= \cap \mathfrak{S}_n(\B)$, then
\[
B \in \bigcap \{ \ker \omega_{z,w} \cap \B \in \mathfrak{S}(\B) : B \in \ker \omega_{z,w} \cap \B \} \subset \ker \omega_{z_0,w_0} \cap \B ,
\]
that is, $\ker \omega_{z_0,w_0} \cap \B \in S(B)$. Moreover, we realize that $S(B) \neq \mathfrak{S}_n(\B)$ by $\ker \omega_{z,Bz} \cap \B \in \mathfrak{S}_n(\B) \setminus S(B)$ whenever $Bz \neq 0$. Now, by the continuity of $\Phi$, we obtain
\[
\Phi ( S^= \cap \mathfrak{S}_n(\A)) \subset \Phi (S)^= \cap \mathfrak{S}_n (\B)
\]
for each $S \subset \mathfrak{S}_n (\A)$. In particular, $\Phi^{-1}(S(B))^= \cap \mathfrak{S}_n (\A) = \Phi^{-1}(S(B)) \neq \mathfrak{S}_n(\A)$. It follows that $\Phi^{-1}(S(B))^= \neq \mathfrak{S}(\A)$, which together with \cite[Lemma 4.1]{Tan23a} implies that there exists an $A \in \A \setminus \{0\}$ such that
\[
A \in \bigcap \{ \Phi^{-1}(\ker \omega_{z,w} \cap \B) : z,w \in \KK \setminus \{0\},~\langle Bz,w \rangle = 0\} .
\]
Since $\Phi^{-1} = \Psi$, we have $\langle AU^{-1}z,V^{-1}w \rangle = 0$ whenever $\langle Bz,w \rangle = 0$. We note that $(V^{-1})^*AU^{-1} \in B(\KK)$ even if $U$ and $V$ are both conjugate-linear. Therefore, Lemma~\ref{sublemma} ensures that $(V^{-1})^*AU^{-1} = cB$ for some $c \in \mathbb{C} \setminus \{0\}$, that is, $V^*BU = c^{-1}A \in \A$. This shows that $V^*\B U \subset \A$.

Conversely, suppose that $V^*\B U \subset \A$. Let $S \subset \mathfrak{S}_n(\A)$, and let $\ker \omega_{x_0,y_0} \cap \A \in S^= \cap \mathfrak{S}_n(\A)$. If
\[
B \in \bigcap \{ \Phi (\ker \omega_{x,y} \cap \A) : \ker \omega_{x,y} \cap \A \in S\}
\]
and $\ker \omega_{x,y} \cap \A \in S$, then $\langle V^*BUx,y \rangle = \langle BUx,Vy \rangle = 0$, which together with $V^*BU \in V^*\B U \subset \A$ implies that $V^*BU \in \bigcap \{ \ker \omega_{x,y} \cap \A : \ker \omega_{x,y} \cap \A \in S\}$. Since $\ker \omega_{x_0,y_0} \cap \A \in S^=$, it follows that $\langle BUx_0,Vy_0 \rangle = \langle V^*BUx_0,y_0 \rangle = 0$ and $B \in \ker \omega_{Ux_0,Vy_0} \cap \B = \Phi (\ker \omega_{x_0,y_0} \cap \A)$. Hence,
\[
\bigcap \{ \Phi (\ker \omega_{x,y} \cap \A) : \ker \omega_{x,y} \cap \A \in S\} \subset \Phi (\ker \omega_{x_0,y_0} \cap \A ) ,
\]
that is, $\Phi (\ker \omega_{x_0,y_0} \cap \A) \in \Phi (S)^= \cap \mathfrak{S}_n (\B)$. This proves that $\Phi (S^= \cap \mathfrak{S}_n(\A)) \subset \Phi (S)^= \cap \mathfrak{S}_n(\B)$. Thus, $\Phi$ is continuous.

Finally, we note that $\Phi$ is a homeomorphism if and only if both $\Phi$ and $\Psi$ are continuous, which occurs if and only if $V^*\B U \subset \A$ and $(V^{-1})^*\A U^{-1} \subset \B$, that is, $V^*\B U = \A$.
\end{proof}
\end{theorem}
The preceding theorem gives typical examples of (continuous) finite-homeomorphisms between $\s_n (\A)$ and $\s_n (\B)$ for irreducible $C^*$-algebras $\A$ and $\B$. Our next task is to show that continuous finite-homeomorphisms between $\s_n (\A)$ and $\s_n (\B)$ are limited to typical ones at least in the infinite-dimensional case. To this end, we consider the relation ``$\rightleftarrows$'' in a local setting. Then, as will be shown in the following lemma, it is well-behaved in the context of normal parts of geometric structure spaces of irreducible $C^*$-algebras.
\begin{lemma}\label{dependent}

Let $\A$ be an irreducible $C^*$-algebra acting on a Hilbert space $\HH$ that acts irreducibly on $\HH$. Suppose that $x_1,x_2,y_1,y_2 \in \HH \setminus \{0\}$. Then, $\ker \omega_{x_1,y_1} \cap \A \rightleftarrows \ker \omega_{x_2,y_2} \cap \A$ in $\s_n (\A)$ if and only if either $\{x_1,x_2\}$ or $\{y_1,y_2\}$ is linearly dependent.
\begin{proof}
First, we note that $\ker \omega_{x_1,y_1} \cap \A = \ker \omega_{x_2,y_2} \cap \A$ if and only if both $\{x_1,x_2\}$ and $\{y_1,y_2\}$ are linearly dependent. Indeed, if $\{x_1,x_2\}$ is linearly independent, then there exists an $A \in \A$ such that $Ax_1=0$ and $Ax_2=y_2$, in which case, $\omega_{x_1,y_1}(A) = 0$ and $\omega_{x_2,y_2}(A) = 1$. Therefore, $\ker \omega_{x_1,y_1} \cap \A \neq \ker \omega_{x_2,y_2} \cap \A$. An argument similar to above also works in the case that $\{y_1,y_2\}$ is linearly independent. Conversely, if both $\{x_1,x_2\}$ and $\{y_1,y_2\}$ are linearly dependent, then it is obvious that $\ker \omega_{x_1,y_1} \cap \A = \ker \omega_{x_2,y_2} \cap \A$.

Now, suppose that $\{x_1,x_2\}$ is linearly dependent and $\{y_1,y_2\}$ is linearly independent. Then, $\ker \omega_{x_2,y_2} \cap \A = \ker \omega_{x_1,y_2} \cap \A$, and
\[
\ker \omega_{x_1,y} \cap \A \in \{ \ker \omega_{x_1,y_1} \cap \A ,\ker \omega_{x_1,y_2} \cap \A \}^= \cap \mathfrak{S}_n(\A)
\]
whenever $y \in \langle \{y_1,y_2\}\rangle \setminus \{0\}$. Hence,
\[
\{ \ker \omega_{x_1,y_1} \cap \A ,\ker \omega_{x_2,y_2} \cap \A \}^= \cap \mathfrak{S}_n (\A) \neq \{ \ker \omega_{x_1,y_1} \cap \A ,\ker \omega_{x_2,y_2} \cap \A \} ,
\]
that is, $\ker \omega_{x_1,y_1} \cap \A  \rightleftarrows \ker \omega_{x_2,y_2} \cap \A$ in $\s_n (\A)$. Similarly, we have the same consequence as above in the other case.

For the converse, suppose that both $\{x_1,x_2\}$ and $\{y_1,y_2\}$ are linearly independent. In this case, $\ker \omega_{x_1,y_1} \cap \A \neq \ker \omega_{x_2,y_2} \cap \A$. It may be assumed that $x_1,x_2,y_1,y_2 \in S_\HH$. Let $x,y \in \HH \setminus \{0\}$ be such that
\[
\ker \omega_{x,y} \cap \A \in \{ \ker \omega_{x_1,y_1} \cap \A ,\ker \omega_{x_1,y_1} \cap \A \}^= \cap \mathfrak{S}_n(\A).
\]
Then, we have $\omega_{x,y}|\A = (c_1\omega_{x_1,y_1}+ c_2\omega_{x_2,y_2})|\A$ for some scalars $c_1,c_2$. If $\{x_1,x_2,x\}$ is linearly independent, then there exists an element $B \in \A$ such that $Bx_j=0$ for $j=1,2$ and $Bx = y$. It follows that
\[
\|y\|^2 = \omega_{x,y}(B) = c_1\omega_{x_1,y_1}(B)+c_2\omega_{x_2,y_2}(B) = 0 ,
\]
which contradicts $y \neq 0$. Hence, $\{x_1,x_2,x\}$ is linearly dependent. By an argument similar to above, we also derive that $\{y_1,y_2,y\}$ is linearly dependent. Therefore, we have $x = a_1x_1+a_2x_2$ and $y = b_1y_1+b_2y_2$ for some scalars $a_1,a_2,b_1,b_2$. It follows that
\begin{align*}
(c_1\omega_{x_1,y_1}+c_2\omega_{x_2,y_2})|\A
&= \omega_{x,y}|\A \\
&= (a_1b_1 \omega_{x_1,y_1} + a_1b_2 \omega_{x_1,y_2} + a_2b_1\omega_{x_2,y_1} + a_2b_2 \omega_{x_2,y_2} )|\A .
\end{align*}
Meanwhile, since $\{x_1,x_2\}$ is linearly independent, there exist $A_1,A_2 \in \A$ such that $A_1x_1=0$, $A_1x_2=y_1-\langle y_1,y_2\rangle y_2$, $A_2x_1=y_2-\langle y_2,y_1 \rangle y_1$ and $A_2x_2=0$. We note that
\[
\omega_{x_1,y_1}(A_1) = \omega_{x_1,y_2}(A_1) = \omega_{x_2,y_2}(A_1) = \omega_{x_1,y_1}(A_2) = \omega_{x_2,y_1}(A_2) = \omega_{x_2,y_2}(A_2) = 0,
\]
while
\[
\omega_{x_2,y_1}(A_1) = \omega_{x_1,y_2}(A_2) = 1-|\langle y_1,y_2 \rangle |^2 >0
\]
since $\{y_1,y_2\}$ is linearly independent. Therefore, $a_2b_1= a_1b_2=0$ by
\begin{align*}
0 
&= (a_1b_1 \omega_{x_1,y_1} + a_1b_2 \omega_{x_1,y_2} + a_2b_1\omega_{x_2,y_1} + a_2b_2 \omega_{x_2,y_2} )(A_1) \\
&= (a_1b_1 \omega_{x_1,y_1} + a_1b_2 \omega_{x_1,y_2} + a_2b_1\omega_{x_2,y_1} + a_2b_2 \omega_{x_2,y_2} )(A_2) .
\end{align*}
Combining this with $x,y \in \HH \setminus \{0\}$, we obtain $(x,y) = (a_1x_1,b_1y_1)$ or $(x,y) = (a_2x_2,b_2y_2)$, which implies that
\[
\ker \omega_{x,y} \cap \A \in \{ \ker \omega_{x_1,y_1} \cap \A ,\ker \omega_{x_1,y_1} \cap \A \} .
\]
Thus,
\[
\{ \ker \omega_{x_1,y_1} \cap \A ,\ker \omega_{x_1,y_1} \cap \A \}^= \cap \mathfrak{S}_n (\A)= \{ \ker \omega_{x_1,y_1} \cap \A ,\ker \omega_{x_1,y_1} \cap \A \} .
\]
This shows that $\ker \omega_{x_1,y_1} \cap \A \not \rightleftarrows \ker \omega_{x_2,y_2} \cap \A$ in $\s_n (\A)$.
\end{proof}
\end{lemma}
In the rest of this section, certain subsets of $\s_n (\A)$ play important roles. We formulate them as follows.
\begin{definition}

Let $(K,c)$ be a close space, and let $S \subset K$. Then, $S$ is called an \emph{$R$-set} if $t_1 \rightleftarrows t_2$ whenever $t_1,t_2 \in S$.
\end{definition}
It can be shown that $R$-sets are preserved under finite-homeomorphisms. 

\begin{lemma}\label{R-sets-preserved}

Let $(K,c)$ and $(L,d)$ be closure spaces, and let $f:K \to L$ be a finite-homeomorphism. Suppose that $S \subset K$. Then, $S$ is an $R$-set if and only if $f(S)$ is.
\begin{proof}
By Lemma~\ref{finite-conti}, for $t_1,t_2 \in K$, we have $t_1\rightleftarrows t_2$ if and only if $f(t_1)\rightleftarrows f(t_2)$. If $S$ is an $R$-set, and if $t_1,t_2 \in S$, then $t_1 \rightleftarrows t_2$, which implies that $f(t_1) \rightleftarrows f(t_2)$. This shows that $f(S)$ is also an $R$-set. Applying the same argument as above to $f^{-1}$, we infer that $S$ is an $R$-set if $f(S)$ is.
\end{proof}
\end{lemma}
In the following lemma, for an irreducible $C^*$-algebra $\A$, we identify $R$-sets in $\s_n (\A)$ maximal with respect to inclusion.
\begin{lemma}\label{R-rep}

Let $\A$ be an irreducible $C^*$-algebra acting on a Hilbert space $\HH$. Suppose that $S \subset \mathfrak{S}_n(\A)$. Then, $S$ is an $R$-set in $\s_n (\A)$ maximal with respect to inclusion if and only if either of the following holds:
\begin{itemize}
\item[{\rm (i)}] There exists an $x_0 \in \HH \setminus \{0\}$ such that $S=\{ \ker \omega_{x_0,y}\cap \A : y \in \HH \setminus \{0\} \}$.
\item[{\rm (ii)}] There exists a $y_0 \in \HH \setminus \{0\}$ such that $S=\{ \ker \omega_{x,y_0}\cap \A : x \in \HH \setminus \{0\} \}$.
\end{itemize}
\begin{proof}
The lemma is obvious if $\dim \HH =1$. We assume that $\dim \HH \geq 2$. Suppose that $S$ is an $R$-set in $\s_n (\A)$ maximal with respect to inclusion. Then, $S$ contains at least two distinct elements. Indeed since each singleton is an $R$-set, it follows that $S \neq \emptyset$. Let $\ker \omega_{x_0,y_0} \cap \A \in S$, and let $y_0^\perp \in \HH \setminus \{0\}$ be such that $\langle y_0,y_0^\perp \rangle = 0$. Then, $\ker \omega_{x_0,y_0} \cap \A \neq \ker \omega_{x_0,y_0^\perp} \cap \A$ and $\{ \ker \omega_{x_0,y_0} \cap \A ,\ker \omega_{x_0,y_0^\perp} \cap \A \}$ is an $R$-set in $\s_n (\A)$ by Lemma~\ref{dependent}. Hence, $S \neq \{ \ker \omega_{x_0,y_0} \cap \A \}$ by the maximality of $S$. Let $x'_0,y'_0 \in \HH \setminus \{0\}$ be such that $\ker \omega_{x'_0,y'_0} \cap \A \in S$ and $\ker \omega_{x_0,y_0} \cap \A \neq \ker \omega_{x'_0,y'_0} \cap \A$. By Lemma~\ref{dependent}, either $\{x_0,x'_0\}$ or $\{y_0,y'_0\}$ is linearly dependent.

Suppose that $\{ x_0,x'_0\}$ is linearly dependent. Then, $\ker \omega_{x'_0,y'_0} \cap \A = \ker \omega_{x_0,y'_0} \cap \A$ and $\{y_0,y'_0\}$ is linearly independent by $\ker \omega_{x_0,y_0} \cap \A \neq \ker \omega_{x'_0,y'_0} \cap \A$. Let $x,y \in \HH \setminus \{0\}$ be such that $\ker \omega_{x,y} \cap \A \in S$. If $\{x_0,x\}$ is linearly independent, then $\{y_0,y\}$ is linearly dependent by Lemma~\ref{dependent}, which implies that $\ker \omega_{x,y} \cap \A = \ker \omega_{x,y_0}\cap \A$. However, then $\ker \omega_{x_0,y'_0} \cap \A \not \rightleftarrows \ker \omega_{x,y_0} \cap \A$ in $\s_n (\A)$ again by Lemma~\ref{dependent}. This contradicts the assumption that $S$ is an $R$-set. Therefore, $\{x_0,x\}$ must be linearly dependent, that is, $\ker \omega_{x,y} \cap \A = \ker \omega_{x_0,y} \cap \A$. This shows that $S \subset \{ \ker \omega_{x_0,y} \cap \A : y \in \HH \setminus \{0\} \}$. Since $\{ \ker \omega_{x_0,y} \cap \A : y \in S_{\HH}\}$ is itself an $R$-set in $\s_n (\A)$ by Lemma~\ref{dependent}, we obtain $S=\{ \ker \omega_{x_0,y} \cap \A : y \in \HH \setminus \{0\} \}$ by the maximality of $S$.

If $\{y_0,y_0'\}$ is linearly dependent, then an argument similar to that in the preceding paragraph shows that $S = \{ \ker \omega_{x,y_0} \cap \A : x \in \HH \setminus \{0\} \}$.

Conversely, let $S = \{ \ker \omega_{x_0,y}\cap \A : y \in \HH \setminus \{0\} \}$, and let $S'$ be an $R$-set in $\s_n (\A)$ that contains $S$. Take an arbitrary pair $y_1,y_2 \in \HH \setminus \{0\}$ with $\langle y_1,y_2 \rangle = 0$. Then, $\ker\omega_{x_0,y_j} \cap \A \in S \subset S'$ for $j=1,2$. If $x,y \in \HH \setminus \{0\}$ satisfies $\ker \omega_{x,y} \cap \A \in S'$, then $\ker \omega_{x,y} \cap \A \rightleftarrows \ker \omega_{x_0,y_j} \cap \A$ in $\s_n (\A)$ for $j=1,2$. If $\{ x_0,x\}$ is linearly independent, then Lemma~\ref{dependent} ensures that $\{ y,y_j \}$ is linearly dependent for $j=1,2$. However, this is impossible since $\{y_1,y_2\}$ is linearly independent. Thus, $\{ x_0,x\}$ must be linearly dependent, in which case, $\ker \omega_{x,y} \cap \A = \ker \omega_{x_0,y} \cap \A \in S$. This proves the maximality of $S$. An argument similar to that in above also works in the other case.
\end{proof}
\end{lemma}
For each nonzero elements $x$ of a Banach space, let $\overline{x}$ denote the one-dimensional subspace generated by $x$. Let $\HH$ be a Hilbert space, and let $x,y \in \HH \setminus \{0\}$. Then, we note that $\ker \omega_{x,y} \cap \A = \ker \omega_{x',y'} \cap \A$ whenever $x' \in \overline{x} \setminus \{0\}$ and $y' \in \overline{y} \setminus \{0\}$. Hence, we have
\[
\{ \ker \omega_{x,y}\cap \A : y \in \HH \setminus \{0\} \} = \{ \ker \omega_{x',y}\cap \A : y \in \HH \setminus \{0\} \}
\]
whenever $\overline{x}=\overline{x'}$, and
\[
\{ \ker \omega_{x,y}\cap \A : y \in \HH \setminus \{0\} \} = \{ \ker \omega_{x,y'}\cap \A : y \in \HH \setminus \{0\} \}
\]
whenever $\overline{y}=\overline{y'}$. Based on this observation, we assign special symbols to these sets.
\begin{definition}

Let $\A$ be an irreducible $C^*$-algebra acting on a Hilbert space $\HH$. Then, let
\begin{align*}
R(\overline{x},\HH) &= \{ \ker \omega_{x,y}\cap \A : y \in \HH \setminus \{0\} \} \\
R(\HH ,\overline{y}) &= \{ \ker \omega_{x,y}\cap \A : x \in \HH \setminus \{0\} \} 
\end{align*}
for each $x,y \in \HH \setminus \{0\}$. 
\end{definition}
We need a partial extension of Lemma~\ref{sublemma}.

\begin{lemma}\label{sublemma-semi}

Let $\HH$ be a Hilbert space, and let $A \in B(\HH)$ be a positive operator with $\dim A(\HH) \geq 2$. Suppose that $B:\HH \to \HH$ is a bounded linear or conjugate linear operator such that $\langle Bx,y \rangle = 0$ whenever $\langle Ax,y \rangle = 0$ and $\langle Bx_0,x_0\rangle \neq 0$ for some $x_0 \in \HH \setminus \{0\}$. Then, $B = cA$ for some $c \in \mathbb{C} \setminus \{0\}$.
\begin{proof}
By Lemma~\ref{sublemma}, it is sufficient to show that $B$ is linear. We note that $\langle Ax_0,x_0 \rangle >0$ from the assumption. Since $A(\HH) \subset A^{1/2}(\HH)$, there exists an element $y_0 \in \HH \setminus \{0\}$ such that $A^{1/2}y_0 \neq 0$ and $\langle A^{1/2}x_0,A^{1/2}y_0 \rangle = 0$. Namely, we have $x_0,y_0 \in \HH \setminus \{0\}$ such that $\langle Ax_0,x_0\rangle \langle Ay_0,y_0 \rangle >0$, $\langle Ax_0,y_0 \rangle = 0$ and $\langle Bx_0,x_0 \rangle \neq 0$.

Let $z \in \mathbb{C} \setminus \{0\}$ and
\[
c = \sqrt{\frac{\langle Ax_0,x_0\rangle}{\langle Ay_0,y_0\rangle}} >0.
\]
We note that $\langle Ax_0,y_0 \rangle = \langle Ay_0,x_0 \rangle = 0$, which implies that
\[
\langle A(x_0+zcy_0),x_0-\overline{z^{-1}}cy_0 \rangle = \langle Ax_0,x_0 \rangle - c^2\langle Ay_0,y_0 \rangle = 0 .
\]
In particular, $\langle A(x_0+cy_0),x_0-cy_0 \rangle = 0$. Combining this with $\langle A(cy_0),x_0 \rangle = \langle A(zcy_0),x_0 \rangle = 0$, we obtain $\langle Bx_0,y_0 \rangle = \langle B(cy_0),x_0 \rangle = \langle B(zcy_0),x_0 \rangle = 0$ and
\[
\langle B(x_0+cy_0),x_0-cy_0 \rangle = \langle B(x_0+zcy_0),x_0-\overline{z^{-1}}cy_0 \rangle = 0 .
\]
It follows that
\[
0 = \langle B(x_0+cy_0),x_0-cy_0 \rangle = \langle Bx_0,x_0 \rangle - \langle B(cy_0),cy_0 \rangle
\]
and
\[
0 = \langle B(x_0+zcy_0),x_0-\overline{z^{-1}}cy_0 \rangle = \langle Bx_0,x_0 \rangle -  \langle B(zcy_0),\overline{z^{-1}}cy_0 \rangle ,
\]
that is,
\[
\langle B(cy_0),cy_0 \rangle = \langle Bx_0,x_0 \rangle = \langle B(zcy_0),\overline{z^{-1}}cy_0 \rangle = \sigma (z)z^{-1} \langle B(cy_0),cy_0 \rangle .
\]
We note that $\langle B(cy_0),cy_0 \rangle = \langle Bx_0,x_0 \rangle \neq 0$. Therefore, $\sigma (z) = z$ for each $z \in \mathbb{C}$. This proves that $B$ is linear.
\end{proof}
\end{lemma}
We present a useful tool for simplifying our argument on (finite-)homeomorphisms between (normal parts of) geometric structure spaces of irreducible $C^*$-algebras.
\begin{lemma}\label{funct-adjoint}

Let $\A$ be a $C^*$-algebra, and let $\Phi^{\A}_a (\ker \rho ) = \ker \rho^*$ for each $\ker \rho \in \mathfrak{S}(\A)$. Then, $\Phi^{\A}_a$ is a homeomorphism on $\mathfrak{S}(\A)$ with $(\Phi^{\A}_a)^{-1}=\Phi^{\A}_a$.
\begin{proof}
Since the mapping $\rho \mapsto \rho^*$ is a conjugate-linear isometric isomorphism on $\A^*$, it follows that $\ker \rho^* \in \mathfrak{S}(\A)$ whenever $\ker \rho \in \mathfrak{S}(\A)$. Moreover, we have
\[
\ker \rho^* = \{ A \in \A : \rho (A^*) = 0 \} = \{ A^* : A \in \ker \rho \}
\]
for each $\rho \in \A$, the mapping $\Phi^{\A}_a:\mathfrak{S}(\A) \to \mathfrak{S}(\A)$ is well-defined. Since $\Phi^{\A}_a \circ \Phi^{\A}_a$ is the identity mapping on $\s (\A)$, the mapping $\Phi^{\A}_a$ is bijective and $(\Phi^{\A}_a)^{-1} = \Phi^{\A}_a$. Now, let $S \subset \s (\A)$. Then, we have
\[
\bigcap \{ \ker \rho^* : \ker \rho \in S \} = \left\{ A^* : A \in \bigcap \{ \ker \rho : \ker \rho \in S \} \right\} ,
\]
which implies that $\Phi^{\A}_a (S^=) \subset \Phi^{\A}_a (S)^=$. This shows that $\Phi^{\A}_a$ (and hence, its inverse) is continuous.
\end{proof}
\end{lemma}
The final ingredient needed for our argument is the following result of Kakutani and Mackey~\cite[Lemma 2 and Corollary]{KM46}. We give its proof for the readers' convenience; see also~\cite[Lemma 2]{FL84}.
\begin{lemma}[Kakutani and Mackey, 1946]\label{auto-conti}

Let $X,Y$ be infinite dimensional complex Banach spaces, and let $T:X \to Y$ be a semilinear bijection. If $T(M)$ is closed whenever $M$ is a closed hyperplane of $X$, then $T$ is linear or conjugate linear.
\begin{proof}
Suppose that $T$ is $\sigma$-linear. Since $\sigma (1)=1$ and $\sigma (i)\in \{-1,1\}$, we have $\sigma (p+iq) = p+iq$ for each $p,q \in \mathbb{Q}$, or $\sigma (p+iq) = p-iq$ for each $p,q \in \mathbb{Q}$. Hence, it is sufficient to show that $\sigma$ is continuous. To this end, suppose to the contrary that $\sigma$ is discontinuous. Since $\sigma$ is additive, there exists a sequence $(c_n)_n \subset \mathbb{C}$ such that $\lim_n c_n = 0$, and $(\sigma (c_n))_n$ does not converge to $0$. We may assume that $|\sigma (c_n)| \geq \varepsilon >0$ for each $n$. Let $q_n \in \mathbb{Q}$ be such that $q_n> 0$, and $| q_n^{-1}|c_n|-1 | < n^{-1}$. Then, $\lim_n q_n^{-1}|c_n| = 1$, $\lim_n q_n = 0$, and $\lim_n |\sigma (q_n^{-1}c_n)| = \lim_n q_n^{-1}|\sigma (c_n)|=+\infty$. Hence, we have a bounded sequence $(d_n)_n \subset \mathbb{C} \setminus \{0\}$ such that $\lim_n |\sigma (d_n)| =+\infty$.

Now, let $(x_n)_n$ be a basic sequence in $X$, and let $(f_n)_n$ be the coordinate functionals of $(x_n)_n$. By the Hahn-Banach theorem, each $f_n \in [(x_n)_n]^*$ extends to an element of $X^*$. Replacing $(x_n,f_n)$ with $(2^n\|f_n\|x_n,2^{-n}\|f_n\|^{-1}f_n)$ if necessary, we have a biorthogonal sequence $((x_n,f_n))_n \subset X \times X^*$ satisfying $\|f_n\|=2^{-n}$. Moreover, taking a subsequence of $(d_n)_n$, it may be assumed that $|\sigma (d_n)|>n\|Tx_n\|$ for each $n$. Since $(d_n)_n$ is bounded, we can define an element of $X^*$ by $f=\sum_n d_n f_n$. It follows that $f(x_n) = d_n$ for each $n$. From this, we know that $f\neq 0$, and $fx_0=1$ for some $x_0 \in X$. Set $y_n = x_0-d_n^{-1}x_n$ for each $n$. Then, $(y_n)_n \subset \ker f$, and 
\[
\|Tx_0-Ty_n\| = \|T(x_0-y_n)\| = \|T(d_n^{-1}x_n)\| = \frac{\|Tx_n\|}{|\sigma (d_n)|}<\frac{1}{n} \to 0 ,
\]
which implies that $Tx_0 \in \overline{T(\ker f)}$. Meanwhile, from the assumption, $T(\ker f)$ is closed as $\ker f$ is a closed hyperplane of $X$. Therefore, $Tx_0 \in T(\ker f)$. However, this leads to $x_0 \in \ker f$ since $T$ is bijective, which contradicts $f(x_0) = 1$. Thus, $\sigma$ is continuous.
\end{proof}
\end{lemma}
Now, we give the main result in this section which describes the form of (continuous) finite-homeomorphisms between $\s_n (\A)$ and $\s_n (\B)$, where $\A$ and $\B$ are irreducible $C^*$-algebras.
\begin{theorem}\label{weak-preserver}

Let $\A$ and $\B$ be irreducible $C^*$-algebras acting on Hilbert spaces $\HH$ and $\KK$, respectively. Suppose that $\Phi : \mathfrak{S}_n(\A) \to \mathfrak{S}_n(\B)$ is a finite-homeomorphism. Then, $\dim \HH = \dim \KK$ provided that either $\HH$ or $\KK$ is finite-dimensional. Moreover, if $\dim \HH =\dim \KK \geq 3$, then there exists a pair of semilinear bijections $U$ and $V$ from $\HH$ onto $\KK$ such that either of the following holds:
\begin{itemize}
\item[{\rm (i)}] $\Phi (\ker \omega_{x,y} \cap \A ) = \ker \omega_{Ux,Vy} \cap \B$ for each $x,y \in \HH \setminus \{0\}$.
\item[{\rm (ii)}] $\Phi (\ker \omega_{x,y} \cap \A ) = \ker \omega_{Vy,Ux} \cap \B$ for each $x,y \in \HH \setminus \{0\}$.
\end{itemize}
If further $\HH$ and $\KK$ are infinite-dimensional and $\Phi$ is continuous, then there exist bounded linear or conjugate-linear isomorphisms $U$ and $V$ that satisfy either (i) or (ii), in which case, $U$ is linear if and only if $V$ is, and $V^*\B U \subset \A$. In particular, $\dim \HH = \dim \KK$. If $\Phi$ is a homeomorphism, then $V^*\B U = \A$.
\begin{proof}
If $\dim \HH = 1$ or $\dim \KK = 1$, then $\A = \mathbb{C}I$ or $\B = \mathbb{C}I$, which implies that
\[
\card (\s(\A)) = \card (\s_n (\A)) = \card (\s_n (\B)) = \card (\s (\B)) = 1.
\]
Hence, $\dim \A = \dim \B = 1$ by \cite[Propositions 4.4 and 4.5]{Tan23a}, in which case, $\dim \HH = \dim \KK = 1$. In the following, we assume that $\dim \HH \geq 2$ and $\dim \KK \geq 2$.

Let $\mathfrak{M}_R (\A)$ and $\mathfrak{M}_R (\B)$ be the collections of all $R$-sets in $\s_n (\A)$ and $\s_n (\B)$ that are maximal with respect to inclusion, respectively. Then, Lemma~\ref{R-rep} ensures that
\begin{align*}
\mathfrak{M}_R(\A) &= \{ R(\overline{x},\HH) ,R(\HH,\overline{y}) : x,y \in \HH \setminus \{0\} \} \\
\mathfrak{M}_R(\B) &= \{ R(\overline{z},\KK) ,R(\KK,\overline{w}) : z,w \in \KK \setminus \{0\} \} .
\end{align*}
Since $\Phi$ is a finite-homeomorphism, it preserves $R$-sets by Lemma~\ref{R-sets-preserved}. Hence, $S \in \mathfrak{M}_R(\A)$ if and only if $\Phi (A) \in \mathfrak{M}_R(\B)$. Fix an arbitrary $x_0 \in \HH \setminus \{0\}$. Then, $\Phi (R(\overline{x_0},\HH)) \in \mathfrak{M}_R(\B)$. If $\Phi (R(\overline{x_0},\HH)) = R(\KK ,\overline{w_0})$ for some $w_0 \in \KK \setminus \{0\}$, we consider $\Phi^{\mathfrak{B}}_a \circ \Phi$ instead of $\Phi$, where $\Phi^{\B}_a$ is a homeomorphism on $\mathfrak{S}(\B)$ defined in Lemma~\ref{funct-adjoint}. We note that $\omega_{z,w}^* = \omega_{w,z}$ for each $z,w \in \KK$, and that $\Phi^{\B}_a (R(\KK,\overline{w_0})) = R(\overline{w_0},\KK)$. It should also be noted that $\Phi^{\B}_a$ maps $\mathfrak{S}_n(\B)$ onto itself. Therefore, it is sufficient to prove the existence of mappings $U$ and $V$ with desired properties under the additional assumption that $\Phi (R(\overline{x_0},\HH)) = R(\overline{z_0},\KK)$ for some $z_0 \in \KK \setminus \{0\}$.

Let $x,y \in \HH \setminus \{0\}$. We claim that there exist $z,w \in \KK \setminus \{0\}$ such that $\Phi (R(\overline{x},\HH)) = R(\overline{z},\KK)$ and $\Phi (R(\HH,\overline{y})) = R(\KK,\overline{w})$. To show this, we note that if $R(\overline{x_0},\HH) \cap R(\overline{x},\HH)\neq \emptyset$ then $\overline{x_0}=\overline{x}$. Indeed, if $\ker \omega_{x',y'} \cap \A \in R(\overline{x_0},\HH) \cap R(\overline{x},\HH)$, then $x' \in \overline{x_0} \cap \overline{x}$, which implies that $\overline{x_0}=\overline{x}$. Now, suppose that $\Phi (R(\overline{x},\HH)) = R(\KK ,\overline{w})$ for some $w \in \KK$. Then, it follows from
\begin{align*}
\Phi (R(\overline{x_0},\HH) \cap R(\overline{x},\HH)) 
&= \Phi (R(\overline{x_0},\HH)) \cap \Phi (R(\overline{x},\HH) \\
&= R(\overline{z_0},\KK) \cap R(\KK,\overline{w}) \\
&= \{ \ker \omega_{z_0,w} \cap \B \}
\end{align*}
that $\overline{x_0}=\overline{x}$. However, this means that $R(\overline{x_0},\HH)=R(\overline{x},\HH)$ and
\[
\Phi (R(\overline{x_0},\HH )) = \{ \ker \omega_{z_0,w} \cap \B\} ,
\]
which is impossible since $R(\overline{x_0},\HH)$ contains infinitely many elements by $\dim \HH \geq 2$. Therefore, $\Phi (R(\overline{x},\HH)) = R(\overline{z},\KK)$ for some $z \in \KK \setminus \{0\}$. Similarly since
\begin{align*}
R(\overline{z_0},\KK) \cap \Phi (R(\HH,\overline{y}))
&= \Phi (R(\overline{x_0},\HH)) \cap \Phi(R(\HH,\overline{y})) \\
&= \Phi (R(\overline{x_0},\HH) \cap R(\HH,\overline{y})) \\
&= \{ \Phi (\ker \omega_{x_0,y} \cap \A) \} ,
\end{align*}
it follows that $\Phi (R(\HH,\overline{y})) = R(\KK,\overline{w})$; for otherwise, $R(\overline{z_0},\KK) = \{ \Phi (\omega_{x_0,y} \cap \A) \}$, which is also impossible since $\dim \KK \geq 2$.

For each $x,y \in \HH \setminus \{0\}$, let $\varphi (\overline{x})$ and $\psi (\overline{y})$ be the unique one-dimensional subspaces of $\KK$ that satisfy $\Phi (R(\overline{x},\HH)) = R(\varphi (\overline{x}),\KK)$ and $\Phi (R(\HH,\overline{y})) = R(\KK,\psi (\overline{y}))$, respectively. Then, both $\varphi$ and $\psi$ are bijection between the collections of one-dimensional subspaces of $\HH$ and $\KK$. Moreover, for each $x,y \in \HH \setminus \{0\}$, we have
\begin{align*}
\Phi (\{ \ker \omega_{x,y} \cap \A\}) 
&= \Phi (R(\overline{x},\HH) \cap R(\HH,\overline{y})) \\
&= \Phi (R(\overline{x},\HH)) \cap \Phi (R(\HH,\overline{y})) \\
&= R(\varphi (\overline{x}),\KK) \cap R(\KK,\psi (\overline{y})) .
\end{align*}
Hence, $\Phi (\ker \omega_{x,y} \cap \A) = \ker \omega_{z,w} \cap \B$ for any choice of $z \in \varphi (\overline{x}) \setminus \{0\}$ and $w \in \psi (\overline{y}) \setminus \{0\}$. For each $x,y \in \HH \setminus \{0\}$, choose $z(x),w(y) \in \KK \setminus \{0\}$ such that $z(x) \in \varphi (\overline{x})$ and $w(y) \in \psi (\overline{y})$. It follows that
\[
\Phi (\ker \omega_{x,y} \cap \A) = \ker \omega_{z(x),w(y)} \cap \B
\]
for each $x,y \in \HH \setminus \{0\}$.

The next step was inspired by \cite[Lemma 1]{FL84}. Let $\mathcal{M}(\HH)$ and $\mathcal{M}(\KK)$ be the collections of all subspaces of $\HH$ and $\KK$. For each $M \in \mathcal{M}(\HH)$, let
\[
\varphi_0 (M) = \langle \{ z(x) : x \in M \setminus \{0\} \} \rangle .
\]
Then, we note that $\varphi_0 (\{0\} ) = \langle \emptyset \rangle = \{0\}$ and $\varphi_0 (\overline{x}) = \varphi (\overline{x})$ for each $x \in \HH \setminus \{0\}$. We claim that
\[
M \setminus \{0\} = \{ x \in \HH \setminus \{0\} : z(x) \in \varphi_0 (M)\}
\]
for each $M \in \mathcal{M}(\HH)$. Indeed, if $x \in M \setminus \{0\}$, then $z(x) \in \varphi_0(M)$ by the definition of $\varphi_0(M)$. Conversely, suppose that $x \in \HH \setminus \{0\}$ and $z(x) \in \varphi_0 (M)$. Then, there exist $x_1,\ldots ,x_n \in M \setminus \{0\}$ and $c_1,\ldots ,c_n \in \mathbb{C}$ such that $z(x) = \sum_{j=1}^n c_jz(x_j)$. Fix an arbitrary $y \in \HH \setminus \{0\}$. Since $\omega_{z(x),w(y)} = \sum_{j=1}^n c_j\omega_{z(x_j),w(y)}$, it follows that $\bigcap_{j=1}^n \ker \omega_{z(x_j),w(y)} \cap \B \subset \ker \omega_{z(x),w(y)} \cap \B$, that is, $\Phi (\ker \omega_{x,y} \cap \A) \in \Phi (\{ \ker \omega_{x_j,y} : j \in \{1,\ldots ,n\} \})^=$. Combining this with
\begin{align*}
\lefteqn{\Phi (\{ \ker \omega_{x_j,y} \cap \A : j \in \{1,\ldots ,n\} \})^= \cap \mathfrak{S}_n(\B)} \\
&= \Phi (\{ \ker \omega_{x_j,y} \cap \A : j \in \{1,\ldots ,n\} \}^= \cap \mathfrak{S}_n(\A)) ,
\end{align*}
we obtain $\ker \omega_{x,y} \cap \A \in \{ \ker \omega_{x_j,y} \cap \A : j \in \{1,\ldots ,n\} \}^=$. Let $E$ be the projection from $\HH$ onto $\langle \{x_1,\ldots ,x_n\} \rangle$. If $x \not \in \langle \{x_1,\ldots ,x_n\} \rangle$, then $(I-E)x \neq 0$. Since $\A$ acts irreducibly on $\HH$, there exists an element $A \in \A$ such that $Ay = (I-E)x$. It follows that $\omega_{x_j,y}(A^*) = 0$ for each $j$ and $\omega_{x,y}(A^*) = \|(I-E)x\|^2 >0$, which contradicts $\bigcap_{j=1}^n \ker \omega_{x_j,y} \cap \A \subset \ker \omega_{x,y} \cap \A$. Therefore, $x \in \langle \{x_1,\ldots ,x_n\} \rangle \subset M$. This proves the claim, which shows that $\varphi_0$ is injective.

Next, let $M_N = \{ x \in \HH \setminus \{0\} : z(x) \in N \} \cup \{0\}$ for each $N \in \mathcal{M}(\KK)$. We claim that $M_N \in \mathcal{M}(\HH)$ and $N = \varphi_0 (M_n)$. To show this, take an arbitrary $x \in \langle M_N \rangle \setminus \{0\}$ and an arbitrary $y \in \HH \setminus \{0\}$. Then, $x \in \langle x_1,\ldots ,x_n \rangle$ for some $x_1,\ldots ,x_n \in M_N \setminus \{0\}$, which implies that
\[
\ker \omega_{x,y} \cap \A \in \{ \ker \omega_{x_j,y} \cap \A : j \in \{1,\ldots ,n\}\}^=  \cap \mathfrak{S}_n(\A).
\]
It follows from
\begin{align*}
\lefteqn{\Phi (\{ \ker \omega_{x_j,y} \cap \A : j \in \{1,\ldots ,n\}\}^= \cap \mathfrak{S}_n(\A))}\\
&= \Phi (\{ \ker \omega_{x_j,y} \cap \A : j \in \{1,\ldots ,n\}\})^= \cap \mathfrak{S}_n(\B)
\end{align*}
that $\ker \omega_{z(x),w(y)} \cap \B \in \{ \ker \omega_{z(x_j),w(y)} \cap \B : j\in \{1,\ldots ,n\} \}^=$. Let $F$ be the projection from $\KK$ onto $\langle z(x_1),\ldots ,z(x_n) \rangle$. If $z(x) \not \in \langle z(x_1),\ldots ,z(x_n) \rangle$, then $(I-F)z(x) \neq 0$. Since $\B$ acts irreducibly on $\KK$, there exists an element $B \in \B$ such that $Bw(y) = (I-F)z(x)$. However, this means that
\[
B^* \in \bigcap_{j=1}^n (\ker \omega_{z(x_j),w(y)} \cap \B) \setminus (\ker \omega_{z(x),w(y)} \cap \B) ,
\]
which contradicts $\ker \omega_{z(x),w(y)} \cap \B \in \{ \ker \omega_{z(x_j),w(y)} \cap \B : j\in \{1,\ldots ,n\} \}^=$. Hence, $z(x) \in \langle z(x_1),\ldots ,z(x_n) \rangle \subset N$, that is, $x \in M_N$. This shows that $\langle M_N \rangle \setminus \{0\} = M_N \setminus \{0\}$, which implies that $M_N = \langle M_N \rangle \in \mathcal{M}(\HH)$. Since $\{ z(x) : x \in M_N \setminus \{0\}\} \subset N$, it follows that $\varphi_0 (M_N) \subset N$. Conversely, if $z \in N \setminus \{0\}$, then $\varphi (\overline{x}) = \overline{z}$ for some $x \in \HH \setminus \{0\}$, in which case, $z(x) = cz$ for some $c \in \mathbb{C} \setminus \{0\}$. In particular, we have $x \in M_N \setminus \{0\}$ and $z = c^{-1}z(x) \in \varphi_0 (M_N)$. Therefore, $N = \varphi_0 (M_N)$. We note that $M_{\KK} = \HH$ and $\varphi_0 (\HH) = \KK$.

By the preceding two paragraphs, we know that $\varphi_0 :\mathcal{M}(\HH) \to \mathcal{M}(\KK)$ is bijective. Moreover, it is in fact an order isomorphism. From this, $\dim \HH =2$ if and only if $\dim \KK=2$. Indeed, if $\dim \HH = 2$, then $\{0\} \subsetneq \overline{x} \subsetneq \HH$ is a maximal chain in $\mathcal{M}(\HH)$ for each $x \in \HH \setminus \{0\}$. Since $\varphi_0$ is an order isomorphism, $\{0\} = \varphi_0 (\{0\}) \subsetneq \varphi_0 (\overline{x}) \subsetneq \varphi_0 (\HH) =\KK$ is also a maximal chain in $\mathcal{M}(\KK)$. This shows that $\dim \KK = 2$. The converse also follows from the same argument as above. Now, suppose that $\dim \HH \geq 3$ (and $\dim \KK \geq 3$). Then, by the fundamental theorem of projective geometry, there exists a semilinear bijection $U:\HH \to \KK$ such that $\varphi_0 (M) = U(M)$ for each $M \in \mathcal{M}(\HH)$. In particular, $\varphi (\overline{x})=\varphi_0(\overline{x}) = U(\overline{x})$ for each $x \in \HH \setminus \{0\}$. Since $U$ preserves linear independence of finite subsets of $\HH$ and $\KK$, it follows that $\HH$ is finite-dimensional if and only if $\KK$ is, and $\dim \HH = \dim \KK$ provided that either $\HH$ or $\KK$ is finite-dimensional.

Suppose that $\HH$ and $\KK$ are infinite-dimensional and $\Phi$ is continuous. Let $N$ be a proper closed subspace of $\KK$, let $y \in \HH \setminus \{0\}$, and let $(x_n)_n$ be a sequence in $\varphi_0^{-1}(N) = M_N$ that converges to $x \in \HH \setminus \{0\}$. It may be assumed that $x_n \neq 0$ for each $n$. Then, we obtain $\ker \omega_{x,y} \cap \A \in \{ \ker \omega_{x_n,y} \cap \A : n \in \mathbb{N}\}^= \cap \mathfrak{S}_n(\A)$, which together with
\[
\Phi (\{ \ker \omega_{x_n,y} \cap \A : n \in \mathbb{N}\}^= \cap \mathfrak{S}_n(\A)) \subset \Phi (\{ \ker \omega_{x_n,y} \cap \A : n \in \mathbb{N}\})^= \cap \mathfrak{S}_n(\B)
\]
implies that $\ker \omega_{z(x),w(y)} \cap \B \in \{ \ker \omega_{z(x_n),w(y)} \cap \B : n \in \mathbb{N} \}^=$. Let $F_N$ be the projection from $\KK$ onto $N$. If $z(x) \not \in N$, then an element $B_N \in \B$ with $B_N w(y) = (I-F_N)z(x)$ satisfies
\[
B_N^* \in \bigcap_n (\ker \omega_{z(x_n),w(y)} \cap \B) \setminus (\ker \omega_{z(x),w(y)} \cap \B)
\]
by $(z(x_n))_n \subset N$, a contradiction. Therefore, $z(x) \in N$, that is, $x \in M_N = \varphi_0^{-1}(N)$. This shows that $U^{-1}(N) = \varphi_0^{-1}(N)$ is closed whenever $N$ is a closed subspace of $\KK$, which together with Lemma~\ref{auto-conti} implies that $U^{-1}$ (and $U$) is linear or conjugate-linear. Considering $U$ as a mapping from $\overline{\HH}$ onto $\KK$ if necessary, we may assume that $U$ is linear. If $g \in \KK^*$, then $\ker g$ is a closed subspace of $\KK$, which implies that $\ker (g \circ U) = U^{-1}(\ker g)$ is a closed subspace of $\HH$ (or $\overline{\HH}$). This means that $g \circ U \in \HH^*$ (or $g \circ U \in \overline{\HH}^*$) whenever $g \in \KK^*$, that is, $U$ is weak-to-weak continuous. Therefore, $U$ is bounded, and $U^{-1}$ is also bounded by the open mapping theorem. Similarly (or applying the same argument as above to the mapping $\Psi = \Phi^{\B}_a \circ \Phi \circ \Phi^{\A}_a$), provided that $\dim \HH \geq 3$ (and $\dim \KK \geq 3$), we can show that there exists a semilinear bijection $V:\HH \to \KK$ (which is bounded linear or conjugate-linear isomorphism if $\HH$ and $\KK$ are infinite-dimensional and $\Phi$ is continuous) such that $\psi (\overline{y}) = V(\overline{y})$ for each $y \in \HH \setminus \{0\}$. Consequently, we have
\[
\Phi (\ker \omega_{x,y} \cap \A ) = \ker \omega_{Ux,Vy} \cap \B
\]
and
\[
(\Phi^{\B}_a \circ \Phi)(\ker \omega_{x,y} \cap \A) = \ker \omega_{Vy,Ux} \cap \B
\]
for each $x,y \in \HH \setminus \{0\}$. This completes the proof of the existence of $U$ and $V$ that satisfy (i) or (ii).

Finally, suppose that both $\HH$ and $\KK$ are infinite-dimensional and $\Phi$ is continuous. In this case, we can choose bounded linear or conjugate-linear isomorphisms $U$ and $V$ that satisfy (i) or (ii). Since $\B$ acts irreducibly on $\KK$, either $K(\KK) \subset \B$ or $K(\KK ) \cap \B = \{0\}$ by \cite[Corollary IV.1.2.5]{Bla06}. If $K(\KK ) \subset \B$, let $B_0$ be any finite-rank projection on $\KK$ with $\dim B_0 (\KK) \geq 2$. In the case that $K(\KK) \cap \B = \{0\}$, take an arbitrary nonzero positive element $B_0 \in \B$. Then, $B_0(\HH)$ is infinite-dimensional since $B_0$ cannot be finite-rank. Let $z_0 \in \KK \setminus \{0\}$ be such that $\langle B_0z_0,z_0 \rangle >0$, and let $S(B_0) = \{ \ker \omega_{z,w} \cap \B \in \mathfrak{S}(\B) : B_0 \in \ker \omega_{z,w} \cap \B \}$. Since $S(B_0)^= \cap \s_n (\B) = S(B_0)$ and $\Phi$ is continuous, we obtain
\[
\Phi^{-1}(S(B_0))^= \cap \mathfrak{S}_n(\A) = \Phi^{-1}(S(B_0))
\]
and in particular
\[
\Phi^{-1}(\ker \omega_{z_0,z_0} \cap \B) \not \in \Phi^{-1}(S(B_0))^= \cap \s_n (\A)
\]
by $\ker \omega_{z_0,z_0} \cap \B \not \in S(B_0)$. From this and $\Phi^{-1}(\ker \omega_{z_0,z_0} \cap \B)\in \s_n (\A)$, there exists an element $A_0 \in \A$ such that
\[
A_0 \in \bigcap \{ \Phi^{-1} (\ker \omega_{z,w} \cap \A) : z,w \in \KK \setminus \{0\},~\langle B_0z,w \rangle = 0\} \setminus \Phi^{-1}(\ker \omega_{z_0,z_0} \cap \B) .
\]
If $\Phi (\omega_{x,y} \cap \A) = \omega_{Ux,Vy} \cap \B$ for each $x,y \in \HH \setminus \{0\}$, then
\[
\Phi^{-1}(\ker \omega_{z,w} \cap B) = \ker \omega_{U^{-1}z,V^{-1}w} \cap \A
\]
for each $z,w \in \KK \setminus \{0\}$, which implies that $\langle A_0U^{-1}z,V^{-1}w \rangle = 0$ whenever $\langle B_0z,w \rangle = 0$ and $\langle A_0U^{-1}z_0,V^{-1}z_0\rangle \neq 0$. Since $(V^{-1})^*A_0U^{-1}$ is a bounded linear or conjugate-linear operator on $\KK$, Lemma~\ref{sublemma-semi} ensures that $(V^{-1})^*A_0U^{-1} = cB_0$ for some $c \in \mathbb{C} \setminus \{0\}$. In particular, $V^*B_0U$ is nonzero and linear. From this, $U$ is linear if and only if $V$ is. Meanwhile, if $\Phi (\ker \omega_{x,y} \cap \A) = \ker \omega_{Vy,Ux} \cap \B$ for each $x,y \in \HH \setminus \{0\}$, then
\[
\Phi^{-1}(\ker \omega_{z,w} \cap \B) = \ker \omega_{U^{-1}w,V^{-1}z} \cap \A
\]
for each $z,w \in \KK \setminus \{0\}$. In this case, we have $\langle (U^{-1})^*A_0^*V^{-1}z,w \rangle = 0$ whenever $\langle B_0z,w \rangle = 0$ and $\langle (U^{-1})^*A_0^*V^{-1}z_0,z_0 \rangle \neq 0$, which together with Lemma~\ref{sublemma-semi} implies that $(U^{-1})^*A_0^*V^{-1} = cB_0$ for some $c \in \mathbb{C}\setminus \{0\}$. In particular, $U^*B_0V$ is nonzero and linear. Therefore, also in this case, $U$ is linear if and only if $V$ is. We note that if $\Phi (\ker \omega_{x,y} \cap \A) = \ker \omega_{Vy,Ux} \cap \B$ for each $x,y \in \HH \setminus \{0\}$, then $(\Phi^{\mathfrak{B}} \circ \Phi)(\ker \omega_{x,y} \cap \A) = \ker \omega_{Ux,Vy}$ for each $x,y \in \HH \setminus \{0\}$. Hence, in either case, the mapping $\ker \omega_{x,y} \cap \A \mapsto \ker \omega_{Ux,Vy} \cap \B$ is continuous if $\Phi$ is. Combining this with Theorem~\ref{continuity}, we have $V^*\B U \subset \mathfrak{A}$ (and $V^*\B U=\A$ if $\Phi$ is a homeomorphism), as desired.
\end{proof}
\end{theorem}


\section{The Naimark problem revisited}\label{Sect:NP}

The main objective of this section is to present a characterization elementary $C^*$-algebras in terms of their geometric structure spaces. In this direction, the separable case has been already completed by Rosenberg~\cite{Ros53} as was mentioned in the introduction, since the spectrum of a $C^*$-algebra can be recognized as a substructure of its geometric structure space. To be precise, a separable $C^*$-algebra $\A$ is elementary if and only if $\hat{\A}$ is a singleton. Meanwhile, by the result of Akemann and Weaver~\cite{AW04}, we cannot show in ZFC that non-separable $C^*$-algebras with singleton spectra are always elementary. Under these circumstances, we try to find a property $P$ for geometric structure spaces of (non-separable) $C^*$-algebras $\A$ such that $\A$ is elementary if and only if $\hat{\A}$ is a singleton and $\A$ has the property $P$ as well as the property $P$ does not imply that $\hat{\A}$ is a singleton. To this end, we first characterize $C^*$-algebras with singleton spectra in terms of their geometric structure spaces. The following theorem gives some conditions equivalent to the reverse inclusion of Corollary~\ref{vector-funct}.
\begin{theorem}\label{class-closed}

Let $\A$ be a $C^*$-algebra, and let $\rho$ be a pure state of $\A$. Then, the following are equivalent:
\begin{itemize}
\item[{\rm (i)}] $C(\ker \rho )^= = C(\ker \rho)$.
\item[{\rm (ii)}] $\mathfrak{S}(\pi_\rho (\A)) = \{ \ker \omega_{x,y} \cap \pi_\rho (\A) : x,y \in S_{\HH_\rho}\}$.
\item[{\rm (iii)}] $\mathfrak{S}(\pi (\A)) = \{ \ker \omega_{x,y} \cap \pi (\A) : x,y \in S_{\HH}\}$ for each $\pi\in \overline{\pi_\rho}$.
\item[{\rm (iv)}] $\mathfrak{S}(\pi (\A)) = \{ \ker \omega_{x,y} \cap \pi (\A) : x,y \in S_{\HH}\}$ for some $\pi \in \overline{\pi_\rho}$.
\end{itemize}
\begin{proof}
Let $\Phi (\ker \tau ) = \ker \pi_\rho^* (\tau )$ for each $\ker \tau \in \mathfrak{S}(\pi_\rho (\A))$. Since $\ker \pi_\rho = J_{C(\ker \rho)}$ by Lemma~\ref{primitive}, it follows from Theorem~\ref{irr-rep} that $\Phi$ is a homeomorphism from $\mathfrak{S}(\pi_\rho (\A))$ onto $C(\ker \rho)^=$. In particular, $\Phi (\mathfrak{S}(\pi_\rho (\A))) = C(\ker \rho)^=$. Meanwhile, we have
\begin{align*}
C(\ker \rho ) 
&= \{ \ker \pi_\rho^* (\omega_{x,y}|\pi_\rho (\A)) : x,y \in S_{\HH_\rho}\} \\
&= \Phi (\{ \ker \omega_{x,y} \cap \pi_\rho (\A) : x,y \in S_{\HH_\rho}\})
\end{align*}
by Lemma~\ref{pure-adjoint}. Hence, we have $C(\ker \rho)^= = C(\ker \rho)$ if and only if $\mathfrak{S}(\pi_\rho (\A)) = \{ \ker \omega_{x,y} \cap \pi_\rho (\A) : x,y \in S_{\HH_\rho}\}$. This proves (i) $\Leftrightarrow$ (ii).

Suppose that (i) holds. Let $\pi:\A \to \HH$ be an irreducible representation of $\A$ that is unitarily equivalent to $\pi_\rho$, and let $U:\HH_\rho \to \HH$ be an isometric isomorphism such that $\pi (A) = U\pi_\rho (A)U^*$ for each $A \in \A$. If $\ker \tau \in \mathfrak{S}(\pi (\A))$, then we obtain $\pi^* (\tau ) \in \ext (B_{\A^*})$ and
\[
J_{C(\ker \rho)} = \ker \pi_\rho = \ker \pi \subset \ker \pi^* (\tau ) ,
\]
which together with Lemmas~\ref{closed} and \ref{pure-adjoint} implies that
\[
\pi^* (\tau ) \in C(\ker \rho )^= = C(\ker \rho ) = \{ \ker \pi_\rho^* (\omega_{x,y}|\pi_\rho (\A)) : x,y \in S_{\HH_\tau}\} .
\]
Hence, $\pi^* (\tau ) = \pi_\rho^* (\omega_{x,y}|\pi_\rho (\A))$ for some $x,y \in S_{\HH_\rho}$, in which case, we have
\[
\tau (\pi (A)) = \omega_{x,y}(\pi_\rho (A)) = \langle \pi_\rho (A)x,y \rangle = \langle U^*\pi (A)Ux,y \rangle
\]
for each $A \in \A$. This shows that $\tau = \omega_{Ux,Uy}|\pi (\A)$. Combining this with Corollary~\ref{vector-funct}, we derive that $\mathfrak{S}(\pi (\A)) = \{ \ker \omega_{x,y} \cap \pi (\A) : x,y \in S_{\HH}\}$. Therefore, (i) $\Rightarrow$ (iii). The implication (iii) $\Rightarrow$ (iv) is obvious.

Finally, suppose that $\mathfrak{S}(\pi (\A)) = \{ \ker \omega_{x,y} \cap \pi (\A) : x,y \in S_{\HH}\}$ for some $\pi \in \overline{\pi_\rho}$. Let $V:\HH \to \HH_\rho$ be an isometric isomorphism such that $\pi (A) = V^*\pi_\rho (A)V$ for each $A \in \A$, and let $\psi (\pi (A)) = \pi_\rho (A)$ for each $A \in \A$. Then, $\psi :\pi (\A) \to \pi_\rho (\A)$ is a $*$-isomorphism. By Theorem~\ref{Tan-310}, setting $\Psi (I) = \psi (I)$ for each $I \in \s (\pi (\A))$, we have a homeomorphism $\Psi : \s (\pi (\A)) \to \s (\pi_\rho (\A))$. Since
\begin{align*}
\Psi (\ker \omega_{x,y} \cap \pi (\A)) = \psi (\ker \omega_{x,y} \cap \pi (\A)) 
&= \ker (\omega_{x,y} \circ \psi^{-1}) \cap \pi_\rho (\A) \\
&= \ker \omega_{Vx,Vy} \cap \pi_\rho (\A) ,
\end{align*}
it follows that
\begin{align*}
\s (\pi_\rho (\A)) = \Psi (\s (\pi (\A))) 
&= \{ \Psi (\ker \omega_{x,y} \cap \pi (\A)) : x,y \in S_{\HH} \} \\
&= \{ \ker \omega_{Vx,Vy} \cap \pi_\rho (\A) : x,y \in S_\HH \} \\
&= \{ \ker \omega_{z,w} \cap \pi_\rho (\A) : z,w \in S_{\HH_\rho} \} .
\end{align*}
Thus, (iv) $\Rightarrow$ (ii).
\end{proof}
\end{theorem}
\begin{corollary}\label{all-closed}

Let $\A$ be a $C^*$-algebra. Then, $C(I)^= = C(I)$ for each $I \in \mathfrak{S}(\A)$ if and only if $\mathfrak{S}(\pi (\A)) = \{ \ker \omega_{x,y} \cap \pi (\A) : x,y \in S_{\HH}\}$ for each irreducible representation $\pi :\A \to B(\HH)$.
\begin{proof}
Suppose that $C(I)^= = C(I)$ for each $I \in \mathfrak{S}(\A)$, and that $\pi :\A \to B(\HH)$ is an irreducible representation of $\A$. Then, there exists a pure state $\rho$ of $\A$ such that $\pi$ and $\pi_\rho$ are unitarily equivalent. Since $C(\ker \rho )^= = C(\ker \rho )$, it follows from Theorem~\ref{class-closed} that $\mathfrak{S}(\pi (\A)) = \{ \ker \omega_{x,y} \cap \pi (\A) : x,y \in S_{\HH}\}$.

Conversely, suppose that $\mathfrak{S}(\pi (\A)) = \{ \ker \omega_{x,y} \cap \pi (\A) : x,y \in S_{\HH}\}$ for each irreducible representation $\pi :\A \to \HH$. Let $I \in \mathfrak{S}(\A)$. Then, by Lemma~\ref{polar-equi}, there exists a pure state $\rho$ such that $C(I) = C(\ker \rho)$. Since $\mathfrak{S}(\pi_\rho (\A)) = \{ \ker \omega_{x,y} \cap \pi_\rho (\A) : x,y \in S_{\HH_\rho}\}$, it follows from Theorem~\ref{class-closed} that $C(I)^= = C(\ker \rho )^= = C(\ker \rho) =C(I)$.
\end{proof}
\end{corollary}
\begin{corollary}\label{simple-Glimm}

Let $\A$ be a $C^*$-algebra. Then, $\hat{\A}$ (or equivalently, $\mathfrak{GS}(\A)$) is a singleton if and only if there exists a faithful irreducible representation $\pi :\A \to B(\HH)$ of $\A$ such that $\mathfrak{S}(\pi (\A)) = \{ \ker \omega_{x,y} \cap \pi (\A ) : x,y \in S_{\HH} \}$.
\begin{proof}
Suppose that $\mathfrak{GS}(\A)$ is a singleton. Then, $C(I) = \mathfrak{S}(\A) \in \mathfrak{C}(\A)$ for any choice of $I \in \mathfrak{S}(\A)$. Hence, by Corollary~\ref{all-closed}, $\mathfrak{S}(\pi (\A)) = \{ \ker \omega_{x,y} \cap \pi (\A ) : x,y \in S_{\HH} \}$ for each irreducible representation $\pi :\A \to B(\HH)$. We note that since $\hat{\A}$ is also a singleton, each irreducible representation $\pi:\A \to B(\HH)$ is faithful.

Conversely, suppose that there exists a faithful irreducible representation $\pi :\A \to B(\HH)$ of $\A$ such that $\mathfrak{S}(\pi (\A)) = \{ \ker \omega_{x,y} \cap \pi (\A ) : x,y \in S_{\HH} \}$. Let $\rho$ be a pure state of $\A$ such that $\pi$ and $\pi_\rho$ are unitarily equivalent. Then, we obtain $C(\ker \rho)^= = C(\ker \rho)$ by Theorem~\ref{class-closed}. Meanwhile, since $\pi$ is faithful, it follows from Lemma~\ref{primitive} that $J_{C(\ker \rho)} = \ker \pi_\rho = \ker \pi = \{0\}$. Combining this with Lemma~\ref{closed}, we infer that $\s (\A) = C(\ker \rho)^= = C(\ker \rho)$. Thus, $\gs (\A)$ is a singleton.
\end{proof}
\end{corollary}
This completes our first task. Next, we introduce an important subcategory of $\mathbb{C}$-$\Gss$.
\begin{definition}

The full subcategory of $\mathbb{C}$-$\Gss$ whose objects are homeomorphs of geometric structure spaces of nontrivial $C^*$-algebras is denoted by $\CGss$.
\end{definition}
\begin{remark}

From the definition, $\CGss$ is a strictly full subcategory of $\mathbb{C}$-$\Gss$. We note that if $(K,c)$ is an object of $\CGss$, then $(K^\g,c^\g)$ and $(K^\p,c^\p)$ are both topologizable by Proposition~\ref{strictly-full}, Theorems~\ref{two-homeo} and \ref{spec-trans}, and Corollary~\ref{top-able}.
\end{remark}
We make use of a characterization of norm-closed ideals of $B(\HH)$ by Luft~\cite[Corollary 6.1]{Luf68}; see also \cite{Gra67}. Recall that, for a cardinal number $\alpha \geq \aleph_0$, an element $A \in B(\HH)$ is said to be \emph{$\alpha$-compact} if $\dim M <\alpha$ whenever $M$ is a closed subspace of $A(\HH)$. Let $K_\alpha (\HH)$ be the set of all $\alpha$-compact operators on $B(\HH)$. Then, $K_{\aleph_0}(\HH)=K(\HH)$ and $K_\alpha (\HH) \subset K_\beta (\HH)$ whenever $\alpha 
 \leq \beta$.
\begin{theorem}[Luft, 1968]\label{Luft}

Let $\HH$ be a Hilbert space. Then, $\mathcal{I}$ is a norm-closed ideal of $B(\HH)$ if and only if $\mathcal{I}=K_\alpha (\HH)$ for some cardinal number $\alpha \geq \aleph_0$.
\end{theorem}
If $U$ and $V$ are bounded conjugate-linear isomorphism from $\HH$ onto $\KK$, and if $A \in K_\alpha (\KK)$, then $V^*AU \in K_\alpha (\HH)$. Indeed, if $M$ is a closed subspace of $V^*AU(\HH)$, then $(V^{-1})^*(M)$ is a closed subspace of $AU(\HH) = A(\KK)$, which implies that $\dim M = \dim (V^{-1})^*(M) <\alpha$. This shows that $V^*AU \in K_\alpha (\HH)$. In summary, we have the identity $V^*K_\alpha (\KK)U = K_\alpha (\HH)$.

The following is the main result in this paper.
\begin{theorem}\label{elementary}

Let $\A$ be a non-separable $C^*$-algebra. Then, $\A$ is elementary if and only if it satisfies the following conditions:
\begin{itemize}
\item[{\rm (i)}] $\hat{\A}$ is a singleton.
\item[{\rm (ii)}] There exists an element $I \in \s (\A)$ with the following property: If $(K,c)$ is an object of $\CGss$, and if there exists a continuous finite-homeomorphism $f:C(t) \to C(I)$ for an element $t \in K$, then $f \circ g \circ f^{-1} $ is a homeomorphism on $C(I)$ whenever $g$ is a continuous finite-homeomorphism on $C(t)$.
\end{itemize}
\begin{proof}
Each elementary $C^*$-algebra has the singleton spectrum by Rosenberg's theorem. We assume that $\A$ is $*$-isomorphic to a norm-closed ideal $\mathcal{I}$ of $B(\HH)$ for some infinite-dimensional Hilbert space $\HH$. Then, there exists a homeomorphism $\Phi :\s (\A) \to \s (\mathcal{I})$ by Theorem~\ref{Tan-310}, and $\mathcal{I}$ acts irreducibly on $\HH$ since it contains $K(\HH)$. Let $I = \Phi^{-1}(\ker \omega_{x,x} \cap \mathcal{I})$ for some $x \in S_\HH$. We note that $\Phi (C(I)) = C(\Phi (I)) = \s_n (\mathcal{I})$ by Lemmas~\ref{finite-conti} and \ref{normal-class}.

Let $(K,c)$ be an object of $\CGss$, and let $f:C(t) \to C(I)$ be a continuous finite-homeomorphism. Then, there exist a $C^*$-algebra $\B$ and a homeomorphism $f_0 :\s (\B) \to K$. Define another mapping $\Phi_0 = \Phi \circ f \circ f_0$ from $C(f_0^{-1}(t))$ onto $\s_n (\mathcal{I})$. Proposition~\ref{subspace-homeo} ensures that $\Phi_0$ is a continuous finite-homeomorphism. If $g$ is a continuous finite-homeomorphism on $C(t)$, then $\Psi = f_0^{-1} \circ g \circ f_0$ is a continuous finite-homeomorphism on $C(f_0^{-1}(t))$. We will know that $f \circ g \circ f^{-1}$ is a homeomorphism on $C(I)$ once it is proved that $\Phi_0 \circ \Psi \circ \Phi_0^{-1} = \Phi \circ (f \circ g \circ f^{-1}) \circ \Phi^{-1}$ is a homeomorphism on $C(\Phi (I))$.

By Lemma~\ref{polar-equi}, there exists a pure state $\rho$ of $\B$ such that $C(f_0^{-1}(t)) = C(\ker \rho)$, in which case, by Theorem~\ref{irr-rep}, the mapping $\Psi_0 : \mathfrak{S}(\pi_\rho (\B)) \to C(\ker \rho)^=$ given by $\Psi (\ker \tau ) = \ker \pi_\rho^* (\tau )$ for each $\ker \tau \in \mathfrak{S}(\pi_\rho (\B))$ is a homeomorphism. Moreover, $\Psi_0 (\mathfrak{S}_n(\pi_\rho (\B))) = C(\ker \rho)$ by Lemma~\ref{pure-adjoint}. Hence, the mapping $\Phi_0 \circ \Psi_0$ gives a continuous finite-homeomorphism from $\mathfrak{S}_n (\pi_\rho (\B))$ onto $\mathfrak{S}_n(\mathcal{I})$. Applying Theorem~\ref{weak-preserver}, we obtain $\dim \HH = \dim \HH_\rho$ and bounded linear or conjugate-linear isomorphisms $U$ and $V$ from $\HH_\rho$ onto $\HH$ such that $U$ is linear if and only if $V$ is, $V^*\mathcal{I}U \subset \pi_\rho (\B)$, and either
\[
(\Phi_0 \circ \Psi_0)(\ker \omega_{z,w} \cap \pi_\rho (\B)) = \ker \omega_{Uz,Vw} \cap \mathcal{I}
\]
for each $z,w \in \HH_\rho \setminus \{0\}$ or
\[
(\Phi_0 \circ \Psi_0)(\ker \omega_{z,w} \cap \pi_\rho (\B)) = \ker \omega_{Vw,Uz} \cap \mathcal{I}
\]
for each $z,w \in \HH_\rho \setminus \{0\}$. Replacing $\Phi_0 \circ \Psi_0$ with $\Phi^{\mathcal{I}}_a \circ \Phi_0 \circ \Psi_0$ if necessary, we may assume that the former case occurs.

Let $\Psi' = \Psi_0^{-1} \circ \Psi \circ \Psi_0$. Then, $\Psi'$ is a continuous finite-homeomorphism on $\pi_\rho (\B)$ and
\[
\Psi' (\mathfrak{S}_n(\pi (\B))) = (\Psi_0^{-1} \circ \Psi)(C(\ker \rho )) = \Psi_0^{-1}(C(\ker \rho)) = \mathfrak{S}_n(\pi (\B)) .
\]
Hence, by Proposition~\ref{subspace-homeo} and Theorem~\ref{weak-preserver}, there exist bounded linear or conjugate-linear isomorphisms $U',V'$ on $\HH_\rho$ such that $U'$ is linear if and only if $V'$ is, $V'^*\pi (\B)U' = \pi (\B)$, and either
\[
\Phi' (\ker \omega_{x,y} \cap \pi_\rho (\B)) = \ker \omega_{U'x,V'y} \cap \pi_\rho (\B)
\]
for each $x,y \in \HH_\rho \setminus \{0\}$ or
\[
\Phi' (\ker \omega_{x,y} \cap \pi_\rho (\B)) = \ker \omega_{V'y,U'x} \cap \pi_\rho (\B)
\]
for each $x,y \in \HH_\rho \setminus \{0\}$. It follows from
\[
\Phi_0 \circ \Psi \circ \Phi_0^{-1} = (\Phi_0 \circ \Psi_0) \circ \Psi' \circ (\Phi_0 \circ \Psi_0)^{-1}
\]
that
\[
(\Phi_0 \circ \Psi \circ \Phi_0^{-1})(\ker \omega_{x,y} \cap \mathcal{I}) = \ker \omega_{UU'U^{-1}x,VV'V^{-1}y} \cap \mathcal{I}
\]
for each $x,y \in \HH \setminus \{0\}$ in the former case, and
\[
(\Phi_0 \circ \Phi \circ \Phi_0^{-1})(\ker \omega_{x,y} \cap \mathcal{I}) = \ker \omega_{UV'V^{-1}y,VU'U^{-1}x} \cap \mathcal{I}
\]
for each $x,y \in \HH \setminus \{0\}$ in the latter case. We note that 
\[
UU'U^{-1},VV'V^{-1},UV'V^{-1},VU'U^{-1}
\]
are all bounded linear or conjugate-linear isomorphisms on $\HH$. Moreover, $UU'U^{-1}$ is linear if and only if $U'$ is linear, while $VV'V^{-1}$ is linear if and only if $V'$ is linear. This shows that $UU'U^{-1}$ is linear if and only if $VV'V^{-1}$ is. Similarly since $U$ is linear if and only $V$ is, it follows that $UV'V^{-1}$ is linear if and only if $VU'U^{-1}$ is. Combining these with the fact that $\mathcal{I}$ is an ideal of $B(\HH)$, in the light of Theorems~\ref{continuity} and \ref{Luft}, we infer that $\Phi_0 \circ \Psi \circ \Phi_0^{-1}$ is a homeomorphism on $C(\Phi (I))$. Thus, $\A$ satisfies (ii) if it is $*$-isomorphic to a norm-closed ideal of some infinite-dimensional $B(\HH)$.

For the converse, suppose that (i) and (ii) hold. By Corollary~\ref{simple-Glimm}, it may be assumed that there $\A$ is an irreducible $C^*$-algebra acting on a Hilbert space $\HH$ and $\mathfrak{S}(\A) = \{ \ker \omega_{x,y} \cap \A : x,y \in \HH \setminus \{0\}\}$. In particular, we have $C(I) = \s (\A) = \s_n (\A)$ for any choice of $I \in \s (\A)$. Let $U$ be a unitary operator on $\HH$. Define mappings $\Phi_0 : \mathfrak{S}_n(B(\HH)) \to \mathfrak{S}(\A)$ and $\Phi_U : \mathfrak{S}(B(\HH)) \to \mathfrak{S}(B(\HH))$ by $\Phi_0 (\ker \omega_{x,y}) = \ker \omega_{x,y} \cap \A$ for each $x,y \in \HH \setminus \{0\}$ and $\Phi_U (\ker \rho ) = \ker \rho U$ for each $\ker \rho \in \mathfrak{S}(B(\HH))$. Then, $\Phi_0$ is a continuous finite-homeomorphism by Theorem~\ref{continuity}. Moreover, $\Phi_U$ is a homeomorphism on $\mathfrak{S}(B(\HH))$ with $\Phi_U (\s_n (B(\HH))) = \s_n (B(\HH))$. Indeed, let $\varphi_{U^*}(A) = U^*A$ for each $A \in B(\HH)$. Then, $\varphi$ is an isometric isomorphism on $B(\HH)$, which together with Theorem~\ref{Tan-310} implies that the mapping $\ker \rho \mapsto \varphi_{U^*}(\ker \rho) = \Phi_U (\ker \rho)$ is a homeomorphism on $\mathfrak{S}(B(\HH))$. We note that $\ker \rho \sim \ker \rho U$ by Theorem~\ref{sim-equiv}, which implies that $\Phi_U (\ker \rho ) \in C(\ker \rho)$ for each $\ker \rho \in \mathfrak{S}(B(\HH))$. Therefore, $\Phi_U (C(I)) = C(\Phi_U (I)) = C(I)$ for each $I \in \s (B(\HH))$. In particular, $\Phi_U$ restricts to $\s_n (B(\HH))$. From the assumption, the mappings $\Phi_0 \circ \Phi_U \circ \Phi_0^{-1}$ and $\Phi_0 \circ \Phi_{U^*} \circ \Phi_0^{-1}$ are homeomorphisms on $\s (\A)$. Since
\[
(\Phi_0 \circ \Phi_U \circ \Phi_0^{-1})(\ker \omega_{x,y} \cap \A ) = \ker \omega_{x,U^*y} \cap \A
\]
and
\[
(\Phi_0 \circ \Phi_{U^*} \circ \Phi_0^{-1})(\ker \omega_{x,y} \cap \A ) = \ker \omega_{x,Uy} \cap \A
\]
for each $x,y \in \HH \setminus \{0\}$, it follows from Theorem~\ref{continuity} that $U\A = \A$ and $U^*\A = \A$. By the self-adjointness of $\A$, we also have $\A U = \A$. From this and the fact that each element of $B(\HH)$ can be represented as a linear combination of at most four unitary elements, it turns out that $B\A = \A B = \A$ for each $B \in B(\HH)$, that is, $\A$ is a norm-closed ideal of $B(\HH)$. If $\A \neq K(\HH)$, then $K(\HH)$ is a proper norm-closed ideal of $\A$. In the light of \cite[Corollary 3.13.8]{Ped18}, each norm-closed ideal of a $C^*$-algebra is the intersection of those primitive ideals that contain it, which ensures that there exists an irreducible representation $\pi$ of $\A$ such that $K(\HH) \subset \ker \pi$. However, this is impossible since each irreducible representation of $\A$ is faithful by (i). Thus, $\A = K(\HH)$.
\end{proof}
\end{theorem}
\begin{remark}

The condition (ii) in the preceding theorem is written in the language of closure spaces, that is, it is preserved under closure space homeomorphisms. To show this, let $\A$ and $\B$ be $C^*$-algebras, and let $\Phi : \s (\A) \to \s (\B)$ be a homeomorphism. If $\A$ satisfies the condition (ii), then there exists an element $I \in \s (\A)$ with the stated property. Let $(K,c)$ be an object of $\CGss$, and if there exists a continuous finite-homeomorphism $f:C(t) \to C(\Phi (I))$ for an element $t \in K$, then $\Phi^{-1} \circ f$ is a continuous finite-homeomorphism from $C(t)$ onto $C(I)$. Therefore, $(\Phi^{-1} \circ f)\circ g \circ (\Phi^{-1} \circ f)^{-1}$ is a homeomorphism on $C(I)$ whenever $g$ is a continuous finite-homeomorphism on $C(t)$. It follows that $f \circ g \circ f^{-1}$ is a homeomorphism on $C(\Phi (I))$. This shows that $\B$ also satisfies the condition (ii).

We cannot show any implication between (i) and (ii) in the preceding theorem. Indeed, as was shown in its proof, if $\HH$ is a non-separable Hilbert space, then each proper norm-closed ideal of $B(\HH)$ other than $K(\HH)$ provides a counterexample to (ii) $\Rightarrow$ (i). Meanwhile, the implication (i) $\Rightarrow$ (ii) is equivalent to the Naimark problem, which cannot be proved in ZFC and have counterexamples, for example, in ZFC $+$ $\diamondsuit_{\aleph_1}$.
\end{remark}
Theorem~\ref{elementary} provides a characterization of elementary $C^*$-algebras in terms of their geometric structure spaces. In particular, \emph{being elementary} can be recognized in the context of closure spaces. Starting with this point, we give characterizations of type I (or CCR, or (sub)homogeneous) $C^*$-algebras. Recall that a $C^*$-algebra $\A$ is said to be \emph{$n$-homogeneous} (or \emph{n-subhomogeneous}) if every irreducible representation of $\A$ is of dimension $n$ (or of dimension less than or equal to $n$).
\begin{theorem}\label{type-I-homogeneous}

Let $\A$ be a $C^*$-algebra. Then, the following hold:
\begin{itemize}
\item[{\rm (i)}] $\A$ is of type I if and only if, for each $I \in \s (\A)$, there exist a Hilbert space $\HH_I$ and a continuous finite-homeomorphism $\Phi_I :C(I) \to \s (K(\HH_I))$.
\item[{\rm (ii)}] $\A$ is CCR if and only if, for each $I \in \s (\A)$, there exist a Hilbert space $\HH_I$ and a homeomorphism $\Phi_I :C(I) \to \s (K(\HH_I))$.
\item[{\rm (iii)}] $\A$ is $n$-homogeneous if and only if, for each $I \in \s (\A)$, there exist a Hilbert space $\HH_I$ with $\dim \HH_I =n$ and a (finite-)homeomorphism $\Phi_I :C(I) \to \s (K(\HH_I))$.
\item[{\rm (iv)}] $\A$ is $n$-subhomogeneous if and only if, for each $I \in \s (\A)$, there exist a Hilbert space $\HH_I$ with $\dim \HH_I \leq n$ and a (finite-)homeomorphism $\Phi_I :C(I) \to \s (K(\HH_I))$.
\end{itemize}
\begin{proof}
Let $\rho$ be a pure state of $\A$. Then, by Lemma~\ref{primitive} and Theorem~\ref{irr-rep}, the mapping $\Phi :\s (\pi_\rho (\A)) \to C(\ker \rho )^=$ given by $\Phi (\ker \tau ) = \ker \pi_\rho^* (\tau)$ is a homeomorphism. Meanwhile, we have $C(\ker \rho ) = \{ \ker \pi_\rho^* (\omega_{x,y}|\pi_\rho (\A)):x,y \in S_{\HH_\rho}\}$ by Lemma~\ref{pure-adjoint}. Hence, $\Phi |\s_n (\pi_\rho (\A))$ is a homeomorphism from $\s_n (\pi_\rho (\A))$ onto $C(\ker \rho )$.

Let $I \in \s (\A)$. Then, by Corollary~\ref{polar-equi}, there exists a pure state $\rho$ of $\A$ such that $C(I) = C(\ker \rho)$. Suppose that $\A$ is of type I, then $K (\HH_\rho ) \subset \pi_\rho (\A)$, which together with Theorem~\ref{continuity} (for the case $U=V=\mathbf{1}$) implies that the mapping $\Psi_\rho : \s_n(\pi_\rho (\A)) \to \s (K(\HH_\rho)$ given by $\Psi_\rho (\ker \omega_{x,y} \cap \pi_\rho (\A)) = \ker \omega_{x,y} \cap K(\HH_\rho )$ is a continuous finite-homeomorphism, which induces a continuous finite-homeomorphism $\Phi_I :C(I) \to K(\HH_\rho )$ by composing with $\Phi^{-1}|C(I)$. If $\A$ is CCR, then $\Psi_I$ and $\Phi_I$ are homeomorphisms by $K(\HH_\rho ) = \pi_\rho (\A)$. In the $n$-homogeneous (or $n$-subhomogemeous) case, we further have $\dim \HH_\rho =n$ (or $\dim \HH_\rho \leq n$).

Conversely, let $\pi :\A \to B(\HH)$ be an irreducible representation of $\A$. Then, there exists a pure state $\rho$ such that $\pi$ and $\pi_\rho$ are unitarily equivalent. Let $U_\rho :\HH \to \HH_\rho$ be an isometric isomorphism. If $\Psi :C(\ker \rho ) \to \s (K(\KK ))$ is a finite-homeomorphism, then so is $\Psi \circ (\Phi |\s_n (\pi_\rho (\A)))$. By Theorem~\ref{weak-preserver}, we obtain $\dim \HH_\rho = \dim \KK$ provided that $\KK$ is finite-dimensional, in which case, $K(\HH_\rho ) = \pi_\rho (\A) = U_\rho \pi (\A)U_\rho^*$ and $K(\HH) = \pi (\A)$. In particular, $\A$ is $n$-homogeneous (or $n$-subhomogeneous) if $\KK$ can be always chosen to satisfy $\dim \KK = n$ (or $\dim \KK \leq n$). If $\KK$ is infinite-dimensional and $\Psi$ is continuous, then Theorem~\ref{weak-preserver} generates bounded linear or conjugate-linear isomorphisms $U$ and $V$ from $\HH_\rho$ onto $\KK$ such that $U$ is linear if and only if $V$ is, and $V^*K(\KK )U \subset \pi_\rho (\A)$. It follows from $K(\HH_\rho ) = V^*K(\HH )U$ that $K(\HH_\rho ) \subset \pi_\rho (\A) = U_\rho \pi (\A) U_\rho^*$. Therefore, $K(\HH ) = U_\rho^* K(\HH_\rho )U_\rho \subset \pi (\A)$. This shows that $\A$ is of type I. If $\Psi$ can be always chosen to be homeomorphism, then $K(\HH_\rho )=V^*K(\HH )U=\pi_\rho (\A)$ and $K(\HH) = \pi (\A)$. Thus, $\A$ is CCR in this case.
\end{proof}
\end{theorem}
What is important in the preceding theorem is that each property of $C^*$-algebras is explained by closure space properties of geometric structure spaces.

\section{Nonlinear classification of $C^*$-algebras}\label{Sect:Cls}

Here, we consider nonlinear classification of $C^*$-algebras in terms of their geometric structure spaces.
\begin{definition}

Let $\mathcal{C}$ be a category. Then, the \emph{core} of $\mathcal{C}$, denoted by $\core (\mathcal{C})$, is the subcategory consisting of all objects of $\mathcal{C}$ but with morphisms only the isomorphisms of $\mathcal{C}$. 
\end{definition}
\begin{definition}

The category of Banach spaces over $\mathbb{K}$ and linear contraction, denoted as $\mathbb{K}$-\textbf{Ban}$_1$, is the category whose objects are Banach spaces over $\mathbb{K}$. The morphisms between Banach spaces $X$ and $Y$ are linear operator $T:X \to Y$ with $\|T\| \leq 1$, and the composition of morphisms is the composition of mappings. The full subcategory of $\mathbb{C}$-\textbf{Ban}$_1$ whose objects are $C^*$-algebras is denoted by $C^*$-\textbf{Ban}$_1$.
\end{definition}
\begin{remark}\label{isometry-remark}

The isomorphisms in $\mathbb{K}$-\textbf{Ban}$_1$ are isometric isomorphisms. If $\A$ and $\B$ are $C^*$-algebras, then $T:\A \to \B$ is an isometric isomorphism if and only if there exist a Jordan $*$-isomorphism $J:\A \to \B$ and a unitary element of the multiplier algebra of $\B$ such that $T(A) = UJ(A)$ for each $A \in \A$, where $J:\A \to \B$ is called a \emph{Jordan $*$-isomorphism} if $J$ is bijective, linear, $J(A^*)=J(A)^*$ and $J(AB+BA) = J(A)J(B)+J(A)J(B)$ for each $A,B \in \A$. The fact that Jordan $*$-isomorphisms between unital $C^*$-algebras are isometric isomorphisms (and preserves commutativity) was shown by Kadison~\cite[Theorem 5]{Kad51}, which was extended naturally by a result of Paterson and Sinclair~\cite[Theorem 2]{PS72}. Namely, each Jordan $*$-isomorphism between (possible non-unital) $C^*$-algebras extends uniquely to a Jordan $*$-isomorphism between their multiplier algebras. Meanwhile, the powerful fact that isometric isomorphisms between $C^*$-algebras generate Jordan $*$-isomorphisms between them was shown in \cite[Theorem 7]{Kad51} for unital $C^*$-algebras and in~\cite[Theorem 1]{PS72} for arbitrary $C^*$-algebras.
\end{remark}
From the viewpoint of nonlinear classification, it is appropriate to define $\s (\{0\})$ as $\emptyset$. Remark that there exists a unique (trivial) closure operator on $\emptyset$. Since morphisms in $\core (\mathbb{K}$-\textbf{Ban}$_1)$ are isometric isomorphisms, we can define an invariant in $\core (\mathbb{K}$-\textbf{Ban}$_1)$ as a consequence of Theorem~\ref{Tan-310}.
\begin{definition}\label{invariant-functor}

Define a $\mathbb{K}$-$\Gss$-valued invariant $\s$ in $\core (\mathbb{K}$-\textbf{Ban}$_1)$ by $X \mapsto \s (X)$ and $T \mapsto \Phi_T$, where $\Phi_T$ is a homeomorphism in Theorem~\ref{Tan-310}.
\end{definition}
A natural question concerning the preceding definition is as follows: Is the invariant $\s$ complete? In this direction, it is already known that there exists a pair of Banach spaces $(X,Y)$ over $\mathbb{K}$ such that $X$ is not (isometrically) isomorphic to $Y$ but $\s (X)$ is homeomorphic to $\s (Y)$. The complex case was settled in uniform algebras, see \cite[Example 4.14]{Tan23b}, while the real case was recently studied in \cite{Tan}. Therefore, the invariant $\s$ is not complete in the whole category $\core (\mathbb{K}$-\textbf{Ban}$_1)$. Meanwhile, it was shown in \cite[Theorem 5.2]{Tan22b} and \cite[Theorem 3.16]{Tan23a} that $\s$ is complete in the context of abelian $C^*$-algebras. At least, there is still no counterexample in $C^*$-algebras. Motivated by these circumstances, we attack the following problem in the rest of this paper, and give some partial answers to it.
\begin{problem}\label{main-p}

Is the invariant $\s$ complete in $\core (C^*$-\textbf{Ban}$_1)$?
\end{problem}
\begin{remark}

By Remark~\ref{isometry-remark}, the preceding problem is equivalently reformulated as follows: If $C^*$-algebras $\A$ and $\B$ have homeomorphic geometric structure spaces, are they Jordan $*$-isomorphic?
\end{remark}
We begin with understanding the effects of the existence of homeomorphisms between geometric structure spaces of $C^*$-algebras.
\begin{theorem}\label{homeo-effect}

Let $\A$ and $\B$ be $C^*$-algebras. Suppose that $\s (\A)$ and $\s (\B)$ are homeomorphic. Then, there exists a homeomorphism $\varphi :\hat{A} \to \hat{B}$ such that $\s (\pi (\A))$ and $\s (\pi' (\B))$ are homeomorphic whenever $\varphi (\overline{\pi}) = \overline{\pi'}$. Moreover, $\A$ is of type I or CCR or $n$-homogeneous or $n$-subhomogeneous if and only if $\B$ has the same corresponding property.
\begin{proof}
Let $\Phi :\s (\A) \to \s (\B)$ be a homeomorphism. Then, by Lemma~\ref{finite-conti}, we have $\Phi (C(I)) = C(\Phi (I))$ and $\Phi (C(I)^=) = C(\Phi (I))^=$ for each $I \in \s (\B)$. In particular, $\Phi |C(I)$ is a homeomorphism from $C(I)$ onto $C(\Phi (I))$ by Proposition~\ref{subspace-homeo}. Hence, the proof of the latter half is completed by combining this with Theorem~\ref{type-I-homogeneous}.

To prove the former half, recall that Theorem~\ref{spec-trans} generates the homeomorphism $\Phi^{\gs} :\gs (\A) \to \gs (\B)$ given by $\Phi^{\gs} (C(I)) = C(\Phi (I))$ for each $I \in \s (\A)$. Moreover, by Theorem~\ref{main-theorem}, natural homeomorphisms $\Phi_\A:\hat{\A} \to \gs (\A)$ and $\Phi_\B:\hat{\B} \to \gs (\B)$ are given by $\Phi_\A (\overline{\pi_\rho}) = C(\ker \rho)$ and $\Phi_\B (\overline{\pi_\tau}) = C(\ker \tau)$ for each pure state $\rho$ and $\tau$ of $\A$ and $\B$, respectively. Therefore, we can define a homeomorphism $\varphi :\hat{\A} \to \hat{\B}$ by $\varphi = \Phi_\B^{-1} \circ \Phi^{\gs} \circ \Phi_\A$. If $\rho$ and $\tau$ are pure states of $\A$ and $\B$, then $\varphi (\overline{\pi_\rho}) = \overline{\pi_\tau}$ if and only if $C(\Phi (\ker \rho)) = C(\ker \tau)$. Now, suppose that $\varphi (\overline{\pi_\rho}) = \overline{\pi_\tau}$. Then, $\Phi (C(\ker \rho)^=) = C(\Phi (\ker \rho))^= = C(\ker \tau)^=$, which implies that $C(\ker \rho)^=$ and $C(\ker \tau)^=$ are homeomorphic. Meanwhile, by Lemma~\ref{primitive} and Theorem~\ref{irr-rep}, $\s (\pi_\rho (\A))$ is homeomorphic to $C(\ker \rho)^=$, while $\s (\pi_\tau (\B))$ is homeomorphic to $C(\ker \tau)^=$. Therefore, $\s (\pi_\rho (\A))$ and $\s (\pi_\tau (\B))$ are homeomorphic. This completes the proof since $\pi (\A)$ and $\pi' (\B)$ are $*$-isomorphic to $\pi_\rho (\A)$ and $\pi_\tau (\B)$, respectively, whenever $\pi \in \overline{\pi_\rho}$ and $\pi' \in \overline{\pi_\tau}$.
\end{proof}
\end{theorem}
To obtain results related to the preceding theorem with higher resolution, we have to focus on each of irreducible representations of $C^*$-algebras. Below we present a basic tool for this.
\begin{lemma}\label{normal-normal}

Let $\A$ and $\B$ be $C^*$-algebras. Suppose that $\A$ acts irreducibly on a Hilbert spaces $\HH$, and that $\s (\A)$ and $\s (\B)$ are homeomorphic. Then, there exist a faithful irreducible representation $\pi$ of $\B$ and a homeomorphism $\Phi :\s (\A) \to \s (\pi (\B))$ such that $\Phi (\s_n (\A)) = \s_n (\pi (\B))$.
\begin{proof}
Let $\Phi : \s (\A) \to \s (\B)$ be a homeomorphism. By Lemma~\ref{normal-class}, we have $\s_n (\A) = C(\ker \omega_{x,y} \cap \A)$ for any choice of $x,y \in \HH \setminus \{0\}$. In particular, since $J_{\s_n (\A)} = \{0\}$, it follows from Lemma~\ref{closed} that $\s (\A) = \s_n (\A)^=$. Hence, $\Phi (\s_n (\A)) = C(\Phi (\ker \omega_{x,y} \cap \A))$ and
\[
C(\Phi (\ker \omega_{x,y} \cap \A))^= = \Phi (\s_n (\A))^= = \Phi (\s_n (\A)^=) = \s (\B) .
\]
Set $I' = \Phi (\ker \omega_{x,y} \cap \A)$ for short. Since $C(I')^= = \s (\B)$, we obtain $J_{C(I')} = \{0\}$ by \cite[Lemma 4.1]{Tan23a}. It follows from Corollary~\ref{polar-equi} that there exists a pure state $\tau$ of $\B$ such that $C(I') = C(\ker \tau)$. We note that $\ker \pi_\tau = J_{C(\ker \tau)} = \{0\}$ by Lemma~\ref{primitive}, that is, $\pi_\tau$ is faithful. Moreover, we derive that
\[
\Phi_{\pi_\tau^{-1}}(\ker \rho) = \pi_\tau^{-1}(\ker \rho) = \ker \pi_\tau^* (\rho)
\]
for each $\ker \rho \in \s (\pi_\tau (\B))$. Since $\Phi_{\pi_\tau}$ is a homeomorphism from $\s (\B)$ onto $\s (\pi_\tau (\B))$ with $\Phi_{\pi_\tau^{-1}}=\Phi_{\pi_\tau}^{-1}$ by Theorem~\ref{Tan-310} and $C(\ker \tau ) = \Phi_{\pi_\tau^{-1}}(\s_n (\pi_\tau (\B)))$ by Lemma~\ref{pure-adjoint}, it follows that the mapping $\Phi_{\pi_\tau} \circ \Phi$ is a homeomorphism from $\s (\A)$ onto $\s (\pi_\tau (\B))$ with $(\Phi_{\pi_\tau} \circ \Phi)(\s_n (\A)) = \s_n (\pi_\tau (\B))$.
\end{proof}
\end{lemma}
Let $\A$ be a $C^*$-algebra, and let $\pi :\A \to B(\HH)$ be an irreducible representation. Recall that $\pi$ is called a GCR representation if $K(\HH) \subset \pi (\A)$. A primitive ideal $\mathcal{I}$ of $\A$ is said to be \emph{$\A$-GCR} if there exists a GCR representation $\pi$ of $\A$ such that $\mathcal{I}=\ker \pi$.
\begin{theorem}\label{CCR-GCR-preserved}

Let $\A$ and $\B$ be are irreducible $C^*$-algebras acting on Hilbert spaces $\HH$ and $\KK$, respectively. Suppose that $\s (\A)$ and $\s (\B)$ are homeomorphic. Then, $K(\HH) = \A$ or $K(\HH)\subsetneq \A$ or $K(\HH) \cap \A = \{0\}$ if and only if $\B$ has the same corresponding property. Moreover, if $K(\HH ) \subset \A$ or $K(\KK) \subset \B$, then $\dim \HH = \dim \KK$.
\begin{proof}
Suppose that $K(\KK) \subset \B$. Then, the identity representation of $\B$ is (primitive and) GCR, and $\{0\}$ is a $\B$-GCR ideal in $\B$. Hence, by \cite[Corollary IV.1.3.6]{Bla06}, any faithful irreducible representation or $\B$ is unitarily equivalent to the identity representation of $\B$. Meanwhile, by Lemma~\ref{normal-normal}, there exist a faithful irreducible representation $\pi :\B \to B(\KK_0)$ and a homeomorphism $\Phi :\s (\A) \to \s (\pi (\B))$ such that $\Phi (\s_n (\A)) = \s_n (\pi (\B))$. Let $W:\KK \to \KK_0$ be an isometric isomorphism such that $\pi (B) = WBW^*$ for each $B \in \B$. We note that $\dim \KK = \dim \KK_0$.

If $\KK$ is finite-dimensional, then $\dim \HH = \dim \KK_0 = \dim \KK$ by Theorem~\ref{weak-preserver}, in which case, $K(\HH)=\A$ and $K(\KK)=\B$. Suppose that $\KK$ is infinite-dimensional. Then, again by Theorem~\ref{weak-preserver}, $\HH$ is necessarily infinite-dimensional and there exist bounded linear or conjugate-linear isomorphisms $U$ and $V$ from $\HH$ onto $\KK_0$ such that $U$ is linear if and only if $V$ is, and $V^*\pi (\B)U = \A$. Hence, it follows from $K(\KK) \subset \B$ that
\[
K(\HH) = V^*WK(\KK)W^*U \subset V^*W\B W^*U = \A .
\]
In particular, if $K(\KK) = \B$, then $K(\HH) = \A$. This shows that $K(\KK) \subset \B$ (or $K(\KK) = \B$) implies that $\dim \HH = \dim \KK$ and $K(\HH) \subset \A$ (or $K(\HH) =\A$). The converse to these implications also follow from an argument similar to above. Combining this with~\cite[IV.1.2.5]{Bla06}, we can conclude that $\A \cap K(\HH) =\{0\}$ if and only if $\B \cap K(\KK) = \{0\}$.
\end{proof}
\end{theorem}
As an application of Theorems~\ref{homeo-effect} and \ref{CCR-GCR-preserved}, we give the first example of a class of non-abelian $C^*$-algebras that are completely determined by their geometric structure spaces. This answers Problem~\ref{main-p} partially and affirmatively.

Let $\A$ be a $C^*$-algebra. Then, $\A$ is said to be \emph{compact} if it is $*$-isomorphic to a $C^*$-subalgebra of $K(\HH)$ for some Hilbert space $\HH$. For compact $C^*$-algebras, we know a reasonable characterization by direct sums of elementary $C^*$-algebras; see~\cite[Theorem 1.4.5]{Arv76}. Namely, a $C^*$-algebra $\A$ is compact if and only if there exists a family of Hilbert spaces $(\HH_\lambda )_{\lambda \in \Lambda}$ such that $\A$ is $*$-isomorphic to $\sum_{\lambda \in \Lambda}K(\HH_\lambda )$, where $\sum_{\lambda \in \Lambda}K(\HH_\lambda )$ denotes the $c_0$-direct sum of the family $(K(\HH_\lambda ))_{\lambda \in \Lambda}$. The following lemma translates this into the language of geometric structure spaces. It should be mentioned that the lemma was due to Huruya~\cite{Hur69} if $\hat{\A}$ is a finite set.
\begin{lemma}\label{CCR-discrete}

Let $\A$ be a $C^*$-algebra. Then, $\A$ is compact if and only if $\A$ is CCR and $\Prim (\A)$ is discrete, in which case, $\pi (\A) = \sum_{\lambda \in \Lambda}\pi_\lambda (\A)$, where $\pi = \sum_{\lambda \in \Lambda}\pi_\lambda$ is the reduced atomic representation of $\A$.
\begin{proof}
Suppose that $\A$ is compact. Then, it may be assumed that $\A = \sum_{\lambda \in \Lambda}K(\HH_\lambda )$ for a family of Hilbert space $(\HH_\lambda )_{\lambda \in \Lambda}$. Let $\pi_\mu$ be a representation of $\A$ given by $\pi_\mu ((A_\lambda )_{\lambda \in \Lambda} ) = A_\mu$ for each $(A_\lambda )_{\lambda \in \Lambda} \in \A$. Since $\pi_\mu (\A) = K(\HH_\mu )$, it follows that $\pi_\mu$ is an irreducible representation of $\A$. Moreover, each irreducible representation $\pi$ of $\A$ is unitarily equivalent to some $\pi_\mu$. Indeed, since $\A = \sum_{\lambda \in \Lambda}K(\HH_\lambda )$, it follows that $\A^*$ is the $\ell_1$-direct sum of the family $(K(\HH_\lambda )^*)_{\lambda \in \Lambda}$. From this, if $\rho$ is a pure state of $\A$ such that $\pi_\rho$ is unitarily equivalent to $\pi$, then there exists a unique index $\mu$ such that $\rho |\A P_\mu \neq 0$, where $P_\mu$ is the projection from the $\ell_2$-direct sum $\HH$ of the family $(\HH_\lambda)_{\lambda \in \Lambda}$ onto the subspace naturally identified with $\HH_\mu$. In particular,  $\rho |\A P_\mu$ is a pure state of $\A P_\mu$ which is the natural embedding of $K(\HH_\mu)$. Hence, we have an element $x \in S_{\HH_\mu}$ such that
\[
\rho ((A_\lambda )_{\lambda \in \Lambda} ) = \langle A_\mu x,x \rangle = \langle \pi_\mu ((A_\lambda )_{\lambda \in \Lambda} )x,x \rangle
\]
for each $(A_\lambda )_{\lambda \in \Lambda} \in \A$. Meanwhile, since $x$ is cyclic in $\pi_\mu (\A)$, the representations $\pi_\rho$ and $\pi_\mu$ are unitarily equivalent. Therefore, $\pi \in \overline{\pi_\mu}$. This shows that $\hat{\A} = \{ \overline{\pi_\lambda} : \lambda \in \Lambda \}$ and $\Prim (\A) = \{ \ker \pi_\lambda : \lambda \in \Lambda \}$. In particular, it follows that $\A$ is CCR. The discreteness of $\Prim (\A)$ follows from $\bigcap_{\lambda \neq \mu}\ker \pi_\lambda = \A P_\mu \not \subset \ker \pi_\mu$ for each $\mu$.

Conversely, suppose that $\A$ is CCR and $\Prim (\A)$ is discrete. Set $\Prim (\A) = \{ \mathcal{I}_\lambda : \lambda \in \Lambda \}$, where $\mathcal{I}_\lambda \neq \mathcal{I}_\mu$ whenever $\lambda \neq \mu$. Then, for each $\lambda \in \Lambda$, there exists an irreducible representation $\pi_\lambda$ of $\A$ such that $\mathcal{I}_\lambda =\ker \pi_\lambda$. Since $\A$ is CCR, by the Glimm-Sakai theorem, the mapping $\pi \mapsto \ker \pi$ is a bijection from $\hat{\A}$ onto $\Prim (\A)$. Therefore, $\hat{\A}=\{ \overline{\pi_\lambda} : \lambda \in \Lambda \}$ and $\Prim (\A) = \{ \ker \pi_\rho : \lambda \in \Lambda \}$. Moreover, $\overline{\pi_\lambda} \neq \overline{\pi_\mu}$ and $\ker \pi_\lambda \neq \ker \pi_\mu$ whenever $\lambda \neq \mu$.

By the discreteness of $\Prim (\A)$, we have $\bigcap_{\lambda \neq \mu}\ker \pi_\lambda \not \subset \ker \pi_\mu$ for each $\mu$. Let $A_\mu \in \bigcap_{\lambda \neq \mu} \ker \pi_\lambda \setminus \ker \pi_\mu$, and let $\mathcal{J}_\mu = \bigcap_{\lambda \neq \mu} \ker \pi_\lambda$ for each $\mu$. Then, $\mathcal{J}_\mu$ is a norm-closed ideal of $\A$ and $0 \neq \pi_\mu (A_\mu ) \in \pi_\mu (\mathcal{J}_\mu )$. Remark that $\pi_\mu$ is faithful on $\mathcal{J}_\mu$ since $\mathcal{J}_\mu \cap \ker \pi_\mu = \{0\}$, which ensures that $\pi_\mu (\mathcal{J}_\mu)$ is a $C^*$-subalgebra of $\pi_\mu (\A) = K(\HH_\mu )$. Meanwhile since $\pi_\mu (A_\mu ) \neq 0$, there exists an $x_0 \in S_{\HH_\mu}$ such that $y_0=\pi_\mu (A_\mu )x_0 \neq 0$. Let $x_1,\ldots ,x_n ,y_1,\ldots ,y_n \in \HH_\mu$ be such that $\{ x_1,\ldots ,x_n \}$ is linearly independent. Then, for each $j$, there exist $A_j,B_j \in \A$ such that $\pi_\mu (A_j)y_0 = y_j$ and $\pi_\mu (B_j)x_i = \delta_{ij}x_0$. It follows that $A_jA_\mu B_j \in \mathcal{J}_\mu$ and $\pi_\mu (A_jA_\mu B_j)x_i = \delta_{ij}y_j$, which implies that $\sum_{j=1}^n \pi_\mu (A_jA_\mu B_j)x_i = y_i$ for each $i$. This shows that $\pi_\mu (\mathcal{J}_\mu )$ also acts irreducibly on $\HH_\mu$. Combining this with $\pi_\mu (\mathcal{J}_\mu ) \subset K(\HH_\mu)$, we obtain $\pi_\mu (\mathcal{J}_\mu ) = K(\HH_\mu )$ by \cite[Proposition IV.1.2.4]{Bla06}.

Let $\pi = \sum_{\lambda \in \Lambda}\pi_\lambda$. Then, $\pi$ is the reduced atomic representation of $\A$ which is known to be faithful. If $P_\mu$ is the projection from the $\ell_2$-direct sum $\sum_{\lambda \in \Lambda}\HH_\lambda$ of the family $(\HH_\lambda )_{\lambda \in \Lambda}$ onto the subspace naturally identified with $\HH_\mu$, and if $(B_\lambda )_{\lambda \in \Lambda} \in \sum_{\lambda \in \Lambda}K(\HH_\lambda)$, then for each $\varepsilon >0$, there exist $\lambda_1,\ldots ,\lambda_n \in \Lambda$ such that $\| \sum_{j=1}^n (B_\lambda )_{\lambda \in \Lambda}P_{\lambda_j} -(B_\lambda )_{\lambda \in \Lambda} \| < \varepsilon$. Moreover, by the preceding paragraph, $B_{\lambda_j} = \pi_{\lambda_j}(A_{\lambda_j})$ for some $A_{\lambda_j} \in \mathcal{J}_{\lambda_j}$, which implies that $\pi (A_{\lambda_j}) = (\pi_{\lambda}(A_{\lambda_j}))_{\lambda \in \Lambda} = (B_\lambda )_{\lambda \in \Lambda }P_{\lambda_j}$ for each $j$. Hence, $\sum_{j=1}^n (B_\lambda )_{\lambda \in \Lambda}P_{\lambda_j} = \pi (\sum_{j=1}^n A_{\lambda_j}) \in \pi (\A)$ and $(B_\lambda )_{\lambda \in \Lambda} \in \pi (\A)$. It follows that $\pi (\A) \supset \sum_{\lambda \in \Lambda}K(\HH_\lambda)$.

We claim that $\pi (\A) = \sum_{\lambda \in \Lambda}K(\HH_\lambda)$. If $\pi (\A) \supsetneq \sum_{\lambda \in \Lambda}K(\HH_\lambda)$, then $\sum_{\lambda \in \Lambda}K(\HH_\lambda)$ is a proper norm-closed ideal of $\pi (\A)$. Since each closed ideal of a $C^*$-algebra is the intersection of those primitive ideals that contain it (\cite[Corollary 3.13.8]{Ped18}), there exists an irreducible representation $\varphi$ of $\pi (\A)$ such that $\sum_{\lambda \in \Lambda}K(\HH_\lambda) \subset \ker \varphi$. We note that $\varphi \circ \pi$ is an irreducible representation of $\A$, which implies that $\varphi \circ \pi \in \overline{\pi_\mu}$ for some $\mu \in \Lambda$. In particular, $\ker (\varphi \circ \pi ) = \ker \pi_\mu$. Meanwhile since $\pi_\mu (\mathcal{J}_\mu ) = K(\HH_\mu )$, there exists an $A \in \mathcal{J}_\mu$ such that $\pi_\mu (A) \neq 0$. However, then $(\varphi \circ \pi)(A) = 0$ by $\pi (A) \in \sum_{\lambda \in \Lambda}K(\HH_\lambda )$. This contradicts $\ker (\varphi \circ \pi ) = \ker \pi_\mu$. Thus, $\pi (\A) = \sum_{\lambda \in \Lambda}K(\HH_\lambda)$.
\end{proof}
\end{lemma}
\begin{remark}

The fact that compact $C^*$-algebras are always CCR was already noted in \cite[the paragraph following Definition 1.5.1]{Arv76}.
\end{remark}
Now, we prove that compact $C^*$-algebras are completely determined by their geometric structure spaces.
\begin{theorem}\label{cpt-op-alg}

Let $\A$ and $\B$ be $C^*$-algebras. Suppose that either $\A$ or $\B$ is compact. Then, $\A$ and $\B$ are $*$-isomorphic if and only if $\s (\A)$ and $\s (\B)$ are homeomorphic.
\begin{proof}
It is sufficient to prove the ``if'' part. Suppose that $\s (\A)$ and $\s (\B)$ are homeomorphic. If either $\A$ or $\B$ is compact, then both of them are compact by Corollary~\ref{spec-preserved}, Theorem~\ref{homeo-effect} and Lemma~\ref{CCR-discrete}. Moreover, we have $\pi (\A) = \sum_{\lambda \in \Lambda} \pi_\lambda (\A)$ for the reduced atomic representation of $\A$. We note that $\hat{\A} = \{ \overline{\pi_\lambda} : \lambda \in \Lambda\}$ by the definition of $\pi$. Meanwhile, Theorem~\ref{homeo-effect} generates a homeomorphism $\varphi :\hat{\A} \to \hat{\B}$ such that $\s (\pi_\lambda (\A)) = \s (\pi'_\lambda (\B))$ whenever $\varphi (\overline{\pi_\lambda}) = \overline{\pi'_\lambda}$. Since $\A$ and $\B$ are CCR, it follows from Theorem~\ref{CCR-GCR-preserved} that irreducible $C^*$-algebras $\pi_\lambda (\A)$ and $\pi'_\lambda (\B)$ are $*$-isomorphic. Now, again by Lemma~\ref{CCR-discrete}, the reduced atomic representation $\pi' = \sum_{\lambda \in \Lambda}\pi'_\lambda$ of $\B$ satisfies $\pi' (\B) = \sum_{\lambda \in \Lambda}\pi'_\lambda (\B)$. Thus, $\A$ and $\B$ are $*$-isomorphic.
\end{proof}
\end{theorem}
\begin{corollary}\label{cpt-invariant}

Let $\A$ and $\B$ be $C^*$-algebras. Suppose that either $\A$ or $\B$ is elementary. Then, $\A$ and $\B$ are $*$-isomorphic if and only if $\s (\A)$ and $\s (\B)$ are homeomorphic.
\end{corollary}
\begin{corollary}\label{finite-dim-invariant}

Let $\A$ and $\B$ be $C^*$-algebras. Suppose that either $\A$ or $\B$ is finite-dimensional. Then, $\A$ and $\B$ are $*$-isomorphic if and only if $\s (\A)$ and $\s (\B)$ are homeomorphic.
\end{corollary}
Recall that the algebras of the form $B(\HH)$ are recognized as \emph{type I factors} in the theory of operator algebras, although $B(\HH)$ is not a type I $C^*$-algebra if $\HH$ is infinite-dimensional: see, for example, \cite[IV.1.1.5]{Bla06}. Hence, the following theorem gives an example of a non-compact $C^*$-algebra that is completely determined by its geometric structure space.
\begin{theorem}\label{full-algebra}

Let $\A$ and $\B$ be $C^*$-algebras. Suppose that either $\A$ or $\B$ is a type I factor. Then, $\A$ and $\B$ are $*$-isomorphic if and only if $\s (\A)$ and $\s (\B)$ are homeomorphic.
\begin{proof}
It is sufficient to prove the ``if'' part. Suppose that $\s (\A)$ and $\s (\B)$ are homeomorphic. It may be assumed that $\A = B(\HH)$ for some Hilbert space $\HH$. Since $\A$ acts irreducibly on $\HH$, by Lemma~\ref{normal-normal}, there exist a faithful irreducible representation $\pi :\B\to B(\KK)$ and a homeomorphism $\Phi :\s (B(\HH)) \to \s (\pi (\B))$ such that $\Phi (\s_n (B(\HH))) = \s_n (\pi (\B))$. By Theorem~\ref{weak-preserver}, we have $\dim \HH = \dim \KK$,  which implies that $\pi (\B) = B(\KK)$ if either $\HH$ or $\KK$ is finite-dimensional. If $\HH$ and $\KK$ are infinite-dimensional, there exists bounded linear or conjugate-linear isomorphisms $U$ and $V$ from $\HH$ onto $\KK$ such that $U$ is linear if and only if $V$ is, and $V\*\pi (\B)U = B(\HH)$. This shows that $\pi (\B) = (V^{-1})^* B(\HH)U^{-1} = B(\KK)$. Thus, $\A$ and $\B$ are $*$-isomorphic via $B(\HH)$ and $\pi (\B)$.
\end{proof}
\end{theorem}
Finally, we introduce two classification problems related to Problem~\ref{main-p}. The first one is based on the structure determined by Birkhoff-James orthogonality. Recall that a (possibly nonlinear) mapping $T$ between Banach spaces is called a Birkhoff-James orthogonality preserver if is bijective and $x \perp_{BJ}y$ is equivalent to $Tx \perp_{BJ}Ty$. A nonlinear equivalence of Banach spaces is naturally defined in terms of these mappings.
\begin{definition}

Let $X$ and $Y$ be Banach spaces. Then, $X$ and $Y$ are \emph{isomorphic with respect to the structure determined by Birkhoff-James orthogonality}, denoted by $X \sim_{BJ}Y$, if there exists a Birkhoff-James orthogonality preserver $T:X \to Y$.
\end{definition}
This notion was first introduced in \cite[Definition 1.3]{Tan22a}. It is known that if $X$ and $Y$ are smooth Banach spaces, then $X \sim_{BJ}Y$ if and only if $X$ is isometrically isomorphic to $Y$ (or $\overline{Y}$ in the complex case); see~\cite[Theorem 3.11]{AGKRZ23} and also~\cite{IT22,Tan22c} for interim results. Moreover, it was shown in \cite[Theorems 5.5 and 5.14]{Tan22b} that spaces of continuous functions and classical sequence spaces are isometrically classified by their structure determined by Birkhoff-James orthogonality. We remark that these results on spaces of continuous functions and classical sequence spaces are also valid in terms of geometric structure spaces; see \cite[Theorems 3.16 and 6.14]{Tan23a}.

Meanwhile, there is a graph theoretical approach to studying the structure determined by Birkhoff-James orthogonality. The following definition was given in \cite{AGKRZ23}.
\begin{definition}

Let $X$ be a Banach space. Then, the \emph{di-orthograph} $\OG (X)$ of $X$ consists of the vertices $V(\OG (X)) = \{ \overline{x} : x \in X \setminus \{0\} \}$ and arrows $E(\OG (X)) = \{ (\overline{x},\overline{y}) : x \perp_{BJ}y \}$, where $\overline{x}$ is the one-dimensional subspace of $X$ generated by $x$, and $(\overline{x},\overline{y}) \in E(\OG (X))$ means that $\OG (X)$ contains the directed arrow $\overline{x} \to \overline{y}$.
\end{definition}
An advantage of di-orthographs is that smoothness of finite-dimensional Banach spaces as well as rotundity of Banach spaces are characterized in terms of them; see~\cite[Lemmas 2.5 and 2.6]{AGKRZ23}. It is worth mentioning that a non-directed version was also considered in \cite{AGKRZ21} earlier.

Using graph isomorphisms between di-orthographs, we can introduce another geometric nonlinear equivalence of Banach spaces. Recall that a mapping $f$ from a directed graph $G$ to another $H$ is called a \emph{graph isomorphism} if it is a bijection from the vertices of $G$ onto that of $H$ and $x \to y$ is equivalent to $f(x) \to f(y)$.
\begin{definition}

Let $X$ and $Y$ be Banach spaces. Then, $X$ and $Y$ said to \emph{isomorphic with respect to di-orthographs}, denoted by $X\sim_{\OG}Y$ if there exists a graph isomorphism $f:\OG (X) \to \OG (Y)$.
\end{definition}
\begin{remark}

Let $X$ be a Banach space. Then, $X$ can be itself viewed as a directed graph by setting $V(X) = X$ and $E(X) = \{ (x,y) : x\perp_{BJ}y\}$. The graph isomorphisms with respect to this graph are just Birkhoff-James orthogonality preservers.
\end{remark}

We note that an isometric isomorphism $T$ from a Banach space $X$ onto another $Y$ induces graph isomorphisms $f_T : \OG (X) \to \OG (Y)$ and $T:X \to Y$, where $f_T$ is given by $f_T (\overline{x}) = \overline{Tx}$ for each $x \in X \setminus \{0\}$. Hence, we can define two invariants in $\core(\mathbb{K}$-\textbf{Ban}$_1)$. Let $\textbf{Quiv}$ denote the category of directed graphs.
\begin{definition}

Define \textbf{Quiv}-valued invariants $\hat{F}$ and $F_{BJ}$ in $\core(\mathbb{K}$-\textbf{Ban}$_1)$ by $\hat{F}(X) = \OG (X)$, $\hat{F}(T)=f_T$, $F_{BJ}(X) = X$, and $F_{BJ}(T)=T$, respectively.
\end{definition}

Now, we formulate two problems related to Problem~\ref{main-p}.
\begin{problem}\label{not-main-p}

Are the invariants $\hat{F}$ and $F_{BJ}$ complete in $\core(C^*$-\textbf{Ban}$_1)$?
\end{problem}
\begin{remark}

Let $X$ be a Banach space, and let $x,y \in X$. Then, $x \perp_{BJ}y$ in $X$ if and only if $x \perp_{BJ}y$ in $\overline{X}$. Moreover, we have $V(\OG (X)) = V(\OG(\overline{X}))$. Hence, $\OG (X) = \OG (\overline{X})$, that is, $X \sim_{\OG}\overline{X}$. From this, we can at most expect from $X \sim_{\OG}Y$ or $X \sim_{BJ}Y$ that $X$ is isometrically isomorphic to \emph{either} $Y$ or $\overline{Y}$ for general complex Banach spaces $X$ and $Y$. Meanwhile, for a $C^*$-algebra $\A$, we always have the isometric isomorphic $A \mapsto A^*$ from $\A$ onto $\overline{\A}$. This justify the formulations of Problems~\ref{main-p} and \ref{not-main-p}.
\end{remark}

\end{document}